\documentclass[11pt,english,12p]{amsart}
\usepackage{amsmath}
\usepackage{amssymb}
\usepackage{tikz}
\usepackage[autostyle]{csquotes}
\usepackage{mathrsfs}
\usepackage{float}
\usepackage{tikz-cd}

\usepackage{pst-node}
\usetikzlibrary{patterns}
\usetikzlibrary{arrows,positioning,shapes,fit,calc}

 \usepackage[all,cmtip]{xy}

\usetikzlibrary{arrows,positioning,shapes,fit,calc}

\newtheorem{thm}{Theorem}[section]
\newtheorem{lem}[thm]{Lemma}

\newtheorem{defn}[thm]{Definition}
\newtheorem{exmp}[thm]{Example}

\newtheorem{remark}[thm]{Remark}

\newcommand{\R}{{\mathbb R}}

\newcommand{\N}{{\mathbb N}}
\newcommand{\C}{{\mathbb C}}
\newcommand{\Z}{{\mathbb Z}}

\def\g{\gamma}

\def\k{\kappa}

\selectfont

\title{  Expository paper on  Clifford algebras, Representations,  and the  octonion algebra.}

\date{\today}

\author[R. Su\'arez]{Ricardo Su\'arez}
\address{Departament of Mathematics, California State University Channel Islands, One University Drive, 93012 Camarillo, CA, United States}
\email{ricardo.suarez532@myci.csuci.edu}

\begin{document}
\maketitle

\normalsize

\begin{abstract}
This paper is meant to be an informative introduction to spinor representations of Clifford algebras.
In this paper we will have a look at Clifford algebras  and the octonion algebra. We begin the paper looking at the quaternion algebra $\mathbb{H}$ and  basic properties that relate Clifford algebras and the well know Pin and Spin groups. We then will look at generalized spinor representations of Clifford algebras, along with many examples. We conclude the paper looking at the octonion algebra $\mathbb{O}$. This paper provides background to constructing representations which can be used to look at elements in the appropriate Pin and Spin groups. 
\end{abstract}
\maketitle

\section{Introduction}

We begin this paper focusing on the well known division algebras $\mathbb{H}$ and $\C$ , and the notion of quadratic spaces whom we will use in the construction of Clifford algebras with signature $(p,q)$, denoted $Cl_{p,q}$. Clifford algebras will be key in establishing relations with the  quadratic spaces and the well known orthogonal and special orthogonal groups; $O(V)$ and $SO(V)$. We then go on to examine spinors and their relation in finding representations for the Clifford algebras over $\R$ , $\C$ and $\mathbb{H}$ . In this paper we talk about generalized spinor representations for Clifford algebras through the method of constructing primitive idempotents ,which we denote $F$ , that we use to build our spinor spaces $S_{p,q}$, this method is explored in detail in [LW]. The Clifford algebra $Cl_{0,7}$ is of particular interest because the octonion algebra $\mathbb{O}$ is contained as the space of para vectors. In this paper we provide some background to Clifford algebras and spin groups, and how to construct spinor representations . We then go on to provide generalized spinor basis over $\R$,$\C$ and $\mathbb{H}$ for all signatures $(p,q)$. We conclude this paper with some comments about the octonion algebra $\mathbb{O}$ and  its automorphism group ,the exceptional lie group $G_2$.

\section{background on Quadratic spaces and orthogonal groups }
\subsection{The Quaternions}
The quaternion algebra $\mathbb{H}$ can be viewed as $\R^4$  with the basis elements  $1,i,j,k$ which all correspond to a real dimension orthogonal to the other 3. Multiplication in $\mathbb{H}$ is known as the quaterntion product $\R^4\times \R^4\rightarrow \R^4$,  with  the relations $ij=k,jk=i,ki=j,kj=-i,ji=-k,ik=-j$ $ijk=-1$  where squares are negative definite , that is $i^2=j^2=k^2=-1$. Thus any element in $\mathbb{H}$ can be written as $q=a+bi+cj+dk$ such that $a,b,c,d\in\R$, where $re(\mathbb{H})=\R$ and $im(\mathbb{H})=\R i\oplus \R j\oplus \R k$. The  \textbf{conjugate} in $q$ denoted $\bar{q}$ is defined as $\bar{q}=a-bi-cj-dk$, such that $q\bar{q}=a^2+b^2+c^2+d^2$.This conjugate induces a bi-linear form and a norm on $\mathbb{H}$. The bi-linear form defined as  $<,>:\mathbb{H}\times \mathbb{H}\rightarrow \R$,such that $<q,r>=re(\bar{q}r)$, and the induced norm being  $||q||=\sqrt{<q,q>}$. $\mathbb{H}$ is a division algebra since for any quaternion $q\in\mathbb{H}^*=\mathbb{H}\setminus{\{0\}}$ , there exists an inverse $q^{-1}$ written as $q^{-1}=\dfrac{\bar{q}}{||q||^2}$. Using this norm we can get the set of units in $\mathbb{H}$ and denote them $\mathbb{H}_{\Delta}=\{x\in\mathbb{H}:||x||=1\}$. For  a usual element $x=a+bi+cj+dk\in\mathbb{H}_{\Delta}$ we have  $||x||=a^2+b^2+c^2+d^2=1$. If we view the unit quaternions as  four-tuples of the form  $x=(a,b,c,d)$, we see that we have a  have the group isomorphism $S^3\cong \mathbb{H}_{\Delta}$. If we restrict our attention to the subset of the unit quaternions that have a negative square, that is  $x^2=-1$, then $x\in im(\mathbb{H})$. Now the set of all such unit quaternions give us , $\{x\in\mathbb{H}:x^2=-1,||x||=1\}\cong S^2$. [Bu].For any $u\in im(\mathbb{H})$ such that $|u|=1$, thus we can write quaternions using the Euler identity , that is $q=re^{u\theta}$. Note that $r=|q|$ and $u$ can be thought of as the imaginary direction of $q$[DM].This fact often is expressed  $\mathbb{H}=\C\oplus \C j$, but this can be done with any imaginary quaternion of norm 1 [DM].

\subsection{Involutions on $\C$ and $\mathbb{H}$}
It is a well known fact that the only finite dimensional $\R$ division algebras are $\R$, $\C$ and $\mathbb{H}$ [Ga]. Since the concept of an involution and anti-involution will be important in Clifford algebras we begin here describing some involutions and anti-involutions for $\C$ and $\mathbb{H}$.
 \begin{defn}
 An \textbf{Involution} on an $\R$-space $A$ is an automorphism $t:A\rightarrow A$ such that $t\circ t=id_A$. An \textbf{anti-involution} is an involution that reverses multiplication in the given algebra $A$, that is $t(ab)=t(b)t(a)$.
 
 \end{defn}
 
On $\C$ we know that conjugation is an involution since $\bar{\bar{z}}=z$ for any $z\in\C$. For $\mathbb{H}$ the conjugation map is an anti-involution where $\bar{q}=re(q)-im(q)$, such that $\bar{qr}=\bar{r}\bar{q}$. In $\mathbb{H}$ we have the main involution $\hat{\ }:\mathbb{H}\rightarrow\mathbb{H}$ , such that $\hat{q}=jqj^{-1}=-jqj$, and the reversion anti-involution $\sim{\ }:\mathbb{H}\rightarrow \mathbb{H}$ is the composition of the conjugation anti-involution and the main involution, that is $\tilde{q}=\hat{\bar{q}}=\bar{\hat{q}}$ [Po]. The reversion map  $\sim:\mathbb{H}\rightarrow \mathbb{H}$ , fixes generators $1$ , $i$ ,$k$, and for $j$ we  have $\tilde{j}= -j$.

\begin{exmp}
If we write a quaternion as $q=z_1+z_2j$ and we conjugate it by $i$ , $iq\bar{i}=i(z_1+z_2j)(-i)=(z_1-z_2j)(-i^2)=z_1-z_2j$. This means that conjugation by $i$ yields a rotation in the $jk$-plane. We would have similar results with conjugation by any other unit quaternion $u\in \mathbb{H}_{\Delta}$. Conjugating by $e^{i\theta}$ gives us $e^{i\theta}z_1e^{-i\theta}+ e^{i\theta}z_2je^{-i\theta}=z_1+e^{i2\theta}j$. That is conjugation by  $e^{i\theta}$ yields a $2\theta$ rotation on the $jk$-plane[DM].

\end{exmp}

\subsection{ $\R$- Quadratic Spaces and orthogonal groups}
The pairing $(V,q)$ will denote an $\R$-linear space $V$ and quadratic form $q$ naturally derived from bi-linear form $B:V\times V\rightarrow \R$, where $q(x)=B(x,x)$. In these quadratic spaces elements $x,y\in V$ are considered \textbf{orthogonal} if $B(x,y)=0$, with the equivalent condition that $q(x+y)=q(x)+q(y)$.A space $V$ will be \textbf{non degenerate} if it satisfies the condition that  $B(x,y)=0$ for all $y\in V\setminus\{0\}$ implies that  $x=0$. We will denote the group of automorphisms of $(V,q)$ by $Aut(V)$, the subset of orthogonal automorphism who form a  group under composition  will be denoted by $O(V)$. An orthogonal map  on a quadratic space $(V,q)$ may be thought of as an automorphism $\phi$, such that $q(x)=q(\phi(x))$ for all $x\in V$, we can equivalently say that $B(x,y)=B(\phi(x),\phi(y))$ for all $x,y\in V$. If an $\R$-vector space $V$ is of dimension $n$ with the quadratic form $q(x)=-\sum_{i=1}^p x_i^2+\sum_{i=p+1}^{p+q} x_i^2$, where $x_i$ are the components and $p+q=n$, then the quadratic vectors space $V$ is an $\R$-orthogonal space denoted $V=\R^{p,q}$ with the orthogonal group $O(p,q)$. The quadratic form is called  be positive definite if $q(x)>0$ for all non zero $x\in \R^{n}$ and negative definite if $q(x)<0$ for all $x\in\R^{n}$ . For the positive  definite case on $\R^n$ we will denote the orthogonal group $O(n)$.
\begin{remark}
Since $O(p,q)\cong O(q,p)$, some authors call the negative definite case $O(n)$, since we care about the positive definite  orthogonal space $\R^{0,7}$, we have chosen to use $O(n)$ for $\R^{0,n}$.
\end{remark}

The group $O(p,q)$ has a subgroup of rotations , that is orientation preserving orthogonal transformations, this group is denoted $SO(p,q)$and is called the \textbf{special orthgonal group}. Note that $SO(\R^{0,n})=SO(n)$. An orthogonal automorphism that reverses orientation will be called an anti-rotation. If $V$ is an orthogonal vector space we have the group isomorphism $O(V)/SO(V)\cong\Z_2\cong S^0=\{\pm 1\}$[Po]. That is for $(V,q)$ we have the short exact sequence 

$$1\rightarrow SO(V)\rightarrow O(V)\rightarrow\Z_2 \rightarrow 1$$.

Where the map from $SO(V)$ to $O(V)$ will be the inclusion map ,and the map from $O(V)$ to $\Z_2$ is the determinant map. We see that $SO(V)$ can be expressed as $SO(V)=\{\phi\in O(V):det\phi=1\}$, then the kernel of the determinant is $SO(V)$ [Ga].We will go on to give a few well known isomorphisms for low dimension special orthogonal groups.

\begin{exmp}
Some well known results for the special orthogonal group in low dimensions are :
\begin{itemize}
\item $SO(1)\cong \Z_2$.
\item $SO(2)\cong \{z\in\C:|z|=1\}\subset \C$.
\item The subset $\mathbb{H}_{\Delta}=\{q\in\mathbb{H}:|q|=1\}$   is a double cover for  $SO(3)$, with the short exact sequence 

$$1\rightarrow \Z_2\rightarrow \mathbb{H}_{\Delta}\xrightarrow\rho SO(3)\rightarrow 1$$

Where $\rho(q):\R^3\rightarrow \R^3$ where $\rho(q)(x)=qxq^{-1}$. When we view $Im(\mathbb{H})$ as $\R^3$ , $\rho(q)(x)\in Im(\mathbb{H})$ for a given $x\in Im(\mathbb{H})$.
\item  $\mathbb{H}_{\Delta}\times\mathbb{H}_{\Delta}$  is a double cover for  $SO(4)$, with the short exact sequence 

$$1\rightarrow \Z_2\rightarrow \mathbb{H}_{\Delta}\times \mathbb{H}_{\Delta}\xrightarrow\phi SO(4)\rightarrow 1$$. Where $\phi(q,r):\R^4\rightarrow \R^4$, where $\phi(q,r)(x)=qxr^{-1}$, for any $x\in\R^4$.
\end{itemize}
\end{exmp}

Any $t\in Aut(\mathbb{H})$ may be expressed as $t(q)=re(q)+\phi(pu(q))$ for some $\phi\in O(3)$. Similarly from any $q\in \mathbb{H}_{\Delta}$  we have  $l_q\in SO(4)$ , where $l_q:\mathbb{H}\rightarrow \mathbb{H}$, such that $l_q(x)=qx$, where we view $\mathbb{H}$ as $\R^4$ with the quaternion product.[Po]

\subsubsection{$\C$-inner product spaces}
We begin with $V$ as a $\C$-vector space, and we will define a $\C$-inner product space as a pair $(V,<,>)$ in the following manner.

\begin{defn}
$(V,<,>)$ is a $\C$ inner product space if $V$ is a $\C$ vector space , with bilinear map $<,>:V\times V\rightarrow \C$, such that for all $x_1,x_2,y_1,y_2\in V$ and $z_1,z_2\in\C$ the following are satisfied :
\begin{itemize}
\item $<z_1x_1+z_2 x_2,y_1>=z_1<x_1,y_1>+z_2<x_2,y_1>$
\item $<x_1,z_1 y_1+z_2y_2>=\bar{z_1}<x_1,y_1>+\bar{z_2}<x_1,y_2>$
\item $<y,x>=\overline{<x,y>}$, for all $x,y\in V$
\item For all non zero $x\in V$ ,$<x,x>$ is positive.

\end{itemize}
\end{defn}

For the $\C^n$ we define the inner product $<z,w>=\sum_{i=1}^n z_i\bar{w}_i$. For an inner product space $(V<,>)$ we naturally induce a norm $||x||=\sqrt{<x,x>}$ , along with the metric $d(x,y)=||x-y||$.

\begin{defn}
Any finite dimensional inner product space $(E,<,>)$ is called a \textbf{Hermitiain space}.

\end{defn}
Now for any Hermitian space we have an orthonormal basis $e_1,...,e_n$ such that  $<e_i,e_j>=\delta_{ij}$,  and

$$<\sum_{i=1}^n z_i e_i,\sum_{i=1}^n w_i e_i>=\sum_{i=1}^n z_i \bar{w_i}$$, thus implying that $(V,<,>)\cong \C^n$ , with the usual inner product. Just like the orthogonal group for $\R$-spaces $V$ , we can assign a group of linear isometries to inner product spaces.
\begin{defn}
The group of linear isometries of $(V,<,>)$ is called the \textbf{unitary group} denoted $U(V)$. Now for $\C^n$ with the usual inner product defined above we denote $U(\C^n)$ with $U(n)$. Now $SU(n):=\{h\in U(n):det\  h=1\}$ , is called the \textbf{special unitary group}.
\end{defn}

\subsubsection{Orthogonal groups as Matrix groups.}
\begin{defn}
A \textbf{Lie Group} is a group $(G,\cdot)$ , is a Manifold $G$ where $G$ is a set with a group operation $\cdot$ so that the maps $\mu:G\times G\rightarrow G$ and $inv:G\rightarrow G$ such that $\mu(g_1,g_2)=g_1\cdot g_2$ and $inv(g)=g^{-1}$ are smooth maps.

\end{defn}

For $V=\R^n$ we can think of the group invertable linear transformation $Gl(\R^n)$  in terms of matrices since for every invertable transformation we can associate an invertable matrix in the usual way. That is we can think of invertable linear transformation in $Gl(\R^n)$  as an invertable matrix in $Gl(n,\R)=\{A\in Mat(n,\R):det(A)\neq 0\}$ (we will use $Gl(n;\R)$ to denote the matrix group). $Gl(n,\R)$ is a Lie group since it can be viewed as an open subset of $Mat(n;\R)$ , thus inheriting a manifold structure, it inherits its group structure from the operation of matrix multiplication. The maps multiplication and inversion can be viewed as a polynomial map and a  quotient of polynomial maps respectively. Moreover $Gl(n,\C)$,$Gl(n,\mathbb{H})$ are also Lie groups, as well as some of the more common closed matrix subgroups of the generalized linear group over a field.[R]

 Since $O(V)$ is itself an invertible linear transformation that fixes the quadratic form, we can think of $O(V)$ as a matrix subgroup of orthogonal matrices in $Gl(n;\R)$. To make this distinction we will denoted it $$O(n;\R)=\{A\in Gl(n;\R):A^tA=I_{n}\}$$.
Where $A^t$ is the transpose of $A$. Since $O(n;\R)$ is a closed subgroup of the well known Lie group $Gl(n;\R)$ it is itself a Lie group. The special linear group can also be viewed as a closed matrix subgroup of $Gl(n;\R)$ and thus as the Lie group $SO(n;\R)=\{A\in O(n;\R):\det A=1\}$.The unitary group is defined as  ;$U(n;\C)=\{A\in Gl(n;\C):A^*A=I\}$, where $A^*=z_{ij}^*=\bar{z_{ji}}$ is the complex conjugate transpose. The unitary group is a closed subgroup of the Lie group $Gl(n;\C)$ and is thus a Lie group.The special unitary group are all the matrices in the unitary group with determinant 1 , $SU(n)=\{A\in U(n):\det\ A=1\}$, this is also a closed subgroup of $U(n)$ and is thus a compact Lie group[Bu]. The symplectic group is  a matrix group with coefficents in the quaternion algebra $\mathbb{H}$, defined as , $Sp(n;\mathbb{H})=\{A\in Gl(n;\mathbb{H}):A^*A=I_n\}$, where $A^*=q_{ij}^*=\bar{q_{ji}}$. The symplectic group can also be viewed in terms of $\C$ defined as 

$$Sp(n,\C)=\{A\in Gl(2n,\C):A^tJA=J\}$$, where $J=\scriptsize{\left[\begin{array}{cccc} 0&-I_n \\ I_n & 0\end{array}\right]}$. The determinant map restricted to these matrix subgroups provides the following projections[Po]; 

\begin{itemize}
\item $O(p,q)\rightarrow \{\pm 1\}$
\item $Gl(n,\R)\rightarrow \R^*$
\item $U(p,q)\rightarrow S^1$
\item $Gl(n,\C)\rightarrow \C^*$
\end{itemize} 

These matrix Lie groups are manifolds in their own right and we provide the dimensions over their respective fields in the following table [Po].

\begin{center}
\begin{tabular}{ |c|c| } 
 \hline
 Lie Group $G$ & $\dim G$\\ 
\hline
 $O(p,q)$ & $\frac{n^2-n}{2}$   \\ 
 
  $O(n,\C)$ & $n(n-1)$   \\ 
  
   $O(n,\mathbb{H})$ & $n(2n-1)$   \\ 
   
    $U(p,q)$ & $n^2$   \\ 
     $Sp(2n,\R)$ & $n(2n+1)$   \\ 
     
      $Sp(2n,\mathbb{C})$ & $2n(2n+1)$   \\

 \hline
\end{tabular}
\end{center}

We conclude this section with the following remark from [R].

\begin{remark}
When we view the special linear groups as linear transformations of determinant 1 the preserve a non degenerate form. For instance the rotation groups preserve the symmetric forms over $\R$ and $\C$, while  the symplectic groups preserve skew symmetric forms over $\R$ and $\C$ respectively , and finally $SU(p,q)$ preserves the hermitian form over $\C$ and $Sp(p,q)$ preserves the Hermitian form over $\mathbb{H}$. Where given $x=(x_1,...,x_n)$ and $y=(y_1,...,y_n)$ , depending of the signature we can view the symmetric form as $\pm x_1y_1\pm x_2y_2+...\pm x_ny_n$, the skew symmetric form as $(x_1y_2-x_2y_1)+...+(x_{2m-1}y_{2m}-x_{2m}y_{2m-1})$ (where $2m=n$), and  finally the Hermitian form being $\pm \bar{x_1}y_1\pm \bar{x_2}y_2+...\pm \bar{x_n}y_n$ .[R]

\end{remark}

\subsubsection{Lie Algebras of Lie groups}
A Lie group $(G,\cdot)$ is by definition a manifold , and thus it naturally has a tangent space $TG$. We will restrict our attention to the linear groups $G$ in this section.

\begin{defn}
Let $\mathfrak{g}$ be the tangent space to $G$ at $1_G$, where $$\mathfrak{g}:=\{A\in Mat(n,k): There\ exists\ a\ C^1\ curve\ , x(t),\ that\ lies\ in\ G\ such\ that\ x(0)=1_G,x'(0)=A\}$$. The tangent space $\mathfrak{g}$ that belongs to the linear group $G$ is called the \textbf{Lie algebra of G}.

\end{defn}

 $\mathfrak{g}$ is an $\R$-vector space closed under the Lie bracket $[X,Y]=XY-YX$ for $X,Y\in\mathfrak{g}$. Now if $G$ is a linear group and $\mathfrak{g}$ is its Lie algebra the exponential map $X\rightarrow e^{X}$ is a map from $\mathfrak{g}$ to $G$ [R]. We conclude this sections with the Lie algebras for some of our linear groups discussed earlier[Wa][R][Ha].

\begin{itemize}

\item For the Linear group $Gl(n,k)$ where $k=\R,\C,\mathbb{H}$ the Lie algebra $\mathfrak{g}=Mat(n,k)$.
\item For a $\C$-vector space ,$V$, the Lie algebra for the automorphism group $Aut(V)$ , is $\mathfrak{g}=End(V)$, where the Lie bracket is defined as $[\phi,\psi]=\phi\circ\psi-\psi\circ\phi$, for $\phi,\psi\in End(V)$.
\item For the orthogonal group $O(n,k)$ , and the special orthogonal group $SO(n,k)$ the Lie algebra is the matrix algebra  of skew symmetric matrices;
$$\mathfrak{o}(n,k)=\mathfrak{so}(n,k)=\{A\in Mat(n,k): X^t=-X\}$$, where $k=\R,\C$.
\item For the unitary group $U(n,\C)$ the corresponding  Lie algebra is $$\mathfrak{u}(n)=\{
A\in Mat(n,\C):A^*=-A\}$$.
\item For the special unitary group $SU(n,\C)$ its corresponding Lie algebra is $$\mathfrak{su}(n)=\{
A\in Mat(n,\C):A^*=-A,tr A=0\}$$.
\item  The symplectic group viewed as a group over $\R$ has the Lie algebra 
$$\mathfrak{sp}(n,\R)=\{X\in Mat(2n,\R):JX^tJ=X, J=\scriptsize{\left[\begin{array}{cccc} 0&-I_n \\ I_n & 0\end{array}\right]}\}$$, while the  symplectic group over $\C$  is defined in the same manner $$\mathfrak{sp}(n,\C)=\{X\in Mat(2n,\C):JX^tJ=X; J=\scriptsize{\left[\begin{array}{cccc} 0&-I_n \\ I_n & 0\end{array}\right]}\}$$. When we view the symplectic group as group over $\mathbb{H}$ its Lie algebra is defined as; $$\mathfrak{sp}(n,\mathbb{H})=\mathfrak{sp}(n,\C)\cap \mathfrak{u}(2n)$$. 

\end{itemize}

\section{Basics of Clifford Algebras}
\subsection{Definition of a Clifford algebra}
For a quadratic space $(V,q)$  we can define a tensor algebra ,$T(V,q)=\bigoplus_{k=0}^{\infty}V^k$, where $V^k=V\oplus ...\oplus V$ (k-times) and $V^0$ is any scalar field and $V^1=V$. Addition in $T(V)$ is component wise addition and multiplication can be defined as a concatenation in the graded components as the  map $m:V^k\times V^m\rightarrow V^{k+m} $ such that $m(x_1\otimes...\otimes x_k,y_1\otimes...\otimes y_m)=x_1\otimes...\otimes x_k\otimes y_1\otimes...\otimes y_m$. 

\begin{exmp}
The exterior algebra $\wedge V$ has a correspondence with the tensor algebra $T(V)$ with a slight modification to the product. The correspondence being $T(V)$ correspond directly with the elements in $\wedge V$ via correspondence $x_1\otimes...\otimes x_k$ with $x_1\wedge ...\wedge x_k$, and for a multiplication of $k$- vector $x$ and $m$-vector $y$ , we have the wedge product which produces $k+m$ vector  $x\wedge y$ such that $x\wedge x=0$ and $x\wedge y=-y\wedge x$. Also $\wedge^0 V=$ is the scalar field for $V$ and $\wedge^1 V=V$.
\end{exmp}

Now having a tensor algebra $T(V)$ , and a quadratic vector space $(V,q)$ we can define the ideal $I_q=<v\otimes v+q(x)1_{T(V)}>$ for $v\in V$. We can define the \textbf{Clifford algebra of $(V,q)$} to be the quotient of the tensor algebra $T(V)$ with the ideal $I_q$ , and denote the Clifford algebra $Cl(V,Q)=T(V)/I_q$. Clifford algebras $Cl(V,q)$ have the natural vector space embedding $(V,q)\rightarrow Cl(V,q)$ , such that the image of the embedding is $V$ itself under the injective canonical projection $\pi:T(V)\rightarrow Cl(V,Q)$[LM]. The Clifford algebra $Cl(V,Q)$ happens to be generated by $V$ and its scalar field , by the relation $||v||^2=-q(v)\cdot 1_{Cl(V,q)}$,  for all $v\in (V,q)$. Allowing us to define the Clifford algebra in terms of $V$, without reference to the  tensor algebra.Clifford algebras $Cl(V,q)$ inherit a natural graded algebra structure, which is  actually a $\Z_2$ grading $Cl(V,q)=Cl(V,q)^+\oplus CL(V,q)^-$ with $Cl(V,q)^{\pm}$ signifying the odd and even parts respectively ( we define this grading in terms of involutions in the next section). We now define the Clifford algebra in terms of an orthonormal basis for $(V,q)$.

\begin{defn} 
\label{definition}
Let $V$ be an $n$-dimensional vector space with basis $\{e_1, ..., e_n\}$ over a field $\mathbb{F}$ such that char$(\mathbb{F}) \neq 2$ and let $q$ be a quadratic form. We define a \textbf{Clifford algebra} $C(V, Q)$ to be an algebra generated by $\{\gamma_1, ..., \gamma_n\}$, subject to the
relation $\gamma_i \gamma_j + \gamma_j \gamma_i = 2 B(e_i, e_j) \cdot 1_C$ for all $i,j = 1, ..., n$, where $1_C$ is the unit in the Clifford algebra. 
\end{defn}

\
Now if we focus on $\R$-orthogonal spaces  $V=\R^{p,q}$, with the usual quadratic form induced from the bi-linear form $B$ such that $B(e_i,e_i)=q(e_i)$ for  and $B(e_i,e_j)=0$ for $i\neq j$ we can define the Clifford algebra in terms of relations with respects to the orthonormal basis $e_1,...e_{p+q}$. Where  the matrix of the non-degenerate symmetric bilinear form $B$ can be expressed as the block matrix
$\scriptsize{\left[\begin{array}{cc} +I_p & 0\\ 0& -I_q\end{array}\right]},$
where $n = p+q$. We call $(p,q)$ the \textbf{signature} of the Clifford algebra and denote it as $Cl_{p,q}$. From Definition \ref{definition}, we can obtain these relations between the generators $\g_1, ..., \g_n$ in $Cl_{p,q}$:

\begin{itemize}

\item $\g_i^2 = 1$ for $i = 1, ..., p$

\item $\g_i^2 = -1$ for $i = p+1, p+2,..., n$

\item $\g_i \g_j = -\g_j \g_i$ for all $i \neq j$.
\end{itemize}

Using the relations above for $\g_1,...,\g_n$, we can form a {canonical basis} of $2^n$ elements, given by the monomials:
\begin{equation}
\label{canonicalbasis}
\{1\} \cup \{\g_{i_1} \cdots \g_{i_k} : i_1 < i_2 < \cdots < i_k, k = 1, ..., n\}.
\end{equation}
For notational purposes, we write the product of generators with their indices together; for instance $\g_{ij}:=\g_i \g_j$. This means that for any $u\in Cl_{p,q}$ we can uniquely write $u=\sum_{i_1<...<i_k}\alpha_{i_1...i_k} \cdot \g_{i_1...i_k}$ for some $\alpha_{i_1...i_k}\in\R$.\\


\begin{exmp}
Some elementary examples of universal Clifford algebras in low dimensions  are of the following;
\begin{itemize}
\item $Cl(0)=\R$.
\item The Clifford algebra $Cl_{0,1}$ is isomorphic to $\C$ , just by identifying $\g_1\mapsto i$.
\item The Clifford algebra $Cl_{0,2}$ is isomorphic to the quaternions $\mathbb{H}$ with $\g_1\mapsto i$,$\g_2\mapsto j$, and $\g_{12}\mapsto k$. This will preserve the multiplicative structure of $\mathbb{H}$. With the $im(\mathbb{H})=\{c_1\g_1+c_2\g_2+c_3\g_{12}:c_1,c_2,c_3\in\R\}$. Notice that we could have looked at $SO(3)$ and $SO(4)$ with respects to the Clifford algebra $Cl_{0,2}$.
\item $Cl_{0,3}$ is a semi simple Clifford algebra  that ends up as the direct sum of two copies of the quaternions $\mathbb{H}$, that is  $Cl_{0,3}\cong\mathbb{H}\oplus \mathbb{H}$. The identification of the basis elements of $Cl_{0,3}$ and  $\mathbb{H}\oplus\mathbb{H}$ is given in the following table.

\begin{center}
\begin{tabular}{ |c|c| } 
 \hline
 $Cl_{0,3}$ & $\mathbb{H}\oplus\mathbb{H}$ \\ 
\hline
 $1$ & $(1,1)$   \\ 
 
  $\g_1$ & $(i,-i)$   \\ 
  
   $\g_2$ & $(j,-j)$   \\ 
   
    $\g_3$ & $(k,-k)$   \\ 
     $\g_{12}$ & $(k,k)$   \\ 
     
      $\g_{13}$ & $(-j,-j)$   \\
       $\g_{23}$ & $(i,i)$   \\ 
        $\g_{123}$ & $(-1,1)$   \\

 \hline
\end{tabular}
\end{center}
\end{itemize}

\end{exmp}

\subsection{Three Involutions for $Cl_{p,q}$ and the quadratic norm $||.||$}

For a given Clifford algebra $Cl(V)$ the set of automorphism $Aut(Cl(V))$ is actually equal to $O(V)$ for that same quadratic space $(V,q)$, that is $Aut(Cl(V))=O(V)$[Har].

 We will define three key involutions for given Clifford algebras [Gar]. 
 \begin{itemize}
 \item The \textbf{principal involution} is the natural  extension to the Clifford algebra $Cl(V)$  ,$\phi_t$ ,of the isometry $t:V\rightarrow V$, such that $t(x)=-x$ for all $x\in V$. Giving us the following diagram ; 
\xymatrix{ V \ar[d]^{\iota} \ar[r]^{t} &V\ar[d]^{\iota}\\  Cl(V)  & Cl(V) \ar[l]^{\phi_{t}}}

Such that $\phi_t(\g_{i_1..i_r})=(-1)^r \g_{i_1...i_r}$. 

The principal involution also defines the natural grading of the Clifford algebra $Cl(V)$ where 
$Cl(V)^+=\{x\in Cl(V):\phi_t(x)=-x\}$, and $Cl(V)^-=\{x\in Cl(V):\phi_t(x)=x\}$. 

So elements on the standard basis thought of as $k$-vectors are in $Cl(V)^+$ if $k$ is odd and $Cl(V)^-$ if $k$ is even.  

\item When we look at the identity  $id:V\rightarrow V$ , we can view the extension $\phi_{id}:Cl(V)\rightarrow Cl(V)^{opp}$  as the \textbf{reversal} giving us to the following diagram ,

\xymatrix{ V \ar[d]^{\iota} \ar[r]^{id} &V\ar[d]^{\iota}\\  Cl(V)^{opp}  & Cl(V) \ar[l]^{\phi_{id}}}

Where multiplication in $Cl(V)^{opp}$ is defined as $x\cdot y=yx$. Thus for a given $r$ -vector $\g_{i_1...i_r}$ we have $\phi_{id}(\g_{i_1...i_r})=\g_{i_r....i_1}=(-1)^{\frac{r(r-1)}{2}}\g_{i_1...i_r}$.  
\item The \textbf{conjugation involution} is the composition of the principal involution and the reversion involution . That is $\bar{x}=\phi_t(\phi_{id}(x))=\phi_{id}(\phi_t(x))$ for a given $x\in Cl(V)$.

Given an element in the basis for $Cl(n)_{p,q}$  we have $\overline{\g_{i_1...i_r}}=(-1)^{\frac{r(r+1)}{2}}\g_{i_1...i_r}$. For an invertable element $\g\in Cl(V)$ , the quadratic norm in $Cl(V)$  is defined to be $||\g||=\bar{\g}\g$.

The quadratic norm $||.||$ has the property that $||x||=||\bar{x}||=||\phi_t(x)||=||\phi_{id}(x)||$ [Ga].

\end{itemize}

\subsection{Groups of Clifford Algebras}

 The group of units of a Clifford algebra is defined as  $Cl(V)^* =\{x\in Cl(V): \exists y\in Cl(V)\ s.t\ xy=1\}$.  The set $Cl(V)^*$ is  a group under multiplication. In fact $Cl(V)^*$ is a lie group of dimension $2^{dim V}$[LM]. We can use the quadratic norm to help us identify whether  in element is in $Cl(V)^*$, via the fact that $x\in Cl(V)^*$ iff $||x||\in Cl(V)^*$[Ga]. We now go on to analyze the Clifford group and $O(V)$.
 \subsubsection{Clifford group}
 
 Given the group of units $Cl(V)^*$, we can define a natural action on the Clifford algebra $Cl(V)$ via the adjoint map associated to a unit $u\in Cl(V)^*$. The adjoint map is an automorphism on $Cl(V)$, where   $Ad_u(x)=uxu^{-1}$, for $x\in Cl(V)$.
 
  Giving us the natural action $Cl(V)$ in  the following sense;
 
 $$Ad_u:Cl(V)^*\rightarrow Aut(Cl(V))$$, via $u\rightarrow Ad_u $[LM]. We can slightly modify this to get the action $$\hat{Ad_u}\rightarrow Aut(Cl(V))$$ such that $\hat{Ad_u}(v)=\phi_t(u)xu^{-1}$ for all $x\in Cl(V)$, and use it to define the the Clifford group.\textbf{The Clifford group}, $\Gamma(V)$, can be defined as the invertable elements that stabilize $V$ under the action $\hat{Ad_u}$[Ga]. That is $$\Gamma (V)=\{\ u\in Cl(V)^*:\hat{Ad_u}(x)\in V,\ \forall x\in V\}$$. The action establishes a correspondence between the invertable elements $u$ and orthogonal transformations in $V$,  the correspondence being that for every element $u\in\Gamma(V)$, we have an induced orthogonal transformation  $\rho_{u}:V\rightarrow V$ , such that $x\rightarrow \hat{Ad_u}(x)$, that is $\rho_{u}\in O(V)$[Po]. This means that we can view the Clifford group of $V$  as the set of all $u\in Cl(V)^*$ such that $\rho_{u}\in O(V)$; that is $\Gamma(V):=\{ x\in Cl(V)^* : \rho_{u}\in O(V)\}$. The map $\Gamma(V)\rightarrow Gl(V)$ . where $u\mapsto \rho_u$ is a natural representation for the Clifford group, whose image is  $O(V)$[Me]. In fact with this natural representation we have the following exact sequence for an $\R$ quadratic space $V$;
 
 $$1\rightarrow \R^* \rightarrow \Gamma(V)\rightarrow O(V)\rightarrow 1$$.
 
 If $(V,q)$ is a quadratic vector space with $q$ is non-degenerate , then every $\sigma\in O(V)$ is represented by an element $u$  in the Clifford group of $V$. The Clifford group also relates to the quadratic norm established in the previous section in the following manner [Po]: 

\begin{itemize}
\item $\g\in\Gamma(V)$, implies that $||\g||\in\R$
\item $||1_{Cl(V)}||=1$
\item For $\g,\g'\in\Gamma(V)$; $||\g\g'||=||\g||\cdot ||\g'||$.
\item For $\g\in\Gamma(V)$, $||\g||\ne 0$ and $||\g^{-1}||=\dfrac{1}{||\g||}$
\end{itemize}

\subsubsection{Pin and Spin groups}
For The Clifford group $\Gamma(V)$ we can view two special  subgroups $Pin(V)$ and $Spin(V)$.  The subgroup $Pin(V)$ is defined as the subgroup of $Cl(V)^*$ that is generated by unit vectors. That is 
$Pin(V)=\{x\in Cl(V)^*:x=u_1....u_k, ||u_i||=\pm 1\}$ . We can view the group $Pin(V)$ as the kernel of the Clifford  norm restricted to the Clifford group , that is $Pin(V)=\ker \{||.||:\Gamma(V)\rightarrow \R^*\}$[Me]. Using the map $\hat{Ad}$ we get the \textbf{twisted adjoint representation} of the Pin group , 
$$\hat{Ad}:Pin(V)\rightarrow O(V)$$, where $\hat{Ad}(u)=\hat{Ad_u}$, which gives us the exact sequence ;

$$1\rightarrow \Z_2\rightarrow Pin(V)\xrightarrow{\hat{Ad}} O(V)\rightarrow 1$$; where the image $\hat{Ad}(Pin(V))$ is a normal subgroup of $O(V)$[LM].If the dimension of $V$ is odd then $\hat{Ad}(Pin(V))=SO(V)$[Tr].


The Spin groups ,$Spin(V)$ ,are very similar to the Pin groups , but they have a correspondence relationship to the group of rotations on $V$, $SO(V)$. $Spin(V)$  is a subgroup of $Pin(V)$ generated by an even number of unit vectors, this implies that we can view $Spin(V)=Pin(V)\cap Cl(V)^-$[Ga]. We can also define it as $Spin(V)=\{x\in Cl(V)^*:x=u_1....u_k, ||u_i||=\pm 1,\ k\ is\ even\}$ if $V$ is a non degenerate orthogonal space. Restricting the adjoint representation to $Spin(V)$ gives us a surjective map $\hat{Ad}: Spin(V)\rightarrow SO(V)$, which  give the short exact sequence; [Har].

$$1\rightarrow \Z_2\rightarrow Spin(V)\xrightarrow{\hat{Ad}}SO(V)\rightarrow 1$$.
Where $\hat{Ad}(Spin(V))$ is a normal subgroup of $O(V)$. In fact $Spin(n)$ defines a universal double cover for $SO(n)$ for $n>2$[Me].In terms of classification we have the following theorem[Po].



\begin{thm}
Let $V$ be a non degenerate $\R$- orthogonal space  with $\dim V\le 5$; then 

$$Spin(V)=\{\g\in Cl(V)^{-}: ||\g||=\pm 1\}$$.

\end{thm}

 We now go on to define the Projective Clifford group in a similar manner to Pin and Spin.

\subsubsection{Projective Clifford groups }
	
When we quotient the Clifford group $\Gamma(V)$	by $\R^*$  we get the group called the\textbf{ projective Clifford group} , denoted $Proj(Cl(V))=\Gamma(V)/\R^*$. The projective Clifford group has the property that it is isomorphic to $O(V)$, via the surjective group map $\Gamma(V)\rightarrow O(V)$. If we restrict the map to the even Clifford group, denoted $\Gamma(V)^{-}=\Gamma(V)\cap Cl(V)^{-}$,  we have the surjective group map $$\Gamma(V)^{-}\rightarrow SO(V)$$, which establishes the isomorphism $\Gamma(V)^{-}/\R^*\cong SO(V)$. The group $Proj(Cl(V))^{-}=\Gamma(V)^{-}/\R^*$ is called the \textbf{even projective Clifford group}. The projective Clifford groups allows to view $O(V)$ and $SO(V)$ in terms of equivalences with the usual correspondence $[\g]\mapsto \rho_{\g}$. For more information on the projective Clifford groups can be found on [Po].

\subsubsection{More on the Pin and Spin groups}

For a given $u\in V$ for a quadratic vector space $V$ we can view $v\rightarrow R_u(v)=-uvu^{-1}$ for a fixed $u\in V$ as a reflection in the hyper plane orthogonal to $u$ [Tr]. In $\R^n$ the reflection $\{\rho_{\vec{v}}:\R^n\rightarrow \R^n\}\in O(n)$ is viewed  as the reflection in the orthogonal complement in $\R^n$ through the origin spanned by the vector $\vec{v}\in \R^n$. Also in  $\R^n$ we can view rotations $\sigma\in SO(V)$ as a composite of an even number of hyperplane reflections [Ab 1]. Now for a give $k$-plane in $V$, $u\in G(k,V)$ , $R_u x=-x$ if $x\in\ span\  u$, and $R_u x=x$ if $x\notin\ span\  u$. More generally if $u\in G(k,V)$ is a unit non degenerate $k$-plane in $V$ then $R_u(x)=\hat{Ad_u}(v)$, for all $x\in V$ [Har]. If we view the $Pin(V)$ strictly as a space generated by  unit vectors then in terms of the Grassmanian algebra, then $Pin(V)$ is generated by $G(1,V)$, and the spin group $Spin(V)$ can be viewed as elements $a\in Cl(V)^*$ such that $a=u_1...u_k$ where each $u_j\in G(2,V)$[Har]. Given a  quadratic vector space $(V,q)$ we can view the  generalized  unit sphere as $S=\{v\in V:q(v)=\pm 1\}$. The groups $Pin(V)$ and $Spin(V)$ can then be seen as the groups  generated by the generalized unit sphere[LM]. $Spin(n)$ is the double cover of $SO(n)$ since  $\rho_{u}(x)=\rho_{-u}(x)$, that is $\pm u$ induce the same map in $SO(n)$. [Ab 1]. We now list some elementary  results for $Spin(n)$ [Po]. 

\begin{itemize}
\item $Spin(1)\cong O(1)$
\item $Spin(2)\cong U(1)=\{z\in \C:|z|=1\}\cong S^1$.
\item  $Spin(3)\cong Sp(1)\cong S^3$
\item $Spin(4)\cong Sp(1)\times Sp(1)\cong S^3\times S^3$.
\item $Spin(5)\cong Sp(2)$
\item $Spin(6)\cong SU(4)$

\end{itemize}
If $V$ is a complex space with the standard bi-linear form we have complex Pin  and Spin groups defined in the usual way, that is $Pin(n,\C)$and $Spin(n,\C)$ respectively, where we have the following results for the complex Spin groups in low dimensions .[Tr]

\begin{itemize}

\item $Spin(2,\C)\cong \C^*$.
\item  $Spin(3,\C)\cong Sl_2(\C)=\{A\in Mat(2,\C): det A=1\}$
\item $Spin(4,\C)\cong Sl(2,\C)\times Sl_2(\C)$.
\item $Spin(5,\C)\cong Sp(4,\C)$
\item $Spin(6,\C)\cong Sl_4(\C)$

\end{itemize}
\subsubsection{Lie Algebra of the Spin group }

The group of units $Cl(V)^*$ of a Clifford algebra is a lie group of dimension $2^{dimV}$, with an associated lie algebra 
$\mathfrak{cl}(V)=Cl(V)$. $\mathfrak{cl}(V)$ is just the standard Clifford algebra of dimension $2^{dimV}$ with the Lie bracket $[x,y]=xy-yx$ for $x,y\in Cl(V)$[LM].  The $Ad$ map defines a natural action on the group of automorphisms in $Cl(V)$. When it comes to the Lie algebra $\mathfrak{cl}(V)$ we define the map 

$$ ad:\mathfrak{cl}(V)\rightarrow Der(Cl(V))$$, such that $ad(y):Cl(V)\rightarrow k$, such that $ad(y)(x)=[y,x]$.This map naturally identifies elements in the Clifford algebra with elements in the derivations of the Clifford algebra[LM]. $Spin(n)$ is a lie group whose dimension is $\frac{(n)(n-1)}{2}$. In fact its the maximal compact subgroup of $Spin(n;\C)$. The Lie algebra of $Spin(V)$ is $\mathfrak{o}(V)$ [Me].

\newpage
\section{Spinors and Representations of Clifford Algebras}
\subsection{Isomorphisms as associative algebras}

The Clifford algebras $Cl_{p,q}$ can be viewed isomorphically, as associative algebras, to matrix algebras. It  should be noted that some of the structure of $Cl_{p,q}$ does not carry to the matrix algebras they are associated with. We state the following modulo 8 classification for Clifford algebras of the form $Cl_{p,q}$.[Po];
\begin{itemize}
\item  $Cl_{p,q}\cong Mat(2^{\frac{n}{2}},\mathbb{\R})$  if $q-p= 0,6 \mod 8$
\item $Cl_{p,q}\cong Mat(2^{\frac{n-1}{2}},\mathbb{\C})$  if $q-p= 1,5 \mod 8$. 
\item $Cl_{p,q}\cong Mat(2^{\frac{n-2}{2}},\mathbb{\mathbb{H}})$  if $q-p= 2,4 \mod 8$.

\item In the cases that $q-p= 3\mod 8$, $Cl_{p,q}$ is semi simple and  as an associative algebra we have ; $Cl_{p,q}\cong Mat(2^{\frac{n-3}{2}},\mathbb{\mathbb{H}})\oplus Mat(2^{\frac{n-3}{2}},\mathbb{\mathbb{H}})$ 

\item  If $q-p= 7\mod 8$, $Cl_{p,q}$ is semi simple and  as associative algebras $Cl_{p,q}\cong Mat(2^{\frac{n-1}{2}},\mathbb{\mathbb{\R}})\oplus Mat(2^{\frac{n-1}{2}},\mathbb{\mathbb{\R}})$

\end{itemize}










\subsection{Spinor spaces and representations}
We begin this section with the simple definition of a representation of an algebra $A$. 
 \begin{defn}
A \textbf{representation} of an algebra $A$ is a vector space $V$, together with a algebra homomorphism $\rho :V\rightarrow End(V)$. A representation is said to be \textbf{minimal} (or of minimal dimension) if there does not exist any faithful representation of $A$ of lower dimension. 
\end{defn}
The isomorphisms stated above allow us to view   $Cl_{p,q}$ as an $\R$-algebra of the endomorphisms of a $D$-linear space $D^m$, where $D$ is one of the division algebras (or direct sum) and $m$ is the corresponding dimension mentioned above in the modulo 8 classification. The space $D^m$ is called the \textbf{Spinor space} of $\R^{p,q}$.If $Cl_{p,q}$ is simple and we have any irreducible representation of the form $\rho:Cl_{p,q}\rightarrow End_D(V)$, then $V$ is a spinor space and elements in $V$ are called spinors.  It is worth noting that the spinor spaces $D^m$ are identifiable with minimal left ideals in $Cl_{p,q}$[Ab]. We will construct a special kind of minimal left ideal in $Cl_{p,q}$ in the following section  that will be denoted as $S_{p,q}$, called the spinor space for $\R^{p,q}$. The representations in $Cl_{p,q}\rightarrow End_D(S_{p,q})$ will be called the Spinor representations for Clifford algebra $Cl_{p,q}$.

\begin{exmp}
In space-time physics the Clifford algebras often chosen are $Cl_{3,1}$ or $Cl_{1,3}$, with the isomoprphisms  $Cl_{3,1}\cong Mat(2,\mathbb{H})$ and $Cl_{1,3}\cong Mat(4,\R)$ as associative algebras. The spinor spaces associated to spaces of matrix algebras are $\mathbb{H}^2$ and $\R^4$. Where the elements in the spinor space $\mathbb{H}^2$ are known as Dirac spinors , while elements in the spinor space $\R^4$ are known as Majorana spinors. The complexification of these spaces is isomorphic to $Mat(4,\C)$ and the elements of the spinor space $\C^4$ are known as Weyl spinors [Ab].
\end{exmp}
We conclude this sections with some lower dimensional spinor spaces for Clifford algebras $Cl_{p,q}$[Ga].

\begin{center}
\begin{tabular}{ |c|c| } 
 \hline
Signature $(p,q)$ & Spinor space \\ 
\hline
 $(0,0)$ &  $\R$  \\ 
 
  $(0,1)$ & $\C$   \\ 
  
   $(0,2)$ & $\mathbb{H}$   \\ 
   
    $(0,3)$ & $\mathbb{H}\oplus \mathbb{H}$   \\ 
     $(0,4)$ & $\mathbb{H}^2$   \\ 
     $(1,0)$ &  $\R\oplus \R$  \\ 
 
  $(2,0)$ & $\R^2$   \\ 
  
   $(3,0)$ & $\C^2$   \\ 
   
    $(4,0)$ & $\mathbb{H}^2$   \\ 
     $(1,1)$ & $\R^2$   \\ 
      $(1,2)$ &  $\C^2$  \\ 
 
  $(1,3)$ & $\mathbb{H}^2$   \\ 
  
   $(2,1)$ & $\R^2\oplus\R^2$   \\ 
   
    $(2,2)$ & $\R^4$   \\ 
     $(3,1)$ & $\C^2$   \\ 
$(5,0)$ &  $\mathbb{H}^2\oplus \mathbb{H}^2$  \\ 
 
  $(0,5)$ & $\C^4$   \\ 
  
   $(6,0)$ & $\mathbb{H}^4$   \\ 
   
    $(0,6)$ & $\R^8$   \\ 
     $(7,0)$ & $\C^8$   \\ 
     $(0,7)$ & $\R^8\oplus \R^8$\\
 \hline
\end{tabular}
\end{center}
\newpage

\subsection{Construction of Spinor spaces via primitive idempotents}
It is a well known fact that minimal left ideals of matrix algebras are generated by \textbf{primitive idempotents}, elements in the algebra whose square is itself and cannot be expressed as the sum of two annihilating idempotents .The minimal ideal generated by such idempotents consists of matrices with the characteristic of having all zero columns expect for one which has a one somewhere and all other zeros.  For any Algebra $A$ the most simple idempotents to construct are of the form $\dfrac{1+a}{2}$ for elements $a\in A$ such that $a^2=1$, which also implies that $\dfrac{1-a}{2}$ is an idempotent[Ab]. For $Cl_{p,q}$  we can choose $k$ commuting involutions ,$\g_{\alpha_1},...,\g_{\alpha_k}$, which give us $2k$ idempotent elements in the Clifford algebra of the form $\dfrac{1\pm \g_{\alpha_j}}{2}$  who decompose the Clifford algebra $Cl_{p,q}$ into a direct sum of left ideals;$$Cl_{p,q}=Cl_{p,q}\dfrac{1+ \g_{\alpha_1}}{2}\oplus Cl_{p,q}\dfrac{1- \g_{\alpha_1}}{2}\oplus...\oplus Cl_{p,q}\dfrac{1+ \g_{\alpha_k}}{2}\oplus Cl_{p,q}\dfrac{1- \g_{\alpha_k}}{2}$$. Facts about this decomposition can be found in [LH],[Ab3]. The idempotents in the  ideal decomposition of $Cl_{p,q}$ are pairwise annihilating and sum to $1$.
\begin{remark}
 The value $k$, for the number of commuting involutions,  is computed as $k= q-r_{q-p}$ , and it depends on Randon -Hurwitz numbers $r_{q-p}$ which have the following properties: $r_0=0,r_1=1,r_2=2,r_3=2,r_j=3$ where $4\le j\le 7$, $r_{i+8}=r_i+4$, $r_{-1}=-1$, and $r_{-i}=1-i+r_{i-2}$[Po].

\end{remark}



From the $k$ commuting involutions we can construct the primitive idempotent $F=(\dfrac{1+\g_{\alpha_1}}{2})...(\dfrac{1+\g_{\alpha_k}}{2})$, which gives us the minimal left ideal $Cl_{p,q}F$, which forms a $Cl_{p,q}$ -module via the left multiplication action $Cl_{p,q}\times Cl_{p,q}F\rightarrow Cl_{p,q}F$.This is a module of dimension $2^{n-k}$, where elements in $Cl_{p,q}F$ can be seen as equivalence classes of order $k$. 

\begin{defn}
The Clifford module $Cl_{p,q}F$ of dimension $2^{n-k}$ constructed from the $k$ commuting involutions that make up the primitive idempotent $F$ will be called the spinor space $S_{p,q}$ .

\end{defn}

\begin{remark}
 If $\g F\in S_{p,q}$ under the Clifford norm we have the property that $||\g F||=\g F\overline{\g F}=\g F\bar{F}\bar{\g}$, now since $F\bar{F}=0$,  we have $||\g F||=0$. Thus the quadratic norm $||.||$ vanishes in the spinor space $S_{p,q}$[Ga].
\end{remark} 
 Given the spinor space $S_{p,q}$ we define the division ring $\mathbb{E}=FS_{p,q}=\{F\g F:\g\in Cl_{p,q}\}$. The non trivial elements with respects to multiplication are the elements of the form $\g F\mathbb{E}=\mathbb{E}\g F$. The division ring has the following isomorphisms;[Lo];[Ab3].

\begin{itemize}

\item $\mathbb{E}\cong \R$ if $q-p=0,6,7 \mod 8$,$\mathbb{E}\cong \C$ if $q-p=1,5\mod 8$,
$\mathbb{E}\cong \mathbb{H}$ if $q-p=2,3,4 \mod 8$

\end{itemize}
 When $q-p=3\mod 4$ ,we have $S_{p,q}\oplus \hat{S_{p,q}}$ with the right $\mathbb{E}\oplus \mathbb{E}$-module structure are called the double spinor spaces for $Cl_{p,q}$.Where $\hat{S_{p,q}}=\{\phi_t(\g F):\g F\in S_{p,q}\}$. More about these division algebras can be found in $[Lo]$. 
\begin{exmp}
For the Clifford algebra $Cl_{3,1}$ we can define $F=\dfrac{1+\g_{1}}{2}\dfrac{1+\g_{24}}{2}$, with the minimal left ideal $S_{3,1}=Cl_{3,1}F$ as our spinor space. In the spinor space $S_{3,1}$  we have the following equalities;
$F=\g_1F=\g_{24}F=\g_{124}F$; $\g_2F=\g_4F=-\g_{12}F=-\g_{14}F$; $\g_3F=-\g_{13}F=-\g_{234}F=\g_{1234}F$, and $\g_{23}F=-\g_{34}F=\g_{123}F=-\g_{134}F$.
 Thus for the spinor space $S_{3,1}$ we have the $\R$-basis $\{F,\g_2F,\g_3F,\g_{23}F\}$, and we can identify $S_{3,1}$ with $\R^4$ as a $\R$-vector space. We can then view the Majorna spinors as elements in $S_{3,1}$  giving us the isomorphism $Cl_{3,1}\cong End_{\mathbb{R}}(S_{3,1})$.For the generators $\g_1,\g_2,\g_3$ and $\g_4$ we have the following associated matrices in $Mat(4,\R)$ ; 
 $$\g_{1}\leftrightarrow \scriptsize{\left[\begin{array}{cccc} 1 & 0 & 0 & 0\\ 0 & -1 & 0 & 0\\ 0 & 0 & -1 & 0 \\ 0 & 0 & 0 & 1 \end{array}\right]}\;\;\g_{2}\leftrightarrow \scriptsize{\left[\begin{array}{cccc} 0 & 1 & 0 & 0\\ 1& 0 & 0 & 0\\ 0 & 0 & 0 & 1 \\ 0 & 0 & 1 & 0 \end{array}\right]}$$

$$\g_{3}\leftrightarrow \scriptsize{\left[\begin{array}{cccc} 0 & 0 & 1 & 0\\ 0 & 0 & 0 & -1\\ 1& 0 & 0 & 0 \\ 0 & -1 & 0 & 0 \end{array}\right]}\;\;
\g_{4}\leftrightarrow \scriptsize{\left[\begin{array}{cccc} 0 & -1 & 0 & 0\\ 1& 0 & 0 & 0\\ 0 & 0 & 0 & -1\\ 0 & 0 & 1 & 0 \end{array}\right]}.$$
\end{exmp}

\begin{exmp}
For the Clifford algebra $Cl_{1,3}$ we can define $F=\dfrac{1+\g_{12}}{2}$, with the minimal left ideal $S_{1,3}=Cl_{1,3}F$ as our spinor space. In the spinor space $S_{1,3}$  we have the following equalities;
$F=\g_{12}F$,$\g_{1}F=\g_{2}F$ ,$\g_{3}F=\g_{123}F$,$\g_{4}F=\g_{124}F$,$\g_{13}F=\g_{23}F$,$\g_{14}F=\g_{24}F$,
$\g_{34}F=\g_{1234}F$, and $\g_{134}F=\g_{234}F$. Thus for the spinor space $S_{1,3}$ we have the $\R$-basis $\{F,\g_1F,\g_3F,\g_4F,\g_{13}F,\g_{14}F,\g_{34}F,\g_{134}F\}$, and we can identify $S_{1,3}$ with $\R^8$ as a vector space. Moreover we can give $S_{1,3}$ an $\mathbb{H}$ structure , by identifying $\g_{3}F$,$\g_{4}F,\g_{34}F$ with $i,j,k$ respectively. Thus we can view $S_{1,3}$ as an  $\mathbb{H}$-space  of the form  $S_{1,3}=\{F,\g_1F\}\otimes_{\R} \mathbb{H}$. Allowing us identify $S_{1,3}$ with $\mathbb{H}^2$ ,that is $S_{1,3}\cong\mathbb{H}^2$ as vector spaces. Thus dirac spinors in $\mathbb{H}^2$ can be viewed as elements in $S_{1,3}$ and we have the isomorphism $Cl_{1,3}\cong End_{\mathbb{H}}(S_{1,3})$ . For the generators $\g_1,\g_2,\g_3$ and $\g_4$ we have the following identifications with matrices in $Mat(2,\mathbb{H})$; 
$$\g_{1}\leftrightarrow \scriptsize{\left[\begin{array}{cccc} 0&1 \\ 1 & 0\end{array}\right]}\;\;
\g_{2}\leftrightarrow \scriptsize{\left[\begin{array}{cccc} 0 & -1 \\ 1 & 0\end{array}\right]}$$

$$\g_{3}\leftrightarrow \scriptsize{\left[\begin{array}{cccc} i & 0 \\ 0 & -i\end{array}\right]}\;\;
\g_{4}\leftrightarrow \scriptsize{\left[\begin{array}{cccc} j & 0 \\ 0 & -j\end{array}\right]},$$

\end{exmp}

\begin{exmp}
If we complexify the Clifford algebra $Cl_{1,3}$ , we have $\C l_{1,3}=\C\otimes_{\R} Cl_{1,3}$, we use the primitive idempotent $F=\dfrac{1}{4}(1+\g_{1234}i+\g_{12}+\g_{34}i)$, to construct the minimal left ideal $\C l_{1,3}F$, which is a spinor space that we will denote $S^{\C}_{1,3}$  that will be the complexification of our spinor space $S_{1,3}$ from the previous example. The equivalences in $S^{\C}_{1,3}$ will be $F=\g_{12}F=\g_{34}iF=\g_{1234}iF$, $\g_1F=\g_{134}iF=\g_{2}F=\g_{234}iF$, $\g_3F=-\g_{4}iF=-\g_{124}iF=\g_{123}F$, $\g_{13}F=-\g_{14}iF=\g_{23}F=-\g_{24}iF$. Giving us the  $\C$- basis $F,\g_1F,\g_3F,\g_{13}F$ for  $S^{\C}_{1,3}$ giving us the canonical vector space isomorphism  $S^{\C}_{1,3}\cong \C^4$. Where $\C^4$ is known as the Weyl spinor space, and elements known as Weyl spinors. We can also get the spinor representation for $\C l_{1,3}$ in $End_{\C}(S^{\C}_{1,3})$. Giving us  have the following identifications for $\g_1,\g_2,\g_3,\g_4$ :
 $$\g_{1}\leftrightarrow \scriptsize{\left[\begin{array}{cccc} 0 & 1 & 0 & 0\\ 1 & 0 & 0 & 0\\ 0 & 0 & 0 & 1 \\ 0 & 0 & 1 & 0 \end{array}\right]}\;\;\g_{2}\leftrightarrow \scriptsize{\left[\begin{array}{cccc} 0 & -1 & 0 & 0\\ 1& 0 & 0 & 0\\ 0 & 0 & 0 & -1 \\ 0 & 0 & 1 & 0 \end{array}\right]}$$

$$\g_{3}\leftrightarrow \scriptsize{\left[\begin{array}{cccc} 0 & 0 & -1 & 0\\ 0 & 0 & 0 & 1\\ 1& 0 & 0 & 0 \\ 0 & -1 & 0 & 0 \end{array}\right]}\;\;
\g_{4}\leftrightarrow \scriptsize{\left[\begin{array}{cccc} 0 & 0 & i & 0\\ 0 & 0 & 0 & -i\\ i& 0 & 0 & 0 \\ 0 & -i & 0 & 0 \end{array}\right]}.$$
\end{exmp}

\subsection{Generalized Spinor spaces and representations }
We begin this section with some well known representations of Clifford Algebras in lower dimensions.

\begin{exmp}

\item $Cl(1)_{0,1}\cong \C$ 

$$\g_1\leftrightarrow i$$.

\item $Cl(2)_{0,2}\cong \mathbb{H}$ 

$$\g_{1}\leftrightarrow i\;\;
\g_{2}\leftrightarrow j\;\;\g_{12}\leftrightarrow k$$

\item $Cl(3)_{1,2}$ 

\subitem Representation of generators in $Mat(2,\C)$ 

$$\g_{1}\leftrightarrow \scriptsize{\left[\begin{array}{cccc} 0 & 1 \\ 1& 0\end{array}\right]}\;\;
\g_{2}\leftrightarrow \scriptsize{\left[\begin{array}{cccc} 0& -1 \\ 1 & 0 \end{array}\right]}\;\;\g_{3}\leftrightarrow \scriptsize{\left[\begin{array}{cccc} i& 0 \\ 0 & i \end{array}\right]}$$

\item $Cl(4)_{2,2}$
\subitem Representation of generators in $Mat(4,\R)$
$$\g_{1}\leftrightarrow \scriptsize{\left[\begin{array}{cccc} 0 & 1 & 0 & 0\\ 1& 0 & 0 & 0\\ 0 & 0 & 0 & 1 \\ 0 & 0 & 1 & 0 \end{array}\right]}\;\;
\g_{2}\leftrightarrow \scriptsize{\left[\begin{array}{cccc} 1 & 0 & 0 & 0\\ 0 & -1 & 0 & 0\\ 0 & 0 & -1 & 0 \\ 0 & 0 & 0 & 1 \end{array}\right]}$$

$$\g_{3}\leftrightarrow \scriptsize{\left[\begin{array}{cccc} 0 & 0 & 1 & 0\\ 0 & 0 & 0 & -1\\ 1& 0 & 0 & 0 \\ 0 & -1 & 0 & 0 \end{array}\right]}\;\;
\g_{4}\leftrightarrow \scriptsize{\left[\begin{array}{cccc} 0 & -1 & 0 & 0\\ 1& 0 & 0 & 0\\ 0 & 0 & 0 & -1\\ 0 & 0 & 1 & 0 \end{array}\right]}.$$

\end{exmp}

The number of involutions $k$ that we use to construct the Spinor spaces $S_{p,q}$  can also be given a modulo $8$-structure  , and this number $k$ depends on congruence classes mod 8 in the following way:
\begin{itemize}
\item $q-p \equiv 0,1,3,5,6\mod 8$, $k=\lfloor{\frac{n}{2}}\rfloor $.
\item $q-p \equiv 2,4\mod 8$, $k=\lfloor{\frac{n}{2}}\rfloor-1$.
\item $q-p \equiv 7\mod 8$, $k=\lfloor{\frac{n}{2}}\rfloor+1$.
\end{itemize}

One can easily see that this calculation follows directly from the isomorphisms established in 4.1.
\begin{defn}
A \textbf{generating set} for a Clifford algebra $Cl_{p,q}$ is  the set of the $k$ commuting involutions established above. The involutions in the generating set are independent insofar that no involution in the generating set is the Clifford product of any of the other involutions in the set.

\end{defn}
The generating sets will provide us with a useful tool to construct the primitive idempotents $F$ needed to construct the general Spinor space $S_{p,q}$ for quadratic space $\R^{p,q}$. We now will provide a detailed account on how to construct the generating sets for a given signature $(p,q)$.
\subsubsection{Building generating sets}

We begin with a Clifford algebra generated from a quadratic space with the signatures $(p,p+a+8l)$, where $l$ is a non negative integer and $a\in \{1,2,3,4,5,6,7\}$. For a fixed $p$ we will define the following set $P_2\subset\N^2$ with $p$ elements as $P_2 := \{(1,p+1),(2,p+2),...,(p,2p)\}$. In a similar manner for a fixed $p$ we will define the set, of $4l$ elements, $P_4\subset \N^4$ in the following manner; $P_4 :=\bigcup_{i=1}^{l}\{(2p+8(i-1)+1,2p+8(i-1)+2,2p+8(i-1)+3,2p+8(i-1)+4),(2p+8(i-1)+1,2p+8(i-1)+2,2p+8(i-1)+5,2p+8(i-1)+6),(2p+8(i-1)+1,2p+8(i-1)+2,2p+ 8(i-1)+7,2p+8(i-1)+8),(2p+8(i-1)+1,2p+8(i-1)+3,2p+8(i-1)+5,2p+8(i-1)+7)\}$	, Where the set $P_4$ is empty if $l=0$.

Similarly for signatures of the form $(q+a+8l,q)$,where $l$ is a non negative integer and $a\in\{1,2,3,4,5,6,7\}$. For a fixed $q,a$ and $l$ we define the following set, of $q$ elements ,$Q_2\subset \N^2$ in the following manner $ Q_2:= \{(1,q+8l+a+1),(2,q +8l+a+2),...,(q,2q+8l+a)\}$. Similarly for a fixed $a$ we define the following subset $Q_4\subset \N^4$, of $4l$ elements in the following manner; $Q_4 :=\bigcup_{i=1}^l\{(q+8(i-1)+1,q+8(i-1)+2,q+8(i-1)+3,q+8(i-1)+4),(q+8(i-1)+1,q+8(i-1)+2,q+8(i-1)+5,q+8(i-1)+6),(q+8(i-1)+1,q+8(i-1)+2,q+ 8(i-1)+7,q+8(i-1)+8),(q+8(i-1)+1,q+8(i-1)+3,q+8(i-1)+5,q+8(i-1)+7)\}$ . Where the set $Q_4$ is empty if $l=0$.

Due to the natural grading of a Clifford algebras we will use the sets $P_2$ (resp $Q_2$) and $P_4$ (resp $Q_4$ ) to generate a subset of $2$ and $4$ vectors respectively. The subsets generated by the sets $P_2$ and $P_4$ (resp $Q_2$ and $Q_4$ ) will be denoted $\{\g_{P_2}\}=\{\g_{\alpha} :\alpha\in P_2\}\subset Cl_{p,p+a+8l}$ and $\{\g_{P_4}\} = \{\g_{\alpha}:\alpha\in P_4\}\subset Cl_{p,p+a+8l}$ ( resp {$\{\g_{Q_2}\}=\{\g_{\alpha} :\alpha\in Q_2\} \subset Cl_{q+a+8l,q}$ and $\{\g_{Q_4}\}=\{\g_{\alpha} :\alpha\in Q_4\} \subset Cl_{q+a+8l,q}$). Moreover in the way the sets $P_2$ and $P_4$ (resp $Q_2$ and $Q_4$ ) are constructed, all elements in $\{\g_{P_2}\}$ and $\{\g_{P_4}\}$  (resp $\{\g_{Q_2}\}$ and $\{\g_{Q_4}\}$ ) are commuting involutions such that all $2$-vectors in $\{\g_{P_2}\}$  commute with $4$-vectors in  $\{\g_{P_4}\}$. (resp all $2$-vectors in $\{\g_{Q_2}\}$  commute with $4$-vectors in  $\{\g_{Q_4}\}$.) With the aid of the sets $P_2$, $P_4$ ($Q_2$ and $Q_4$ respectively) we construct the following table of generating sets for all $Cl_{p,q}$.

\begin{table}[ht]
\begin{center}
\begin{tabular}{|c|c|c|p{8 cm}|c|}
\hline 

Signature & $q-p\, (\mod 8)$ & $k$  &  Generating set \\ 
\hline 
\hline 
$(p,p+8l)$  & $0,1,2$ & $p+4l$ & $\{\g_{P_{2}}\}\cup \{\g_{P_{4}}\}$
\\ 

\hline 
$(p,p+8l+3)$  & $3$ & $p+4l+1$ & $\{\g_{P_{2}}\}\cup \{\g_{P_{4}}\}\cup \{\g_{2p+8l+1,2p+8l+2,2p+8l+3}\}$
\\ 
\hline 
$(p,p+8l+4)$  & $4$ & $p+4l+1$ & $\{\g_{P_{2}}\}\cup \{\g_{P_{4}}\}\cup \{\g_{2p+8l+1,2p+8l+2,2p+8l+3,2p+8l+4}\}$
\\ 
\hline 
$(p,p+8l+5)$  & $5$ & $p+4l+2$ & $\{\g_{P_{2}}\}\cup \{\g_{P_{4}}\}\cup\{\g_{2p+8l+1,2p+8l+2,2p+8l+3,2p+8l+4},$ 

$\g_{2p+8l+1,2p+8l+2,2p+8l+5}\}$
\\ 
\hline 
$(p,p+8l+6)$  & $6$ & $p+4l+3$ & $\{\g_{P_{2}}\}\cup \{\g_{P_{4}}\}\cup \{\g_{2p+8l+1,2p+8l+2,2p+8l+4}, $

$\g_{2p+8l+2,2p+8l+3,2p+8l+5}, \g_{2p+8l+3,2p+8l+4,2p+8l+6}\}$
\\ 
\hline 
$(p,p+8l+7)$  & $7$ & $p+4l+4$ & $\{\g_{P_{2}}\}\cup \{\g_{P_{4}}\}\cup \{\g_{{2p+8l+1,2p+8l+2,2p+8l+4}},$

$\g_{{2p+8l+2,2p+8l+3,2p+8l+5}},$

$\g_{{2p+8l+3,2p+8l+4,2p+8l+6}},\g_{{2p+8l+4,2p+8l+5,2p+8l+7}}\}$
\\ 
\hline 
\end{tabular} 
\end{center}
\caption{ Generating sets for Clifford algebras of signature $(p,q)$, $p\leq q$.}
\label{tab:table1}
\end{table}

\begin{table}[ht]
\begin{center}
\begin{tabular}{|c|c|c|p{8 cm}|c|}
\hline 

Signature & $q-p\, (\mod 8)$ & $k$  &  Generating set \\ 
\hline 
\hline 
$(q+8l,q)$  & $0$ & $q+4l$ & $\{\g_{Q_{2}}\}\cup \{\g_{Q_{4}}\}$
\\ 
\hline 
$(q+8l+1,q)$  & $4,5,6,7$ & $q+4l+1$ & $\{\g_{Q_{2}}\}\cup \{\g_{Q_{4}}\}\cup \{\g_{q+8l+1}\}$
\\ 

\hline 
$(q+8l+5,q)$  & $3$ & $q+4l+2$ & $\{\g_{Q_{2}}\}\cup \{\g_{Q_{4}}\}\cup\{\g_{q+8l+1,q+8l+2,q+8l+3,q+8l+4},\g_{q+8l+5}\}$
\\ 
\hline 
$(q+8l+6,q)$  & $2$ & $q+4l+2$ & $\{\g_{Q_{2}}\}\cup \{\g_{Q_{4}}\}\cup\{\g_{q+8l+1,q+8l+2,q+8l+3,q+8l+4},$

$\g_{q+8l+1,q+8l+2,q+8l+5, q+8l+6}\}$
\\ 
\hline 
$(q+8l+7,q)$  & $1$ & $q+4l+3$ & $\{\g_{Q_{2}}\}\cup \{\g_{Q_{4}}\}\cup\{\g_{q+8l+1,q+8l+2,q+8l+3,q+8l+4},$

$\g_{q+8l+1,q+8l+2,q+8l+5,q+8l+6}, \g_{q+8l+7}\}$
\\ 
\hline 
\end{tabular} 
\end{center}
\caption{ Generating sets for Clifford algebras of signature $(p,q)$, $p\geq q$.}
\label{tab:table2}
\end{table}
\newpage 
\subsubsection{General $\R$-spinor basis for $S_{p,q}$.}

The $\R$ dimension of the spinor space $S_{p,q}$ constructed from our $k$ involutions will naturally be $2^{n-k}$. Thus we have a natural $\R$ basis for $S_{p,q}$ for each of our generating sets established in our previous section. Each of the  $2^{n-k}$ basis elements of $S_{p,q}$ can be thought of as an equivalence class  of order $k$. To be more precise about what these identified equivalence classes contains, we use the fact that our generating sets are composed of $k$-vectors, where $0\le k\le 4$, and that for any generator $\g_I$ in the generating set we have $\g_IF=F$.Notice that $P_2$(resp $Q_2$), which generates the set $\g_{P_2}$(resp $\g_{Q_2}$) ,is composed of pairs in $\N^2$, also notice that for a given point $(i,j)\in P_2$ (resp $(i,j)\in Q_2$), where $i\in\{1,...,p\}$ (resp  $i\in\{1,...,q\}$ ) and $j\in\{p + 1,...,2p\}$(resp $j \in\{q + 8l + a + 1,...,2q + 8l + a \}$), we have $\g_iF=\pm \g_jF$ .This implies that for any point $(i,j)\in P_2$ (resp $(i,j)\in Q_2$) we have $\g_iF=\pm \g_jF$ , depending on the signature of the Clifford algebra, in $S_{p,q}$. This implies that the set $\hat{P_2} = \{1,...,p\}\subset \N$ (resp $\hat{Q_2}= \{1,...,q\}\subset \N$) will yield the set  $\{\g_{\hat{P}_2}\}=\{\g_{\alpha}F :\alpha \in \hat{P}_2 \}$ ( resp $\{\g_{\hat{Q}_2}\}=\{\g_{\alpha}F :\alpha \in \hat{Q}_2 \}$ )Where the set $\{\g_{\hat{P}_2}\}$ (resp $\{\g_{\hat{Q}_2}\}$) contains a total of $p$ elements (resp $q$ elements) in the $\R$-basis for $S_{p,q}$.

For the $4$-vectors that make up our generating sets identified with the sets $P_4,Q_4\subset \N^4$, we establish similar equivalences amongst the points that make up $P_4$ and $Q_4$ , these equivalences will be denoted by the sets $\hat{P_4}$,$\hat{Q_4}$ . Thus any point $(n_1,n_2,n_3,n_4)\in P_4$ (resp $(n_1,n_2,n_3,n_4)\in Q_4$ ) will correspond to a $4$-vector in the canonical basis with equivalences; $\g_{(n_1,n_2,n_3,n_4)}F=\pm F$ , $\g_{(n_1)}F= \pm \g_{(n_2,n_3,n_4)}F$ , $\g_{(n_2)}F= \pm \g_{(n_1,n_3,n_4)}F$ ,$\g_{(n_3)}F=\pm \g_{(n_1,n_2,n_4)}F$ ,$\g_{(n_4)}F=\pm \g_{(n_1,n_2,n_3)}F$ . Where the sign is dependent of whether the signature of the Clifford algebra is positive or negative definite. Now if we examine the construction of the set $P_4$ we will notice that for a fixed $i$ we have the following 4 points in $\N^4$ ;
$\big<(2p + 8(i-1) + 1,2p + 8(i-1) + 2,2p + 8(i-1) + 3,2p + 8(i-1) + 4), (2p + 8(i-1) + 1,2p + 8(i-1) + 2,2p + 8(i-1) + 5,2p + 8(i-1) + 6), (2p + 8(i-1) + 1,2p + 8(i-1) + 2,2p + 8(i-1) + 7,2p + 8(i-1) + 8), (2p + 8(i-1) + 1,2p + 8(i-1) + 3,2p + 8(i-1) + 5,2p + 8(i-1) + 7)\big>$.

 Using the equivalence relation mentioned above , for each fixed $i$ we get the following $4$ points in $\N$, $\{2p+8(i-1)+ 1,2p + 8(i-1) + 2,2p + 8(i-1) + 3,2p + 8(i-1) + 5\}$, to resemble equivalence classes in the $\R$-basis of $S_{p,q}$ .Thus from the initial set $P_4\subset \N^4$ used to construct the generating sets , we  have the  set $\hat{P}_4 =\bigcup _{i=1}^l\{(2p+8(i-1)+1),(2p+8(i-1)+ 2),(2p+8(i-1)+3),(2p+8(i-1)+5)\}\subset \N$ ,  that we use to construct the $\R$-basis for  $S_{p,q}$ , where $\hat{P}_4$ is empty if $l = 0$. Using the same method for the set $Q_4\subset \N^4$  we define the set $\hat{Q}_4 =\bigcup_{i=1}^{l} \{(q+8i+1),(q+8i+2),(q+8i+3),(q+8i+5)\} \subset \N$ where the set is empty if $l = 0$. Since our idempotent $F$ is made up primarily of elements in $P_2,P_4$ (resp $Q_2,Q_4$ ) the $\R$-spinor basis for $S_{p,q}$ can be constructed primarily using $\hat{P}_2,\hat{P}_4$ (resp $\hat{Q}_2,\hat{Q}_4$) where we will generate the sets $\{\g_{\hat{P}_2}\}=\{\g_{\alpha}F :\alpha \in \hat{P}_2 \}$,$\{\g_{\hat{P}_4}\}=\{\g_{\alpha}F :\alpha \in \hat{P}_4 \}$ ( resp $\{\g_{\hat{Q}_2}\}=\{\g_{\alpha}F :\alpha \in \hat{Q}_2 \}$,$\{\g_{\hat{Q}_4}\}=\{\g_{\alpha}F :\alpha \in \hat{Q}_4 \}$). Notice that the sets $\{\g_{\hat{P}_2}\},\{\g_{\hat{P}_4}\}$ ( resp  $\{\g_{\hat{Q}_2}\},\{\g_{\hat{Q}_4}\})$  are all composed of one vectors. With the natural grading of Clifford algebras and a fixed $k$ we define the set of $k$-vectors in the following manner ; 
$\{\g_{\hat{P}_{\mu}}\}^k=\{\g_{\alpha}F :\alpha=(l_{i_1},...,l_{i_k}), l_{i_j} \in \hat{P}_{\mu},l_{i_1}<...<l_{i_k} \}$ (resp  $\{\g_{\hat{Q}_{\mu}}\}^k=\{\g_{\alpha}F :\alpha=(l_{i_1},...,l_{i_k}), l_{i_j} \in \hat{Q}_{\mu},l_{i_1}<...<l_{i_k} \}$ )
  where $\mu$  is either $2$ or $4$. The set $\{\g_{\hat{P}_{\mu}}\}^0$ (resp $\{\g_{\hat{Q}_{\mu}}\}^0$) will be the identity element $F$ in $S_{p,q}$. With the construction of $\{\g_{\hat{P}_{\mu}}\}^k$ (resp $\{\g_{\hat{Q}_{\mu}}\}^k$ )we will generate an $\R$ basis for  $S_{p,q}$. The constructions of the $\R$ -basis of  $S_{p,q}$ with respects to the appropriate signatures are given by tables 3 and 4. 
  
\begin{remark}  

 With the construction of sets $\{\g_{\hat{P}_2}\},\{\g_{\hat{P}_4}\}$ ( resp  $\{\g_{\hat{Q}_2}\},\{\g_{\hat{Q}_4}\}$  the union 
 $\{\g_{\hat{P}_2}\}\bigcup \{\g_{\hat{P}_4}\}$ ( resp  $\{\g_{\hat{Q}_2}\}\bigcup \{\g_{\hat{Q}_4}\}$ , will consist of a total of $p+4l$ elements in the $\R$-basis for signatures of the form $(p,p + 8l + a)$ (resp will consist of a total of $q+4l$ elements in the $\R$-basis for signatures of the form $(q + 8l + a,q)$). When $a$ is non trivial, the remaining elements in the $\R$-basis of $S_{p,p+8l+a}$ (resp $S_{q+8l+a,q}$ ) will be generated by a subset $A$ of the remaining one vectors $\g_{2p+8l+1}F,...,\g_{2p+8l+a}F$ (resp  $\g_{q+8l+1}F,...,\g_{q+8l+a}F$). With $\{\g_{\hat{P}_2}\}\bigcup \{\g_{\hat{P}_4}\}\bigcup A$ ( resp  $\{\g_{\hat{Q}_2}\}\bigcup \{\g_{\hat{Q}_4}\}\bigcup A$) and the natural grading of Clifford algebras, the union of the sets of all $k$-vectors $\{\{\g_{\hat{P}_2}\}\bigcup \{\g_{\hat{P}_4}\}\bigcup A\}^k$ ( resp  $\{\{\g_{\hat{Q}_2}\}\bigcup \{\g_{\hat{Q}_4}\}\bigcup A\}^k$)will give us the $\R$-basis for the spinor space $S_{p,q}$.

\end{remark}
\begin{remark}
We will denote all remaining elements $\g_{\alpha}F\in A$ from the previous remark , as $\hat{\g_{\alpha}}$.
\end{remark}
\newpage
\begin{table}[ht]
\begin{center}
\begin{tabular}{|p{5 cm}|p{9.5 cm}|c}
\hline 

Signature & $\R$-basis for $S_{p,q}$   &  $\dim_{\R}S_{p,q}$ \\ 
\hline 
\hline 
$(p,p+8l)$  & $\bigcup_{k=0}^{p+4l} \{\{\g_{\hat{P}_2}\}\bigcup \{\g_{\hat{P}_4}\}\}^k$  & $2^{p+4l}$
\\ 

\hline 
$(p,p+8l+1)$  & $\bigcup_{k=0}^{p+4l+1}\{\{\g_{\hat{P}_2}\}\bigcup \{\g_{\hat{P}_4}\}\bigcup \{\hat{\g}_{2p+8l+1}\}\}^k$  & $2^{p+4l+1}$
\\ 
\hline 
$(p,p+8l+2)$, 

$(p,p+8l+3)$ & $\bigcup_{k=0}^{p+4l+2}\{\{\g_{\hat{P}_2}\}\bigcup \{\g_{\hat{P}_4}\}\bigcup \{\hat{\g}_{2p+8l+1},\hat{\g}_{2p+8l+2}\}\}^k$  &$2^{p+4l+2}$
\\ 
\hline 
$(p,p+8l+4)$
 $(p,p+8l+5)$,
 $(p,p+8l+6)$, 
 $(p,p+8l+7)$          &   $\bigcup_{k=0}^{p+4l+3}\{\{\g_{\hat{P}_2}\}\bigcup \{\g_{\hat{P}_4}\}\bigcup \{\hat{\g}_{2p+8l+1},\hat{\g}_{2p+8l+2},\hat{\g}_{2p+8l+3}\}\}^k$ & $2^{p+4l+3}$
\\ 

\hline 
\end{tabular} 
\end{center}
\caption{  $\R$-basis for Clifford algebras of signature$ (p,p + a + 8l)$
}
\label{tab:table3}
\end{table}

\begin{table}[ht]
\begin{center}
\begin{tabular}{|p{4.75 cm}|p{10.23 cm}|c}
\hline 

Signature & $\R$-basis for $S_{p,q}$   &  $\dim_{\R}S_{p,q}$ \\ 
\hline 
\hline 
$(q+8l,q)$ ,
$(q+8l+1,q)$ & $\bigcup_{k=0}^{q+4l} \{\{\g_{\hat{Q}_2}\}\bigcup \{\g_{\hat{Q}_4}\}\}^k$  & $2^{q+4l}$
\\ 

\hline 
$(q+8l+2,q)$  & $\bigcup_{k=0}^{q+4l+1}\{\{\g_{\hat{Q}_2}\}\bigcup \{\g_{\hat{Q}_4}\}\bigcup \{\hat{\g}_{q+8l+2}\}\}^k$  & $2^{q+4l+1}$
\\ 
\hline 
$(q+8l+3,q)$ & $\bigcup_{k=0}^{q+4l+2}\{\{\g_{\hat{Q}_2}\}\bigcup \{\g_{\hat{Q}_4}\}\bigcup \{\hat{\g}_{q+8l+2},\hat{\g}_{q+8l+3}\}\}^k$  &$2^{q+4l+2}$
\\ 
\hline 
$(q+8l+4,q)$ & $\bigcup_{k=0}^{q+4l+3}\{\{\g_{\hat{Q}_2}\}\bigcup \{\g_{\hat{Q}_4}\}\bigcup \{\hat{\g}_{q+8l+2},\hat{\g}_{q+8l+3},\hat{\g}_{q+8l+4}\}\}^k$  &$2^{q+4l+3}$
\\
\hline 
$(q+8l+5,q)$ & $\bigcup_{k=0}^{q+4l+3}\{\{\g_{\hat{Q}_2}\}\bigcup \{\g_{\hat{Q}_4}\}\bigcup \{\hat{\g}_{q+8l+1},\hat{\g}_{q+8l+2},\hat{\g}_{q+8l+3}\}\}^k$  &$2^{q+4l+3}$
\\
\hline 
$(q+8l+6,q)$,

$(q+8l+7,q)$ & $\bigcup_{k=0}^{q+4l+4}\{\{\g_{\hat{Q}_2}\}\bigcup \{\g_{\hat{Q}_4}\}\bigcup \{\hat{\g}_{q+8l+1},\hat{\g}_{q+8l+2},\hat{\g}_{q+8l+3},\hat{\g}_{q+8l+5}\}\}^k$  &$2^{q+4l+4}$
\\

\hline 
\end{tabular} 
\end{center}
\caption{Table 4. $\R$-basis for Clifford algebras of signature $(q + a + 8l,q)$.  }
\label{tab:table4}
\end{table}
\newpage
From the given $\R$-basis  $S_{p,q}$ we can proceed to find the $\C$ and $\mathbb{H}$ basis for the appropriate signatures $(p,q)$, so that our spinor spaces $S_{p,q}$ correspond with the isomorphisms established in section 4.1. To do this however we will need to identify basis elements in $S_{p,q}$ with generator $i$ in $\C$, and with generators $i,j,k$ in $\mathbb{H}$. 

\subsection{$\C$-Basis for $S_{p,q}$}

When we have an associative $\R$-algebra with a unit element ,we know by [Po,88] that any $2$-dimensional sub-algebra generated by an element $e_0$ such that $e_0^2=-1$ is isomorphic to $\C$.From tables 3 and 4 we can see all the $\R$-basis constructions for $S_{p,q}$.  Thus for the Clifford algebras where $q-p\equiv 1,5 \mod 8$   we will go from the $\R$-basis identified with $S_{p,q}$ to a $\C$-basis. To do this will make use of the fact that negative definitive elements in the spinor basis behave just like $i$ in $\C$. When we view $S_{p,q}$ as an $\R$-vector space , we can easily make it into a $\C$-vector space by identifying a basis element that is negative definite with the generator $i$ in $\C$. A identification of the sort would allow us to view $S_{p,q}$ as a $\C$-vector space instead of just an $\R$-vector space. Thus for the $\R$-basis for $S_{p,q}$ where $q-p\equiv 1,5 \mod 8$ we will choose negative definite spinor basis element $\g_{\beta}F$ to identify with $i$, thus establishing a canonical vector space isomorphism between $S_{p,q}$ and  $\C\otimes S_{p,q-1}$. The construction of the $\C$ -basis for signatures $(p,q)$ where $q-p\equiv 1,5 \mod 8$  are given by table 5.

\begin{table}[ht]
\begin{center}
\begin{tabular}{|p{5 cm}|p{10 cm}|p{2cm}}
\hline 

Signature & $\C$-basis for $S_{p,q}$   &  $\dim_{\C}S_{p,q}$ \\ 
\hline 
\hline 
$(p,p+8l+1)$  & $\bigcup_{k=0}^{p+4l} \{\{\g_{\hat{P}_2}\}\bigcup \{\g_{\hat{P}_4}\}\}^k\otimes \{F,\hat{\g}_{2p+8l+1}\}$, where we identify $\hat{\g_{2p+8l+1}}\leftrightarrow i$ & $2^{p+4l}$
\\ 

\hline 
$(p,p+8l+5)$  & $\bigcup_{k=0}^{p+4l+2} \{\{\g_{\hat{P}_2}\}\bigcup \{\g_{\hat{P}_4}\}\bigcup \{\g_{2p+8l+1},\g_{2p+8l+2}\}\}^k\otimes \{F,\hat{\g}_{2p+8l+3}\}$, where we identify $\hat{\g_{2p+8l+3}}\leftrightarrow i$   & $2^{p+4l+1}$
\\ 
\hline

$(q+8l+3,q)$ & $\bigcup_{k=0}^{q+4l+1} \{\{\g_{\hat{Q}_2}\}\bigcup \{\g_{\hat{Q}_4}\}\bigcup \{\hat{\g}_{q+8l+2}\}\}^k\setminus\{\hat{\g}_{1,q+8l+2}\}\otimes \{F,\hat{\g}_{1,q+8l+2}\}$, where we identify $\hat{\g_{1,q+8l+2}}\leftrightarrow i$  &$2^{2+4l+1}$
\\ 
\hline

$(q+8l+7,q)$ & $\bigcup_{k=0}^{q+4l+3} \{\{\g_{\hat{Q}_2}\}\bigcup \{\g_{\hat{Q}_4}\}\bigcup \{\hat{\g}_{q+8l+1},\hat{\g}_{q+8l+2},\hat{\g}_{q+8l+3}\}\}^k\setminus\{\hat{\g}_{1,q+8l+5}\}\otimes \{F,\hat{\g}_{1,q+8l+5}\}$, where we identify $\hat{\g_{1,q+8l+5}}\leftrightarrow i$  & $2^{2+4l+1}$
\\ 

\hline 
\end{tabular} 
\end{center}
\caption{  $\C$-basis for $S_{p,q}$ where $q-p\equiv 1,5\mod 8$
}
\label{tab:table5}
\end{table}

\subsection{$\mathbb{H}$-basis for $S_{p,q}$}

From [Po,88] we know that for any associative $\R$-algebra with a unity element $1$ , any $4$-dimensional sub algebra generated by a set $\{e_0,e_1\}$ of mutually anti commuting elements such that $e_0^2=e_1^2=-1$ is isomorphic to $\mathbb{H}$.  The Spinor representations of $Cl_{p,q}$ where $q-p\equiv 2,3,4 \mod 8$  involve an $\mathbb{H}$- spinor space.  To find an $\mathbb{H}$-basis for $S_{p,q}$ we will have to find three negative definite elements in the spinor basis $\g_{\alpha}F,\g_{\beta}F,\g_{\zeta}F$ that have following multiplicative structure :

$$ \g_{\alpha}F\g_{\beta}F=\g_{\zeta}F,\g_{\beta}F\g_{\zeta}F=\g_{\alpha}F,\g_{\zeta}F\g_{\alpha}F=\g_{\beta}F,\g_{\beta}F\g_{\alpha}F=-\g_{\zeta}F,\g_{\alpha}F\g_{\zeta}F=-\g_{\beta}F,\g_{\zeta}F\g_{\beta}F=-\g_{\alpha}F.$$
 When we find these three basis elements we identify them with the generators $i,j,k$ in $\mathbb{H}$ in the following manner; $\g_{\alpha}F\leftrightarrow i$ , $\g_{\beta}F\leftrightarrow j$, $\g_{\zeta}F\leftrightarrow k$. It then follows that the $\R$-basis in  $S_{p,q}$ can be expressed as an $\mathbb{H}$-basis. An identification of this sort will allow us to establish the vector space isomorphism  $S_{p,q}\cong \mathbb{H}\otimes S_{p,q-3}$. The construction of the $\mathbb{H}$-basis for signatures $(p,q)$ where $q-p\equiv 2,3,4\mod 8$  are given by table 6.

\begin{table}[ht]
\begin{center}
\begin{tabular}{|p{5 cm}|p{10.5 cm}|p{2cm}}
\hline 

Signature & $\mathbb{H}$-basis for $S_{p,q}$   &  $\dim_{\mathbb{H}}S_{p,q}$ \\ 
\hline 
\hline 
$(p,p+8l+2)$, 

 $(p,p+8l+3)$ & $\bigcup_{k=0}^{p+4l} \big[\{\{\g_{\hat{P}_2}\}\bigcup \{\g_{\hat{P}_4}\}\bigcup\{\hat{\g}_{2p+8l+1},\hat{\g}_{2p+8l+2}\}\}\setminus\{\hat{\g}_{2p+1},\hat{\g}_{2p+2}\}\big]^k\otimes \{F,\hat{\g}_{2p+1},\hat{\g}_{2p+2},\hat{\g_{2p+1,2p+2}}\}$, where we identify $\hat{\g_{2p+1}}\leftrightarrow i$,$\hat{\g_{2p+2}}\leftrightarrow j,\hat{\g_{2p+1,2p+2}}\leftrightarrow k$ & $2^{p+4l}$
\\ 

\hline 
, 

 $(p,p+8l+4)$ & $\bigcup_{k=0}^{p+4l+2} \big[\{\{\g_{\hat{P}_2}\}\bigcup \{\g_{\hat{P}_4}\}\bigcup\{\hat{\g}_{2p+8l+1},\hat{\g}_{2p+8l+2}\}\}\setminus\{\hat{\g}_{2p+1},\hat{\g}_{2p+2}\}\big]^k\otimes \{F,\hat{\g}_{2p+1},\hat{\g}_{2p+2},\hat{\g_{2p+1,2p+2}}\}$, where we identify $\hat{\g_{2p+1}}\leftrightarrow i$,$\hat{\g_{2p+2}}\leftrightarrow j,\hat{\g_{2p+1,2p+2}}\leftrightarrow k$ & $2^{p+4l+2}$
\\ 

\hline 
 
$(q+8l+4,q)$ & $\bigcup_{k=0}^{q+4l+1} \big[\{\{\g_{\hat{Q}_2}\}\bigcup \{\g_{\hat{Q}_4}\}\bigcup\{\hat{\g}_{q+8l+2},\hat{\g}_{q+8l+3},\hat{\g}_{q+8l+4}\}\}\setminus\{\hat{\g}_{q+2},\hat{\g}_{q+3}\}\big]^k\otimes \{F,\hat{\g}_{q+1,q+2},\hat{\g}_{q+2,q+3},\hat{\g_{q+1,q+3}}\}$, where we identify $\hat{\g_{q+1,q+2}}\leftrightarrow i$,$\hat{\g_{q+2,q+3}}\leftrightarrow j,\hat{\g_{q+1,q+3}}\leftrightarrow k$ & $2^{q+4l+1}$
\\ 

\hline

$(q+8l+5,q)$ & $\bigcup_{k=0}^{q+4l+1} \big[\{\{\g_{\hat{Q}_2}\}\bigcup \{\g_{\hat{Q}_4}\}\bigcup\{\hat{\g}_{q+8l+1},\hat{\g}_{q+8l+2},\hat{\g}_{q+8l+3}\}\}\setminus\{\hat{\g}_{q+2},\hat{\g}_{q+3}\}\big]^k\otimes \{F,\hat{\g}_{q+1,q+2},\hat{\g}_{q+2,q+3},\hat{\g_{q+1,q+3}}\}$, where we identify $\hat{\g_{q+1,q+2}}\leftrightarrow i$,$\hat{\g_{q+2,q+3}}\leftrightarrow j,\hat{\g_{q+1,q+3}}\leftrightarrow k$ & $2^{q+4l+1}$
\\ 
\hline
$(q+8l+6,q)$ & $\bigcup_{k=0}^{q+4l+2} \big[\{\{\g_{\hat{Q}_2}\}\bigcup \{\g_{\hat{Q}_4}\}\bigcup\{\hat{\g}_{q+8l+1},\hat{\g}_{q+8l+2},\hat{\g}_{q+8l+3},\hat{\g}_{q+8l+5}\}\}\setminus\{\hat{\g}_{q+2},\hat{\g}_{q+3}\}\big]^k\otimes \{F,\hat{\g}_{q+1,q+2},\hat{\g}_{q+2,q+3},\hat{\g_{q+1,q+3}}\}$, where we identify $\hat{\g_{q+1,q+2}}\leftrightarrow i$,$\hat{\g_{q+2,q+3}}\leftrightarrow j,\hat{\g_{q+1,q+3}}\leftrightarrow k$ & $2^{q+4l+2}$
\\ 
\hline 
\end{tabular} 
\end{center}
\caption{  $\mathbb{H}$-basis for $S_{p,q}$ where $q-p\equiv 2,3,4\mod 8$
}
\label{tab:table6}
\end{table}

With a firm understanding of the how the $\R$, $\C$ and  $\mathbb{H}$ spinor spaces look for $Cl_{p,q}$ we can now view our spinor representations for $Cl_{p,q}$ in the endomorphism algebra of these spinor spaces, or as matrices in the appropriate matrix algebra. 

\section{Spinor Representations for Clifford algebras}

\subsection{$Cl_{0,7}$}
We begin this section with the spinor representation for $Cl_{0,7}$, where we have the  spinor space isomorphism $S_{0,7}\oplus \hat{S_{0,7}}\cong \R^8\oplus \R^8$. We can use our table of generating sets and use the elements $\g_{124},\g_{235},\g_{346},\g_{457}$ to construct the primitive idempotent$F=\dfrac{1+\g_{124}}{2}\dfrac{1+\g_{235}}{2}\dfrac{1+\g_{346}}{2}\dfrac{1+\g_{457}}{2}$.  Giving us the left ideal $S_{0,7}$ that we will use to construct the spinor space $S_{0,7}\oplus \hat{S_{0,7}}$ to get our representations for $Cl_{0,7}$, where $\hat{S_{0,7}}=\{\phi_t(\g_{\alpha}F):\g_{\alpha}F\in S_{0,7}\}$ is the space that corresponds to the principal involution. As we can see by our previous section , or computationally , that we we have the $\R$ basis $F,\g_1F,\g_{2}F,\g_{3}F,\g_{12}F,\g_{13}F,\g_{23}F,\g_{123}F$ for the half spinor space $S_{0,7}$  giving us the  vector space isomorphism $S_{0,7}\cong \R^8$. We can get the odd (positive )spinor representation in $Mat(8,\R)$. (Since $Cl_{0,7}^+\cong End_{\R}(S_{0,7})$. Giving us the following :

 $$\g_1 \leftrightarrow \scriptsize{\left[\begin{array}{cccccccc} 0 & -1 & 0 & 0 & 0 & 0 & 0 & 0\\ 1& 0 & 0 & 0 & 0 & 0 & 0 & 0\\ 0 & 0 & 0 & 0 & -1 & 0 & 0 & 0\\ 0 & 0 & 0 & 0 & 0 & -1 & 0 & 0\\ 0 & 0 & 1 & 0 & 0 & 0 & 0 & 0\\0 & 0 & 0 & 1 & 0 & 0 & 0 & 0\\ 0 & 0 & 0 & 0 & 0 & 0 & 0 & -1\\ 0 & 0 & 0 & 0 & 0 & 0 & 1 & 0 \end{array}\right]},\g_2 \leftrightarrow \scriptsize{\left[\begin{array}{cccccccc} 0 & 0 & -1 & 0 & 0 & 0 & 0 & 0\\ 0 & 0 & 0 & 0 & 1 & 0 & 0 & 0\\ 1 & 0 & 0 & 0 & 0 & 0 & 0 & 0\\ 0 & 0 & 0 & 0 & 0 & 0 & -1 & 0\\ 0 & -1 & 0 & 0 & 0 & 0 & 0 & 0\\0 & 0 & 0 & 0 & 0 & 0 & 0 & 1\\ 0 & 0 & 0 & 1 & 0 & 0 & 0 & 0\\ 0 & 0 & 0 & 0 & 0 & -1 & 0 & 0 \end{array}\right]},$$
 
 $$\g_3 \leftrightarrow \scriptsize{\left[\begin{array}{cccccccc} 0 & 0 & 0 & -1 & 0 & 0 & 0 & 0\\ 0& 0 & 0 & 0 & 0 & 1 & 0 & 0\\ 0 & 0 & 0 & 0 & 0 & 0 & 1 & 0\\ 1 & 0 & 0 & 0 & 0 & 0 & 0 & 0\\ 0 & 0 & 0 & 0 & 0 & 0 & 0 & -1\\0 & -1 & 0 & 0 & 0 & 0 & 0 & 0\\ 0 & 0 & -1 & 0 & 0 & 0 & 0 & 0\\ 0 & 0 & 0 & 0 & 1 & 0 & 0 & 0 \end{array}\right]},\g_4 \leftrightarrow \scriptsize{\left[\begin{array}{cccccccc} 0 & 0 & 0 & 0 & 1 & 0 & 0 & 0\\ 0 & 0 & 1 & 0 & 0 & 0 & 0 & 0\\ 0 & -1 & 0 & 0 & 0 & 0 & 0 & 0\\ 0 & 0 & 0 & 0 & 0 & 0 & 0 & -1\\ -1& 0 & 0 & 0 & 0 & 0 & 0 & 0\\0 & 0 & 0 & 0 & 0 & 0 & -1 & 0\\ 0 & 0 & 0 & 0 & 0 & 1 & 0 & 0\\ 0 & 0 & 0 & 1 & 0 & 0 & 0 & 0 \end{array}\right]},$$
 
 $$\g_5 \leftrightarrow \scriptsize{\left[\begin{array}{cccccccc} 0 & 0 & 0 & 0 & 0 & 0 & 1 & 0\\ 0& 0 & 0 & 0 & 0 & 0 & 0 & -1\\ 0 & 0 & 0 & 1 & 0 & 0 & 0 & 0\\ 0 & 0 & -1 & 0 & 0 & 0 & 0 & 0\\ 0 & 0 & 0 & 0 & 0 & -1 & 0 & 0\\0 & 0 & 0 & 0 & 1 & 0 & 0 & 0\\ -1& 0 & 0 & 0 & 0 & 0 & 0 & 0\\ 0 & 1 & 0 & 0 & 0 & 0 & 0 & 0 \end{array}\right]},\g_6 \leftrightarrow \scriptsize{\left[\begin{array}{cccccccc} 0 & 0 & 0 & 0 & 0 & 0 & 0 & -1\\ 0 & 0 & 0 & 0 & 0 & 0 & -1 & 0\\ 0 & 0 & 0 & 0 & 0 & 1 & 0 & 0\\ 0 & 0 & 0 & 0 & -1 & 0 & 0 & 0\\ 0 & 0 & 0 & 1 & 0 & 0 & 0 & 0\\0 & 0 & -1 & 0 & 0 & 0 & 0 & 0\\ 0 & 1 & 0 & 0 & 0 & 0 & 0 & 0\\ 1 & 0 & 0 & 0 & 0 & 0 & 0 & 0 \end{array}\right]},$$
 
 $$\g_7 \leftrightarrow \scriptsize{\left[\begin{array}{cccccccc} 0 & 0 & 0 & 0 & 0 & -1 & 0 & 0\\ 0 & 0 & 0 & -1 & 0 & 0 & 0 & 0\\ 0 & 0 & 0 & 0 & 0 & 0 & 0 & -1\\ 0 & 1 & 0 & 0 & 0 & 0 & 0 & 0\\ 0 & 0 & 0 & 0 & 0 & 0 & -1 & 0\\1 & 0 & 0 & 0 & 0 & 0 & 0 & 0\\ 0 & 0 & 0 & 0 & 1 & 0 & 0 & 0\\ 0 & 0 & 1 & 0 & 0 & 0 & 0 & 0 \end{array}\right]},$$

We then can see that using $\hat{S}_{0,7}$ we can get the negative spinor representation, whose elements can also be seen in $Mat(8,\R)$, since $Cl_{0,7}^-\cong End_{\R}(\hat{S}_{0,7})$. Where the matrices of the basis elements would be nothing more than the negative of these matrices. Their direct sum of would give us our spinor representation  for $Cl_{0,7}$.

\subsection{ Spinor Representations for other Clifford algebras}
We conclude this section with two examples of spinor representions where the spinor space is a $\C$ and an $\mathbb{H}$-space. 
\subsubsection{The spinor representation for $Cl_{2,3}$}
 To find the spinor representation  for $Cl_{2,3}$ we use the generating set $\{\g_{13},\g_{24}\}$  the construction of the primitive idempotent  $F$, where we will use Table 3 to find the $\R$-basis for $S_{2,3}$. Thus as an $\R$-basis we have ;
 $$S_{2,3}= span_{\R}\{F,\g_1F,\g_2F,\g_5F,\g_{12}F,\g_{15}F,\g_{25}F,\g_{125}F\}\cong\R^8$$
 . Using Table 5 the $\R$-basis  can be expressed as a $\C$-basis in the following manner:
$$S_{2,3}=span_{\C}\{P^+,\g_1P^+,\g_2P^+,\g_{12}P^+\}\cong\C^4$$ , where $\g_{5}F$ is identified with the generator  $i$ in $\C$.

To find representations for $Cl_{2,3}$ in $End_{\C}(S_{2,3})$ we make the usual identification with elements in the Clifford algebra and the left multiplication endomorphism; resulting in the following identification with matrices in $Mat(4,\C)$.

$$\g_{1}\leftrightarrow \scriptsize{\left[\begin{array}{cccc} 0 & 1 & 0 & 0\\ 1& 0 & 0 & 0\\ 0 & 0 & 0 & 1 \\ 0 & 0 & 1 & 0 \end{array}\right]}\;\;
\g_{2}\leftrightarrow \scriptsize{\left[\begin{array}{cccc} 0 & 0 & 1 & 0\\ 0 & 0 & 0 & -1\\ 1 & 0 & 0 & 0 \\ 0 & -1 & 0 & 0 \end{array}\right]}$$

$$\g_{3}\leftrightarrow \scriptsize{\left[\begin{array}{cccc} 0 & -1 & 0 & 0\\ 1 & 0 & 0 & 0\\ 0& 0 & 0 & -1 \\ 0 & 0 & 1 & 0 \end{array}\right]}\;\;
\g_{4}\leftrightarrow \scriptsize{\left[\begin{array}{cccc} 0 & 0 & -1 & 0\\ 0& 0 & 0 & 1\\ 1 & 0 & 0 & 0\\ 0 & -1 & 0 & 0 \end{array}\right]},$$

$$\g_{5}\leftrightarrow \scriptsize{\left[\begin{array}{cccc} i & 0 & 0 & 0\\ 0 & -i & 0 & 0\\ 0& 0 & -i & 0 \\ 0 & 0 & 0 & i \end{array}\right]}.$$
\vspace{.1 in}

\subsubsection{The Spinor representation for $Cl_{4,0}$}

To find the spinor  representation for $Cl_{4,0}$   we use the primitive idempotent $F=\dfrac{1+\g_1}{2}$ so that $S_{4,0}$ is spanned ,as an $\R$-basis, by $span_{\R} \{F,\g_2F,\g_3F,\g_4F,\g_{23}F,\g_{24}F,\g_{34}F,\g_{234}F\}\cong \R^8$. In $S_{4,0}$ the multiplicative structure of generators $\g_{23}F,\g_{24}F,\g_{34}F$  behave like generators $i,j,k\in\mathbb{H}$. Thus we identify $\g_{23}F\leftrightarrow i,\g_{34}F\leftrightarrow j,\g_{34}F\leftrightarrow k$. Notice that the $\R$ -basis for $S_{4,0}$ can be expressed as $\{F,\g_{2}F\}\otimes \{F,\g_{23}F,\g_{34}F,\g_{24}F\}$. Using the identification with $\mathbb{H}$ we can view $S_{0,4}$ as an $\mathbb{H}$-space. To find representations for $Cl_{0,4}$ in $End_{\mathbb{H}}(S_{0,4})$ we make the usual identification with elements in the Clifford algebra and the left multiplication endomorphism; resulting in the following identification with matrices in $Mat(2,\mathbb{H})$:

$$\g_{1}\leftrightarrow \scriptsize{\left[\begin{array}{cccc} 1&0 \\ 0 & -1\end{array}\right]}\;\;
\g_{2}\leftrightarrow \scriptsize{\left[\begin{array}{cccc} 0 & 1 \\ 1 & 0\end{array}\right]}$$

$$\g_{3}\leftrightarrow \scriptsize{\left[\begin{array}{cccc} 0 & -i \\ i & 0\end{array}\right]}\;\;
\g_{4}\leftrightarrow \scriptsize{\left[\begin{array}{cccc} 0 & -k \\ k & 0\end{array}\right]},$$

\subsection{Other Spinor representations}

We will now provide some calculated spinor representations for Clifford algebras using the same method as the examples above.
\subsubsection{Examples of $\R$ -spinor representations}
\begin{itemize}
\item $Cl_{0,6}$ 
\subitem Generating set: $\{\g_{124},\g_{235},\g_{346}\}$ 
\subitem $S_{0,6}=span_{\R}\{F,\g_1F,\g_2F,\g_{3}F,\g_{4}F,\g_{5}F,\g_{6}F,\g_{13}F\}\cong \R^8,$

\subitem Spinor representation in $Mat(8,\R)$: 
$$\g_1 \leftrightarrow \scriptsize{\left[\begin{array}{cccccccc} 0 & -1 & 0 & 0 & 0 & 0 & 0 & 0\\ 1& 0 & 0 & 0 & 0 & 0 & 0 & 0\\ 0 & 0 & 0 & 0 & 1 & 0 & 0 & 0\\ 0 & 0 & 0 & 0 & 0 & 0 & 0 & -1\\ 0 & 0 & -1 & 0 & 0 & 0 & 0 & 0\\0 & 0 & 0 & 0 & 0 & 0 & 1 & 0\\ 0 & 0 & 0 & 0 & 0 & -1 & 0 & 0\\ 0 & 0 & 0 & 1 & 0 & 0 & 0 & 0 \end{array}\right]},\g_2 \leftrightarrow \scriptsize{\left[\begin{array}{cccccccc} 0 & 0 & -1 & 0 & 0 & 0 & 0 & 0\\ 0 & 0 & 0 & 0 & -1 & 0 & 0 & 0\\ 1 & 0 & 0 & 0 & 0 & 0 & 0 & 0\\ 0 & 0 & 0 & 0 & 0 & 1 & 0 & 0\\ 0 & 1 & 0 & 0 & 0 & 0 & 0 & 0\\0 & 0 & 0 & -1 & 0 & 0 & 0 & 0\\ 0 & 0 & 0 & 0 & 0 & 0 & 0 & -1\\ 0 & 0 & 0 & 0 & 0 & 0 & 1 & 0 \end{array}\right]},$$

$$\g_3 \leftrightarrow \scriptsize{\left[\begin{array}{cccccccc} 0 & 0 & 0 & -1 & 0 & 0 & 0 & 0\\ 0 & 0 & 0 & 0 & 0 & 0 & 0 & 1\\ 0 & 0 & 0 & 0 & 0 & -1 & 0& 0\\ 1 & 0 & 0 & 0 & 0 & 0 & 0 & 0\\ 0 & 0 & 0 & 0 & 0 & 0 & 1 & 0\\0 & 0 & 1 & 0 & 0 & 0 & 0 & 0\\ 0 & 0 & 0 & 0 & -1 & 0 & 0 & 0\\ 0 & -1 & 0 & 0 & 0 & 0 & 0 & 0 \end{array}\right]},\g_4 \leftrightarrow \scriptsize{\left[\begin{array}{cccccccc} 0 & 0 & 0 & 0 & -1 & 0 & 0 & 0\\ 0& 0 & 1 & 0 & 0 & 0 & 0 & 0\\ 0 & -1 & 0 & 0 & 0 & 0 & 0 & 0\\ 0 & 0 & 0 & 0 & 0 & 0 & -1 & 0\\ 1 & 0 & 0 & 0 & 0 & 0 & 0 & 0\\0 & 0 & 0 & 0 & 0 & 0 & 0 & -1\\ 0 & 0 & 0 & 1 & 0 & 0 & 0 & 0\\ 0 & 0 & 0 & 0 & 0 & 1 & 0 & 0 \end{array}\right]}$$

$$\g_5 \leftrightarrow \scriptsize{\left[\begin{array}{cccccccc} 0 & 0 & 0 & 0 & 0 & -1 & 0 & 0\\ 0 & 0 & 0 & 0 & 0 & 0 & -1 & 0\\ 0 & 0 & 0 & 1 & 0 & 0 & 0 & 0\\ 0 & 0 & -1 & 0 & 0 & 0 & 0 & 0\\ 0 & 0 & 0 & 0 & 0 & 0 & 0 & 1\\1 & 0 & 0 & 0 & 0 & 0 & 0 & 0\\ 0 & 1 & 0 & 0 & 0 & 0 & 0 & 0\\ 0 & 0 & 0 & 0 & -1 & 0 & 0 & 0 \end{array}\right]},\g_6 \leftrightarrow \scriptsize{\left[\begin{array}{cccccccc} 0 & 0 & 0 & 0 & 0 & 0 & -1 & 0\\ 0 & 0 & 0 & 0 & 0 & 1 & 0 & 0\\ 0 & 0 & 0 & 0 & 0 & 0 & 0 & 1\\ 0 & 0 & 0 & 0 & 1 & 0 & 0 & 0\\ 0 & 0 & 0 & -1 & 0 & 0 & 0 & 0\\0 & -1 & 0 & 0 & 0 & 0 & 0 & 0\\ 1 & 0 & 0 & 0 & 0 & 0 & 0 & 0\\ 0 & 0 & -1 & 0 & 0 & 0 & 0 & 0 \end{array}\right]}.$$

\item $Cl_{3,3}$ 
\subitem Generating set: $\{\g_{14},\g_{25},\g_{36}\}$ 
\subitem $S_{3,3}=span_{\R}\{F,\g_1F,\g_2F,\g_{3}F,\g_{12}F,\g_{13}F,\g_{23}F,\g_{123}F\}\cong\R^8$.

\subitem Spinor representation  in $Mat(8,\R)$: 
$$\g_1 \leftrightarrow \scriptsize{\left[\begin{array}{cccccccc} 0 & 1 & 0 & 0 & 0 & 0 & 0 & 0\\ 1& 0 & 0 & 0 & 0 & 0 & 0 & 0\\ 0 & 0 & 0 & 0 & 1 & 0 & 0 & 0\\ 0 & 0 & 0 & 0 & 0 & 1 & 0 & 0\\ 0 & 0 & 1 & 0 & 0 & 0 & 0 & 0\\0 & 0 & 0 & 1 & 0 & 0 & 0 & 0\\ 0 & 0 & 0 & 0 & 0 & 0 & 0 & 1\\ 0 & 0 & 0 & 0 & 0 & 0 & 1 & 0 \end{array}\right]},\g_2 \leftrightarrow \scriptsize{\left[\begin{array}{cccccccc} 0 & 0 & 1 & 0 & 0 & 0 & 0 & 0\\ 0 & 0 & 0 & 0 & -1 & 0 & 0 & 0\\ 1 & 0 & 0 & 0 & 0 & 0 & 0 & 0\\ 0 & 0 & 0 & 0 & 0 & 0 & -1 & 0\\ 0 & 1 & 0 & 0 & 0 & 0 & 0 & 0\\0 & 0 & 0 & 0 & 0 & 0 & 0 & -1\\ 0 & 0 & 0 & 1 & 0 & 0 & 0 & 0\\ 0 & 0 & 0 & 0 & 0 & -1 & 0 & 0 \end{array}\right]},$$

$$\g_3 \leftrightarrow \scriptsize{\left[\begin{array}{cccccccc} 0 & 0 & 0 & 1 & 0 & 0 & 0 & 0\\ 0 & 0 & 0 & 0 & 0 & -1 & 0 & 0\\ 0 & 0 & 0 & 0 & 0 & 0 & -1 & 0\\ 1 & 0 & 0 & 0 & 0 & 0 & 0 & 0\\ 0 & 0 & 0 & 0 & 0 & 0 & 0 & 1\\0 & -1 & 0 & 0 & 0 & 0 & 0 & 0\\ 0 & 0 & -1 & 0 & 0 & 0 & 0 & 0\\ 0 & 0 & 0 & 0 & 1 & 0 & 0 & 0 \end{array}\right]},\g_4 \leftrightarrow \scriptsize{\left[\begin{array}{cccccccc} 0 & -1 & 0 & 0 & 0 & 0 & 0 & 0\\ 1& 0 & 0 & 0 & 0 & 0 & 0 & 0\\ 0 & 0 & 0 & 0 & -1 & 0 & 0 & 0\\ 0 & 0 & 0 & 0 & 0 & -1 & 0 & 0\\ 0 & 0 & 1 & 0 & 0 & 0 & 0 & 0\\0 & 0 & 0 & 1 & 0 & 0 & 0 & 0\\ 0 & 0 & 0 & 0 & 0 & 0 & 0 & -1\\ 0 & 0 & 0 & 0 & 0 & 0 & 1 & 0 \end{array}\right]}$$

$$\g_5 \leftrightarrow \scriptsize{\left[\begin{array}{cccccccc} 0 & 0 & -1 & 0 & 0 & 0 & 0 & 0\\ 0 & 0 & 0 & 0 & 1 & 0 & 0 & 0\\ 1 & 0 & 0 & 0 & 0 & 0 & 0 & 0\\ 0 & 0 & 0 & 0 & 0 & 0 & -1 & 0\\ 0 & -1 & 0 & 0 & 0 & 0 & 0 & 0\\0 & 0 & 0 & 0 & 0 & 0 & 0 & 1\\ 0 & 0 & 0 & 1 & 0 & 0 & 0 & 0\\ 0 & 0 & 0 & 0 & 0 & -1 & 0 & 0 \end{array}\right]},\g_6 \leftrightarrow \scriptsize{\left[\begin{array}{cccccccc} 0 & 0 & 0 & -1 & 0 & 0 & 0 & 0\\ 0 & 0 & 0 & 0 & 0 & 1 & 0 & 0\\ 0 & 0 & 0 & 0 & 0 & 0 & 1 & 0\\ 1 & 0 & 0 & 0 & 0 & 0 & 0 & 0\\ 0 & 0 & 0 & 0 & 0 & 0 & 0 & -1\\0 & -1 & 0 & 0 & 0 & 0 & 0 & 0\\ 0 & 0 & -1 & 0 & 0 & 0 & 0 & 0\\ 0 & 0 & 0 & 0 & 1 & 0 & 0 & 0 \end{array}\right]}.$$

\item $Cl_{4,2}$ 
\subitem Generating set: $\{\g_{15},\g_{26},\g_{3}\}$ 
\subitem $S_{4,2}=span_{\R}\{F,\g_1F,\g_2F,\g_{4}F,\g_{12}F,\g_{14}F,\g_{24}F,\g_{124}F\}\cong\R^8$.

\subitem Spinor representation in $Mat(8,\R)$: 
$$\g_1 \leftrightarrow \scriptsize{\left[\begin{array}{cccccccc} 0 & 1 & 0 & 0 & 0 & 0 & 0 & 0\\ 1& 0 & 0 & 0 & 0 & 0 & 0 & 0\\ 0 & 0 & 0 & 0 & 1 & 0 & 0 & 0\\ 0 & 0 & 0 & 0 & 0 & 1 & 0 & 0\\ 0 & 0 & 1 & 0 & 0 & 0 & 0 & 0\\0 & 0 & 0 & 1 & 0 & 0 & 0 & 0\\ 0 & 0 & 0 & 0 & 0 & 0 & 0 & 1\\ 0 & 0 & 0 & 0 & 0 & 0 & 1 & 0 \end{array}\right]},\g_2 \leftrightarrow \scriptsize{\left[\begin{array}{cccccccc} 0 & 0 & 1 & 0 & 0 & 0 & 0 & 0\\ 0 & 0 & 0 & 0 & -1 & 0 & 0 & 0\\ 1 & 0 & 0 & 0 & 0 & 0 & 0 & 0\\ 0 & 0 & 0 & 0 & 0 & 0 & 1 & 0\\ 0 & -1 & 0 & 0 & 0 & 0 & 0 & 0\\0 & 0 & 0 & 0 & 0 & 0 & 0 & -1\\ 0 & 0 & 0 & 1 & 0 & 0 & 0 & 0\\ 0 & 0 & 0 & 0 & 0 & -1 & 0 & 0 \end{array}\right]},$$

$$\g_3 \leftrightarrow \scriptsize{\left[\begin{array}{cccccccc} 1 & 0 & 0 & 0 & 0 & 0 & 0 & 0\\ 0 & -1 & 0 & 0 & 0 & 0 & 0 & 0\\ 0 & 0 & -1 & 0 & 0 & 0 & 0 & 0\\ 0 & 0 & 0 & -1 & 0 & 0 & 0 & 0\\ 0 & 0 & 0 & 0 & 1 & 0 & 0 & 0\\0 & 0 & 0 & 0 & 0 & 1 & 0 & 0\\ 0 & 0 & 0 & 0 & 0 & 0 & 1 & 0\\ 0 & 0 & 0 & 0 & 0 & 0 & 0 & -1 \end{array}\right]},\g_4 \leftrightarrow \scriptsize{\left[\begin{array}{cccccccc} 0 & 0 & 0 & 1 & 0 & 0 & 0 & 0\\ 0& 0 & 0 & 0 & 0 & -1 & 0 & 0\\ 0 & 0 & 0 & 0 & 0 & 0 & -1 & 0\\ 1 & 0 & 0 & 0 & 0 & 0 & 0 & 0\\ 0 & 0 & 0 & 0 & 0 & 0 & 0 & 1\\0 & -1 & 0 & 0 & 0 & 0 & 0 & 0\\ 0 & 0 & -1 & 0 & 0 & 0 & 0 & 0\\ 0 & 0 & 0 & 0 & 1 & 0 & 0 & 0 \end{array}\right]}$$

$$\g_5 \leftrightarrow \scriptsize{\left[\begin{array}{cccccccc} 0 & -1 & 0 & 0 & 0 & 0 & 0 & 0\\ 1 & 0 & 0 & 0 & 0 & 0 & 0 & 0\\ 0& 0 & 0 & 0 & -1 & 0 & 0 & 0\\ 0 & 0 & 0 & 0 & 0 & -1 & 0 & 0\\ 0 & 0 & 1 & 0 & 0 & 0 & 0 & 0\\0 & 0 & 0 & 1 & 0 & 0 & 0 & 0\\ 0 & 0 & 0 & 0 & 0 & 0 & 0 & -1\\ 0 & 0 & 0 & 0 & 0 & 0 & 1 & 0 \end{array}\right]},\g_6 \leftrightarrow \scriptsize{\left[\begin{array}{cccccccc} 0 & 0 & -1 & 0 & 0 & 0 & 0 & 0\\ 0 & 0 & 0 & 0 & 1 & 0 & 0 & 0\\ 1 & 0 & 0 & 0 & 0 & 0 & 0 & 0\\ 0 & 0 & 0 & 0 & 0 & 0 & -1 & 0\\ 0 & -1 & 0 & 0 & 0 & 0 & 0 & 0\\0 & 0 & 0 & 0 & 0 & 0 & 0 & 1\\ 0 & 0 & 0& 1 & 0 & 0 & 0 & 0\\ 0 & 0 & 0 & 0 & 0 & -1 & 0 & 0 \end{array}\right]}.$$

\item $Cl_{8,0}$ 
\subitem Generating set: $\{\g_{1234},\g_{1256},\g_{1278},\g_{1357}\}$ 
\subitem $S_{8,0}=span_{\R}\{1,\g_1,\g_2,\g_3,\g_5,\g_{12},\g_{13},\g_{15},\g_{23},\g_{25},\g_{35},\g_{123},\g_{125},\g_{135},\g_{235},\g_{1235}\}\otimes\{F\}\cong\R^{16}$.

\subitem Spinor represention in $Mat(16,\R)$, 
$$\g_1 \leftrightarrow \scriptsize{\left[\begin{array}{cccccccccccccccc} | & | &| & | & | & | & | & |&| & | &| & | & | & | & | & | \\e_{2}^T & e_{1}^T & e_{6}^T & e_{7}^T& e_{8}^T & e_{3}^T & e_{4}^T & e_{5}^T &e_{12}^T & e_{13}^T & e_{14}^T & e_{9}^T& e_{10}^T & e_{11}^T & e_{16}^T & e_{15}^T\\| & | &| & | & | & | & | & |&| & | &| & | & | & | & | & |\end{array}\right]} ,$$

$$\g_2 \leftrightarrow \scriptsize{\left[\begin{array}{cccccccccccccccc} | & | &| & | & | & | & | & |&| & | &| & | & | & | & | & | \\e_{3}^T & -e_{6}^T & e_{1}^T & e_{9}^T& e_{10}^T & -e_{2}^T & -e_{12}^T & -e_{13}^T &e_{4}^T & e_{5}^T & e_{15}^T & -e_{7}^T& -e_{8}^T & -e_{16}^T & e_{11}^T & -e_{14}^T\\| & | &| & | & | & | & | & |&| & | &| & | & | & | & | & |\end{array}\right]} ,$$

$$\g_3 \leftrightarrow \scriptsize{\left[\begin{array}{cccccccccccccccc} | & | &| & | & | & | & | & |&| & | &| & | & | & | & | & | \\e_{4}^T & -e_{7}^T & -e_{9}^T & e_{1}^T& e_{11}^T & e_{12}^T & -e_{2}^T & -e_{14}^T &-e_{3}^T & -e_{15}^T & e_{5}^T & e_{6}^T& e_{16}^T & -e_{8}^T & -e_{10}^T & e_{13}^T\\| & | &| & | & | & | & | & |&| & | &| & | & | & | & | & |\end{array}\right]} ,$$

$$\g_4 \leftrightarrow \scriptsize{\left[\begin{array}{cccccccccccccccc} | & | &| & | & | & | & | & |&| & | &| & | & | & | & | & | \\-e_{12}^T & e_{9}^T & -e_{7}^T & e_{6}^T& -e_{16}^T & e_{4}^T & -e_{3}^T & e_{15}^T &e_{2}^T & -e_{14}^T & e_{13}^T & -e_{1}^T& e_{11}^T & -e_{10}^T & e_{8}^T & -e_{5}^T\\| & | &| & | & | & | & | & |&| & | &| & | & | & | & | & |\end{array}\right]} ,$$

$$\g_5 \leftrightarrow \scriptsize{\left[\begin{array}{cccccccccccccccc} | & | &| & | & | & | & | & |&| & | &| & | & | & | & | & | \\e_{5}^T & -e_{8}^T & -e_{10}^T & -e_{11}^T& e_{1}^T & e_{13}^T & e_{14}^T & -e_{2}^T &e_{15}^T & -e_{3}^T & -e_{4}^T & -e_{16}^T& e_{6}^T & e_{7}^T & e_{9}^T & -e_{12}^T\\| & | &| & | & | & | & | & |&| & | &| & | & | & | & | & |\end{array}\right]} ,$$

$$\g_6 \leftrightarrow \scriptsize{\left[\begin{array}{cccccccccccccccc} | & | &| & | & | & | & | & |&| & | &| & | & | & | & | & | \\-e_{13}^T & e_{10}^T & -e_{8}^T & e_{16}^T& e_{6}^T & e_{5}^T & -e_{15}^T & -e_{3}^T & e_{14}^T & e_{2}^T & -e_{12}^T & -e_{11}^T& -e_{1}^T & e_{9}^T & -e_{7}^T & e_{4}^T\\| & | &| & | & | & | & | & |&| & | &| & | & | & | & | & |\end{array}\right]} ,$$

$$\g_7 \leftrightarrow \scriptsize{\left[\begin{array}{cccccccccccccccc} | & | &| & | & | & | & | & |&| & | &| & | & | & | & | & | \\  -e_{14}^T & e_{11}^T & -e_{16}^T& -e_{8}^T & e_{7}^T & e_{15}^T & e_{5}^T & -e_{4}^T &- e_{13}^T & e_{12}^T & e_{2}^T& e_{10}^T & -e_{9}^T & -e_{1}^T & e_{6}^T & -e_{3}^T\\| & | &| & | & | & | & | & |&| & | &| & | & | & | & | & |\end{array}\right]} ,$$

$$\g_8 \leftrightarrow \scriptsize{\left[\begin{array}{cccccccccccccccc} | & | &| & | & | & | & | & |&| & | &| & | & | & | & | & | \\-e_{15}^T & e_{16}^T & e_{11}^T & -e_{10}^T& e_{9}^T & -e_{14}^T & e_{13}^T & -e_{12}^T &e_{5}^T & -e_{4}^T & e_{3}^T & -e_{8}^T& e_{7}^T & -e_{6}^T & -e_{1}^T & e_{2}^T\\| & | &| & | & | & | & | & |&| & | &| & | & | & | & | & |\end{array}\right]} ,$$

\item $Cl_{5,3}$ 
\subitem Generating set: $\{\g_{16},\g_{27},\g_{38},\g_{4}\}$ 
\subitem $S_{5,3}=span_{\R}\{1,\g_1,\g_2,\g_3,\g_5,\g_{12},\g_{13},\g_{15},\g_{23},\g_{25},\g_{35},\g_{123},\g_{125},\g_{135},\g_{235},\g_{1235}\}\otimes\{F\}\cong \R^{16}$.

\subitem Spinor representation of in $Mat(16,\R)$
$$\g_1 \leftrightarrow \scriptsize{\left[\begin{array}{cccccccccccccccc} | & | &| & | & | & | & | & |&| & | &| & | & | & | & | & | \\e_{2}^T & e_{1}^T & e_{6}^T & e_{7}^T& e_{8}^T & e_{3}^T & e_{4}^T & e_{5}^T &e_{12}^T & e_{13}^T & e_{14}^T & e_{9}^T& e_{10}^T & e_{11}^T & e_{16}^T & e_{15}^T\\| & | &| & | & | & | & | & |&| & | &| & | & | & | & | & |\end{array}\right]} ,$$

$$\g_2 \leftrightarrow \scriptsize{\left[\begin{array}{cccccccccccccccc} | & | &| & | & | & | & | & |&| & | &| & | & | & | & | & | \\e_{3}^T & -e_{6}^T & e_{1}^T & e_{9}^T& e_{10}^T & -e_{2}^T & -e_{12}^T & -e_{13}^T &e_{4}^T & e_{5}^T & e_{15}^T & -e_{7}^T&  -e_{8}^T & -e_{16}^T & e_{11}^T & -e_{14}^T\\| & | &| & | & | & | & | & |&| & | &| & | & | & | & | & |\end{array}\right]} ,$$

$$\g_3\leftrightarrow \scriptsize{\left[\begin{array}{cccccccccccccccc} | & | &| & | & | & | & | & |&| & | &| & | & | & | & | & | \\e_{4}^T & -e_{7}^T & -e_{9}^T & e_{1}^T& e_{11}^T & e_{12}^T & -e_{2}^T & -e_{14}^T &-e_{3}^T & -e_{15}^T & e_{5}^T & e_{6}^T& e_{16}^T & -e_{8}^T & -e_{10}^T & e_{13}^T\\| & | &| & | & | & | & | & |&| & | &| & | & | & | & | & |\end{array}\right]} ,$$

$$\g_4 \leftrightarrow \scriptsize{\left[\begin{array}{cccccccccccccccc} | & | &| & | & | & | & | & |&| & | &| & | & | & | & | & | \\e_{1}^T & -e_{2}^T & -e_{3}^T & -e_{4}^T& e_{5}^T & e_{6}^T & e_{7}^T & -e_{8}^T &e_{9}^T & -e_{10}^T & -e_{11}^T & -e_{12}^T& e_{13}^T & e_{14}^T & e_{15}^T & -e_{16}^T\\| & | &| & | & | & | & | & |&| & | &| & | & | & | & | & |\end{array}\right]} ,$$

$$\g_5 \leftrightarrow \scriptsize{\left[\begin{array}{cccccccccccccccc} | & | &| & | & | & | & | & |&| & | &| & | & | & | & | & | \\e_{5}^T & -e_{8}^T & -e_{10}^T & -e_{11}^T& e_{1}^T & e_{13}^T & e_{14}^T & -e_{2}^T &e_{15}^T & -e_{3}^T & -e_{4}^T & -e_{16}^T& e_{6}^T & e_{7}^T & e_{9}^T & -e_{12}^T\\| & | &| & | & | & | & | & |&| & | &| & | & | & | & | & |\end{array}\right]} ,$$

$$\g_6 \leftrightarrow \scriptsize{\left[\begin{array}{cccccccccccccccc} | & | &| & | & | & | & | & |&| & | &| & | & | & | & | & | \\e_{2}^T & -e_{1}^T & e_{6}^T & e_{7}^T& e_{8}^T & -e_{3}^T & -e_{4}^T & -e_{5}^T &e_{12}^T & e_{13}^T & e_{14}^T & -e_{9}^T& -e_{10}^T & -e_{11}^T & e_{16}^T & -e_{15}^T\\| & | &| & | & | & | & | & |&| & | &| & | & | & | & | & |\end{array}\right]} ,$$

$$\g_7 \leftrightarrow \scriptsize{\left[\begin{array}{cccccccccccccccc} | & | &| & | & | & | & | & |&| & | &| & | & | & | & | & | \\e_{3}^T & -e_{6}^T & -e_{1}^T & e_{9}^T& e_{10}^T & e_{2}^T & -e_{12}^T & -e_{13}^T &-e_{4}^T & -e_{5}^T & e_{15}^T & e_{7}^T& e_{8}^T & -e_{16}^T & -e_{11}^T & e_{14}^T\\| & | &| & | & | & | & | & |&| & | &| & | & | & | & | & |\end{array}\right]} ,$$

$$\g_8 \leftrightarrow \scriptsize{\left[\begin{array}{cccccccccccccccc} | & | &| & | & | & | & | & |&| & | &| & | & | & | & | & | \\e_{4}^T & -e_{7}^T & -e_{9}^T & -e_{1}^T& e_{11}^T & e_{12}^T & e_{2}^T & -e_{14}^T &e_{3}^T & -e_{15}^T & -e_{5}^T & -e_{6}^T& e_{16}^T & e_{8}^T & e_{10}^T & -e_{13}^T\\| & | &| & | & | & | & | & |&| & | &| & | & | & | & | & |\end{array}\right]} ,$$

\item $Cl_{4,4}$ 
\subitem Generating set: $\{\g_{15},\g_{26},\g_{37},\g_{48}\}$ 
\subitem $S_{4,4}=span_{\R}\{1,\g_1,\g_2,\g_3,\g_4,\g_{12},\g_{13},\g_{14},\g_{23},\g_{24},\g_{34},\g_{123},\g_{124},\g_{134},\g_{234},\g_{1234}\}\otimes\{F\}\cong\R^{16}$.

\subitem Spinor representation in $Mat(16,\R)$:
$$\g_1 \leftrightarrow \scriptsize{\left[\begin{array}{cccccccccccccccc} | & | &| & | & | & | & | & |&| & | &| & | & | & | & | & | \\e_{2}^T & e_{1}^T & e_{6}^T & e_{7}^T& e_{8}^T & e_{3}^T & e_{4}^T & e_{5}^T &e_{12}^T & e_{13}^T & e_{14}^T & e_{9}^T& e_{10}^T & e_{11}^T & e_{16}^T & e_{15}^T\\| & | &| & | & | & | & | & |&| & | &| & | & | & | & | & |\end{array}\right]} ,$$
$$\g_2 \leftrightarrow \scriptsize{\left[\begin{array}{cccccccccccccccc} | & | &| & | & | & | & | & |&| & | &| & | & | & | & | & | \\e_{3}^T & -e_{6}^T & e_{1}^T & e_{9}^T& e_{10}^T & -e_{2}^T & -e_{12}^T & -e_{13}^T & e_{4}^T & e_{5}^T & e_{15}^T & -e_{7}^T& -e_{8}^T & -e_{16}^T & e_{11}^T & -e_{14}^T\\| & | &| & | & | & | & | & |&| & | &| & | & | & | & | & |\end{array}\right]} ,$$

$$\g_3\leftrightarrow \scriptsize{\left[\begin{array}{cccccccccccccccc} | & | &| & | & | & | & | & |&| & | &| & | & | & | & | & | \\e_{4}^T & -e_{7}^T & -e_{9}^T & e_{1}^T& e_{11}^T & e_{12}^T & -e_{2}^T & -e_{14}^T &-e_{3}^T & -e_{15}^T & e_{5}^T & e_{6}^T& e_{16}^T & -e_{8}^T & -e_{10}^T & e_{13}^T\\| & | &| & | & | & | & | & |&| & | &| & | & | & | & | & |\end{array}\right]} ,$$

$$\g_4 \leftrightarrow \scriptsize{\left[\begin{array}{cccccccccccccccc} | & | &| & | & | & | & | & |&| & | &| & | & | & | & | & | \\e_{5}^T & -e_{8}^T & -e_{10}^T & -e_{11}^T& e_{1}^T & e_{13}^T & e_{14}^T & -e_{2}^T &e_{15}^T & -e_{3}^T & -e_{4}^T & -e_{16}^T& e_{6}^T & e_{7}^T & e_{9}^T & -e_{12}^T\\| & | &| & | & | & | & | & |&| & | &| & | & | & | & | & |\end{array}\right]} ,$$

$$\g_5 \leftrightarrow \scriptsize{\left[\begin{array}{cccccccccccccccc} | & | &| & | & | & | & | & |&| & | &| & | & | & | & | & | \\e_{2}^T & -e_{1}^T & e_{6}^T & e_{7}^T& e_{8}^T & -e_{3}^T & -e_{4}^T & -e_{5}^T &e_{12}^T & e_{13}^T & e_{14}^T & -e_{9}^T& -e_{10}^T & -e_{11}^T & e_{16}^T & -e_{15}^T\\| & | &| & | & | & | & | & |&| & | &| & | & | & | & | & |\end{array}\right]} ,$$

$$\g_6 \leftrightarrow \scriptsize{\left[\begin{array}{cccccccccccccccc} | & | &| & | & | & | & | & |&| & | &| & | & | & | & | & | \\e_{3}^T & -e_{6}^T & -e_{1}^T & e_{9}^T& e_{10}^T & e_{2}^T & -e_{12}^T & -e_{13}^T &-e_{4}^T & -e_{5}^T & e_{15}^T & e_{7}^T& e_{8}^T & -e_{16}^T & -e_{11}^T & e_{14}^T\\| & | &| & | & | & | & | & |&| & | &| & | & | & | & | & |\end{array}\right]} ,$$

$$\g_7\leftrightarrow \scriptsize{\left[\begin{array}{cccccccccccccccc} | & | &| & | & | & | & | & |&| & | &| & | & | & | & | & | \\e_{4}^T & -e_{7}^T & -e_{9}^T & -e_{1}^T& e_{11}^T & e_{12}^T & e_{2}^T & -e_{14}^T & e_{3}^T & -e_{15}^T & -e_{5}^T & -e_{6}^T& e_{16}^T & e_{8}^T & e_{10}^T & -e_{13}^T\\| & | &| & | & | & | & | & |&| & | &| & | & | & | & | & |\end{array}\right]} ,$$

$$\g_8 \leftrightarrow \scriptsize{\left[\begin{array}{cccccccccccccccc} | & | &| & | & | & | & | & |&| & | &| & | & | & | & | & | \\e_{5}^T & -e_{8}^T & -e_{10}^T & -e_{11}^T& -e_{1}^T & e_{13}^T & e_{14}^T & e_{2}^T &e_{15}^T & e_{3}^T & e_{4}^T & -e_{16}^T& -e_{6}^T & -e_{7}^T & -e_{9}^T & e_{12}^T\\| & | &| & | & | & | & | & |&| & | &| & | & | & | & | & |\end{array}\right]} ,$$

\item $Cl_{1,7}$ 
\subitem Generating set: $\{\g_{12},\g_{346},\g_{457},\g_{568}\}$ 
\subitem $S_{1,7}=span_{\R}\{1,\g_1,\g_3,\g_4,\g_5,\g_{13},\g_{14},\g_{15},\g_{34},\g_{35},\g_{45},\g_{134},\g_{135},\g_{145},\g_{345},\g_{1345}\}\otimes\{F\}\cong \R^{16}$.

\subitem Spinor representation in $Mat(16,\R)$:

$$\g_1 \leftrightarrow \scriptsize{\left[\begin{array}{cccccccccccccccc} | & | &| & | & | & | & | & |&| & | &| & | & | & | & | & | \\e_{2}^T & e_{1}^T & e_{6}^T & e_{7}^T& e_{8}^T & e_{3}^T & e_{4}^T & e_{5}^T &e_{12}^T & e_{13}^T & e_{14}^T & e_{9}^T& e_{10}^T & e_{11}^T & e_{16}^T & e_{15}^T\\| & | &| & | & | & | & | & |&| & | &| & | & | & | & | & |\end{array}\right]} ,$$

$$\g_2 \leftrightarrow \scriptsize{\left[\begin{array}{cccccccccccccccc} | & | &| & | & | & | & | & |&| & | &| & | & | & | & | & | \\e_{2}^T & -e_{1}^T & e_{6}^T & e_{7}^T& e_{8}^T & -e_{3}^T & -e_{4}^T & -e_{5}^T &e_{12}^T & e_{13}^T & e_{14}^T & -e_{9}^T& -e_{10}^T & -e_{11}^T & e_{16}^T & -e_{15}^T\\| & | &| & | & | & | & | & |&| & | &| & | & | & | & | & |\end{array}\right]} ,$$
$$\g_3 \leftrightarrow \scriptsize{\left[\begin{array}{cccccccccccccccc} | & | &| & | & | & | & | & |&| & | &| & | & | & | & | & | \\e_{3}^T & -e_{6}^T & -e_{1}^T & e_{9}^T& e_{10}^T & e_{2}^T & -e_{12}^T & -e_{13}^T &-e_{4}^T & -e_{5}^T & e_{15}^T & e_{7}^T& e_{8}^T & -e_{16}^T & -e_{11}^T & e_{14}^T\\| & | &| & | & | & | & | & |&| & | &| & | & | & | & | & |\end{array}\right]} ,$$

$$\g_4\leftrightarrow \scriptsize{\left[\begin{array}{cccccccccccccccc} | & | &| & | & | & | & | & |&| & | &| & | & | & | & | & | \\e_{4}^T & -e_{7}^T & -e_{9}^T & -e_{1}^T& e_{11}^T & e_{12}^T & e_{2}^T & -e_{14}^T & e_{3}^T & -e_{15}^T & -e_{5}^T & -e_{6}^T& e_{16}^T & e_{8}^T & e_{10}^T & -e_{13}^T\\| & | &| & | & | & | & | & |&| & | &| & | & | & | & | & |\end{array}\right]} ,$$

$$\g_5 \leftrightarrow \scriptsize{\left[\begin{array}{cccccccccccccccc} | & | &| & | & | & | & | & |&| & | &| & | & | & | & | & | \\e_{5}^T & -e_{8}^T & -e_{10}^T & -e_{11}^T& -e_{1}^T & e_{13}^T & e_{14}^T & e_{2}^T &e_{15}^T & e_{3}^T & e_{4}^T & -e_{16}^T& -e_{6}^T & -e_{7}^T & -e_{9}^T & e_{12}^T\\| & | &| & | & | & | & | & |&| & | &| & | & | & | & | & |\end{array}\right]} ,$$

$$\g_6 \leftrightarrow \scriptsize{\left[\begin{array}{cccccccccccccccc} | & | &| & | & | & | & | & |&| & | &| & | & | & | & | & | \\-e_{9}^T & e_{12}^T & -e_{4}^T & e_{3}^T& e_{15}^T & e_{7}^T & -e_{6}^T & -e_{16}^T &e_{1}^T & e_{11}^T & -e_{10}^T & -e_{2}^T& -e_{14}^T & e_{13}^T & -e_{5}^T & e_{8}^T\\| & | &| & | & | & | & | & |&| & | &| & | & | & | & | & |\end{array}\right]} ,$$

$$\g_7 \leftrightarrow \scriptsize{\left[\begin{array}{cccccccccccccccc} | & | &| & | & | & | & | & |&| & | &| & | & | & | & | & | \\-e_{11}^T & e_{14}^T & e_{15}^T & -e_{5}^T& e_{4}^T & -e_{16}^T & e_{8}^T & -e_{7}^T &e_{10}^T & -e_{9}^T & e_{1}^T & -e_{13}^T& e_{12}^T & -e_{2}^T & -e_{3}^T & e_{6}^T\\| & | &| & | & | & | & | & |&| & | &| & | & | & | & | & |\end{array}\right]} ,$$

$$\g_8 \leftrightarrow \scriptsize{\left[\begin{array}{cccccccccccccccc} | & | &| & | & | & | & | & |&| & | &| & | & | & | & | & | \\e_{15}^T & -e_{16}^T & e_{11}^T & -e_{10}^T& e_{9}^T & -e_{14}^T & e_{13}^T & -e_{12}^T &-e_{5}^T & e_{4}^T & -e_{3}^T & e_{8}^T& -e_{7}^T & e_{6}^T & -e_{1}^T & e_{2}^T\\| & | &| & | & | & | & | & |&| & | &| & | & | & | & | & |\end{array}\right]} ,$$

\item $Cl_{0,8}$ 
\subitem Generating set: $\{\g_{1234},\g_{1256},\g_{1278},\g_{1357}\}$ 
\subitem $S_{0,8}=span_{\R}\{1,\g_1,\g_2,\g_3,\g_5,\g_{12},\g_{13},\g_{15},\g_{23},\g_{25},\g_{35},\g_{123},\g_{125},\g_{135},\g_{235},\g_{1235}\}\otimes\{F\}\cong \R^{16}$.

\subitem Spinor representation in $Mat(16,\R)$
$$\g_1 \leftrightarrow \scriptsize{\left[\begin{array}{cccccccccccccccc} | & | &| & | & | & | & | & |&| & | &| & | & | & | & | & | \\e_{2}^T & -e_{1}^T & e_{6}^T & e_{7}^T& e_{8}^T & -e_{3}^T & -e_{4}^T & -e_{5}^T &e_{12}^T & e_{13}^T & e_{14}^T & -e_{9}^T& -e_{10}^T & -e_{11}^T & e_{16}^T & -e_{15}^T\\| & | &| & | & | & | & | & |&| & | &| & | & | & | & | & |\end{array}\right]} ,$$

$$\g_2 \leftrightarrow \scriptsize{\left[\begin{array}{cccccccccccccccc} | & | &| & | & | & | & | & |&| & | &| & | & | & | & | & | \\e_{3}^T & -e_{6}^T & -e_{1}^T & e_{9}^T& e_{10}^T & e_{2}^T & -e_{12}^T & -e_{13}^T &-e_{4}^T & -e_{5}^T & e_{15}^T & e_{7}^T& e_{8}^T & -e_{16}^T & -e_{11}^T & e_{14}^T\\| & | &| & | & | & | & | & |&| & | &| & | & | & | & | & |\end{array}\right]} ,$$

$$\g_3\leftrightarrow \scriptsize{\left[\begin{array}{cccccccccccccccc} | & | &| & | & | & | & | & |&| & | &| & | & | & | & | & | \\e_{4}^T & -e_{7}^T & -e_{9}^T & -e_{1}^T& e_{11}^T & e_{12}^T & e_{2}^T & -e_{14}^T & e_{3}^T & -e_{15}^T & -e_{5}^T & -e_{6}^T& e_{16}^T & e_{8}^T & e_{10}^T & -e_{13}^T\\| & | &| & | & | & | & | & |&| & | &| & | & | & | & | & |\end{array}\right]} ,$$

$$\g_4\leftrightarrow \scriptsize{\left[\begin{array}{cccccccccccccccc} | & | &| & | & | & | & | & |&| & | &| & | & | & | & | & | \\e_{12}^T & e_{9}^T & -e_{7}^T & e_{6}^T& e_{16}^T & -e_{4}^T & e_{3}^T & e_{15}^T & -e_{2}^T & -e_{14}^T & e_{13}^T & -e_{1}^T& -e_{11}^T & e_{10}^T & -e_{8}^T & -e_{5}^T\\| & | &| & | & | & | & | & |&| & | &| & | & | & | & | & |\end{array}\right]} ,$$

$$\g_5 \leftrightarrow \scriptsize{\left[\begin{array}{cccccccccccccccc} | & | &| & | & | & | & | & |&| & | &| & | & | & | & | & | \\e_{5}^T & -e_{8}^T & -e_{10}^T & -e_{11}^T& -e_{1}^T & e_{13}^T & e_{14}^T & e_{2}^T &e_{15}^T & e_{3}^T & e_{4}^T & -e_{16}^T& -e_{6}^T & -e_{7}^T & -e_{9}^T & e_{12}^T\\| & | &| & | & | & | & | & |&| & | &| & | & | & | & | & |\end{array}\right]} ,$$

$$\g_6\leftrightarrow \scriptsize{\left[\begin{array}{cccccccccccccccc} | & | &| & | & | & | & | & |&| & | &| & | & | & | & | & | \\e_{13}^T & e_{10}^T & -e_{8}^T & -e_{16}^T& e_{6}^T & -e_{5}^T & -e_{15}^T & e_{3}^T & e_{14}^T & -e_{2}^T & -e_{12}^T & e_{11}^T& -e_{1}^T & -e_{9}^T & e_{7}^T & e_{4}^T\\| & | &| & | & | & | & | & |&| & | &| & | & | & | & | & |\end{array}\right]} ,$$

$$\g_7\leftrightarrow \scriptsize{\left[\begin{array}{cccccccccccccccc} | & | &| & | & | & | & | & |&| & | &| & | & | & | & | & | \\e_{14}^T & e_{11}^T & e_{16}^T & -e_{8}^T& e_{7}^T & e_{15}^T & -e_{5}^T & e_{4}^T & -e_{13}^T & e_{12}^T & -e_{2}^T & -e_{10}^T& e_{9}^T & -e_{1}^T & -e_{6}^T & -e_{3}^T\\| & | &| & | & | & | & | & |&| & | &| & | & | & | & | & |\end{array}\right]} ,$$

$$\g_8\leftrightarrow \scriptsize{\left[\begin{array}{cccccccccccccccc} | & | &| & | & | & | & | & |&| & | &| & | & | & | & | & | \\ e_{15}^T&-e_{16}^T & e_{11}^T & -e_{10}^T & e_{9}^T& -e_{14}^T & e_{13}^T & -e_{12}^T & -e_{5}^T & e_{4}^T & -e_{3}^T & e_{8}^T & -e_{7}^T& e_{6}^T & -e_{1}^T & e_{2}^T \\| & | &| & | & | & | & | & |&| & | &| & | & | & | & | & |\end{array}\right]} ,$$

\end{itemize}

\subsubsection{$\C$-spinor representations}

\begin{itemize}
\item $Cl_{1,2}$ 
\subitem Generating set $\{\g_{12}\}$ 
\subitem $S_{1,2}=span_{\R}\{F,\g_1F,\g_3F,\g_{13}F\}\cong \R^4$
\subitem  $S_{1,2}=span_{\C}\{F,\g_1F\}\cong \C^2$. Identifying $\g_{3}F$ with $i$.
\subitem Spinor representation in $Mat(2,\C)$ 

$$\g_{1}\leftrightarrow \scriptsize{\left[\begin{array}{cccc} 0 & 1 \\ 1& 0\end{array}\right]}\;\;
\g_{2}\leftrightarrow \scriptsize{\left[\begin{array}{cccc} 0& -1 \\ 1 & 0 \end{array}\right]}\;\;\g_{3}\leftrightarrow \scriptsize{\left[\begin{array}{cccc} i& 0 \\ 0 & i \end{array}\right]}$$

\item $Cl_{0,5}$ 
\subitem Generating set $\{\g_{1234},\g_{125}\}$ 
\subitem $S_{0,5}=span_{\R}\{F,\g_1F,\g_3F,\g_{5}F,\g_{13}F,\g_{15}F,\g_{35}F,\g_{135}F\}\cong \R^8$,
\subitem  $S_{0,5}=span_{\C}\{F,\g_1F,\g_3F,\g_{13}F\}\cong \C^4$. Identifying $i$ with  $\g_{5}F$, .
\subitem Spinor representation in $Mat(4,\C)$ 

$$\g_{1}\leftrightarrow \scriptsize{\left[\begin{array}{cccc} 0 & -1 & 0 & 0\\ 1& 0 & 0 & 0\\ 0 & 0 & 0 & -1 \\ 0 & 0 & 1 & 0 \end{array}\right]}\;\;
\g_{2}\leftrightarrow \scriptsize{\left[\begin{array}{cccc} 0& i & 0 & 0\\ i & 0 & 0 & 0\\ 0 & 0 & 0 & i \\ 0 & 0 & i & 0 \end{array}\right]}$$

$$\g_{3}\leftrightarrow \scriptsize{\left[\begin{array}{cccc} 0 & 0 & -1 & 0\\ 0 & 0 & 0 & 1\\ 1& 0 & 0 & 0 \\ 0 & -1 & 0 & 0 \end{array}\right]}\;\;
\g_{4}\leftrightarrow \scriptsize{\left[\begin{array}{cccc} 0 & 0 & -i & 0\\ 0& 0 & 0 & i\\ -i & 0 & 0 & 0\\ 0 & i & 0 & 0 \end{array}\right]},$$

$$\g_{5}\leftrightarrow \scriptsize{\left[\begin{array}{cccc} i & 0 & 0 & 0\\ 0 & -i & 0 & 0\\ 0& 0 & -i & 0\\ 0 & 0 & 0 & i\end{array}\right]}.$$

\item $Cl_{4,1}$ 
\subitem Generating set $\{\g_{15},\g_{2}\}$ 
\subitem $S_{4,1}=span_{\R}\{F,\g_1F,\g_3F,\g_{4}F,\g_{13}F,\g_{14}F,\g_{34}F,\g_{134}F\}\cong\R^8$
\subitem  $S_{4,1}=span_{\C}\{F,\g_1F,\g_4F,\g_{14}F\}\cong \C^4$ ,identifying $i$ with $\g_{13}F$.
\subitem Spinor representation in $Mat(4,\C)$

$$\g_{1}\leftrightarrow \scriptsize{\left[\begin{array}{cccc} 0 & 1 & 0 & 0\\ 1& 0 & 0 & 0\\ 0 & 0 & 0 & 1 \\ 0 & 0 & 1 & 0 \end{array}\right]}\;\;
\g_{2}\leftrightarrow \scriptsize{\left[\begin{array}{cccc} 1& 0 & 0 & 0\\ 0 & -1 & 0 & 0\\ 0 & 0 & -1 & 0 \\ 0 & 0 & 0 & 1 \end{array}\right]}$$

$$\g_{3}\leftrightarrow \scriptsize{\left[\begin{array}{cccc} 0 & -i & 0 & 0\\ i & 0 & 0 & 0\\ 0& 0 & 0 & -i \\ 0 & 0 & i & 0 \end{array}\right]}\;\;
\g_{4}\leftrightarrow \scriptsize{\left[\begin{array}{cccc} 0 & 0 & 1 & 0\\ 0& 0 & 0 & -1\\ 1 & 0 & 0 & 0\\ 0 & -1 & 0 & 0 \end{array}\right]},$$

$$\g_{5}\leftrightarrow \scriptsize{\left[\begin{array}{cccc} 0 & -1 & 0 & 0\\ 1 & 0 & 0 & 0\\ 0& 0 & 0 & -1 \\ 0 & 0 & 1 &0 \end{array}\right]}.$$

\item $Cl_{7,0}$ 
\subitem Generating set $\{\g_{1234},\g_{1256},\g_{127}\}$ 
\subitem $S_{7,0}=span_{\R}\{1,\g_1,\g_2,\g_3,\g_{5},\g_{12},\g_{13},\g_{15},\g_{23},\g_{25},\g_{35},\g_{123},\g_{125},\g_{135},\g_{235},\g_{1235}\}\otimes\{F\}\cong\R^{16}$.
\subitem  $S_{7,0}=span_{\C}\{P^+,\g_1P^+,\g_2P^+,\g_{3}P^+,\g_{12}P^+,\g_{13}P^+,\g_{23}P^+,\g_{123}P^+\}\cong\C^8$. Identifying $i$ with $\g_{15}F$.
\subitem Spinor representation in $Mat(8,\C)$.

$$\g_1 \leftrightarrow \scriptsize{\left[\begin{array}{cccccccc} 0 & 1 & 0 & 0 & 0 & 0 & 0 & 0\\ 1& 0 & 0 & 0 & 0 & 0 & 0 & 0\\ 0 & 0 & 0 & 0 & 1 & 0 & 0 & 0\\ 0 & 0 & 0 & 0 & 0 & 1 & 0 & 0\\ 0 & 0 & 1 & 0 & 0 & 0 & 0 & 0\\0 & 0 & 0 & 1 & 0 & 0 & 0 & 0\\ 0 & 0 & 0 & 0 & 0 & 0 & 0 & 1\\ 0 & 0 & 0 & 0 & 0 & 0 & 1 & 0 \end{array}\right]},\g_2 \leftrightarrow \scriptsize{\left[\begin{array}{cccccccc} 0 & 0 & 1 & 0 & 0 & 0 & 0 & 0\\ 0 & 0 & 0 & 0 & -1 & 0 & 0 & 0\\ 1 & 0 & 0 & 0 & 0 & 0 & 0 & 0\\ 0 & 0 & 0 & 0 & 0 & 0 & 1 & 0\\ 0 & -1 & 0 & 0 & 0 & 0 & 0 & 0\\0 & 0 & 0 & 0 & 0 & 0 & 0 & -1\\ 0 & 0 & 0 & 1 & 0 & 0 & 0 & 0\\ 0 & 0 & 0 & 0 & 0 & -1 & 0 & 0 \end{array}\right]},$$

$$\g_3 \leftrightarrow \scriptsize{\left[\begin{array}{cccccccc} 0 & 0 & 0 & 1 & 0 & 0 & 0 & 0\\ 0 & 0 & 0 & 0 & 0 & -1 & 0 & 0\\ 0 & 0 & 0 & 0 & 0 & 0 & -1 & 0\\ 1 & 0 & 0 & 0 & 0 & 0 & 0 & 0\\ 0 & 0 & 0 & 0 & 0 & 0 & 0 & 1\\0 & -1 & 0 & 0 & 0 & 0 & 0 & 0\\ 0 & 0 & -1 & 0 & 0 & 0 & 0 & 0\\ 0 & 0 & 0 & 0 & 1 & 0 & 0 & 0 \end{array}\right]},\g_4 \leftrightarrow \scriptsize{\left[\begin{array}{cccccccc} 0 & 0 & 0 & 0 & 0 & 0 & 0 & -1\\ 0& 0 & 0 & 0 & 0 & 0 & 1 & 0\\ 0 & 0 & 0 & 0 & 0 & -1 & 0 & 0\\ 0 & 0 & 0 & 0 & 1 & 0 & 0 & 0\\ 0 & 0 & 0 & 1 & 0 & 0 & 0 & 1\\0 & 0 & -1 & 0 & 0 & 0 & 0 & 0\\ 0 & 1 & 0 & 0 & 0 & 0 & 0 & 0\\ -1 & 0 & 0 & 0 & 0 & 0 & 0 & 0 \end{array}\right]}$$

$$\g_5 \leftrightarrow \scriptsize{\left[\begin{array}{cccccccc} 0 & -i& 0 & 0 & 0 & 0 & 0 & 0\\ i& 0 & 0 & 0 & 0 & 0 & 0 & 0\\ 0 & 0 & 0 & 0 & -i & 0 & 0 & 0\\ 0 & 0 & 0 & 0 & 0 & -i & 0 & 0\\ 0 & 0 & i & 0 & 0 & 0 & 0 & 0\\0 & 0 & 0 & i & 0 & 0 & 0 & 0\\ 0 & 0 & 0 & 0 & 0 & 0 & 0 & -i\\ 0 & 0 & 0 & 0 & 0 & 0 & i & 0 \end{array}\right]},\g_6 \leftrightarrow \scriptsize{\left[\begin{array}{cccccccc} 0 & 0 & -i & 0 & 0 & 0 & 0 & 0\\ 0 & 0 & 0 & 0 & i & 0 & 0 & 0\\ i & 0 & 0 & 0 & 0 & 0 & 0 & 0\\ 0 & 0 & 0 & 0 & 0 & 0 & -i & 0\\ 0 & -i & 0 & 0 & 0 & 0 & 0 & 0\\0 & 0 & 0 & 0 & 0 & 0 & 0 & i\\ 0 & 0 & 0 & i & 0 & 0 & 0 & 0\\ 0 & 0 & 0 & 0 & 0 & -i & 0 & 0 \end{array}\right]},$$
$$\g_{7} \leftrightarrow \scriptsize{\left[\begin{array}{cccccccc} 1 & 0 & 0 & 0 & 0 & 0 & 0 & 0\\ 0& -1 & 0 & 0 & 0 & 0 & 0 & 0\\ 0 & 0 & -1 & 0 & 0 & 0 & 0 & 0\\ 0 & 0 & 0 & -1 & 0 & 0 & 0 & 0\\ 0 & 0 & 0 & 0 & 1 & 0 & 0 & 0\\0 & 0 & 0 & 0 & 0 & 1 & 0 & 0\\ 0 & 0 & 0 & 0 & 0 & 0 & 1 & 0\\ 0 & 0 & 0 & 0 & 0 & 0 & 0 & -1 \end{array}\right]},$$

\item $Cl_{5,2}$ 
\subitem Generating set $\{\g_{16},\g_{27},\g_{3}\}$ 
\subitem $S_{7,0}=span_{\R}\{1,\g_1,\g_2,\g_4,\g_{5},\g_{12},\g_{14},\g_{15},\g_{24},\g_{25},\g_{45},\g_{124},\g_{125},\g_{145},\g_{245},\g_{1245}\}\otimes\{F\}\cong \R^{16}$.
\subitem  $S_{5,2}=span_{\C}\{F,\g_1F,\g_2F,\g_{4}F,\g_{12}F,\g_{14}F,\g_{24}F,\g_{124}F\}\cong \C^8$. Identifying $i$ with $\g_{15}F$.
\subitem Spinor representation in $Mat(8,\C)$

$$\g_1 \leftrightarrow \scriptsize{\left[\begin{array}{cccccccc} 0 & 1 & 0 & 0 & 0 & 0 & 0 & 0\\ 1& 0 & 0 & 0 & 0 & 0 & 0 & 0\\ 0 & 0 & 0 & 0 & 1 & 0 & 0 & 0\\ 0 & 0 & 0 & 0 & 0 & 1 & 0 & 0\\ 0 & 0 & 1 & 0 & 0 & 0 & 0 & 0\\0 & 0 & 0 & 1 & 0 & 0 & 0 & 0\\ 0 & 0 & 0 & 0 & 0 & 0 & 0 & 1\\ 0 & 0 & 0 & 0 & 0 & 0 & 1 & 0 \end{array}\right]},\g_2 \leftrightarrow \scriptsize{\left[\begin{array}{cccccccc} 0 & 0 & 1 & 0 & 0 & 0 & 0 & 0\\ 0 & 0 & 0 & 0 & -1 & 0 & 0 & 0\\ 1 & 0 & 0 & 0 & 0 & 0 & 0 & 0\\ 0 & 0 & 0 & 0 & 0 & 0 & 1 & 0\\ 0 & -1 & 0 & 0 & 0 & 0 & 0 & 0\\0 & 0 & 0 & 0 & 0 & 0 & 0 & -1\\ 0 & 0 & 0 & 1 & 0 & 0 & 0 & 0\\ 0 & 0 & 0 & 0 & 0 & -1 & 0 & 0 \end{array}\right]},$$

$$\g_3 \leftrightarrow \scriptsize{\left[\begin{array}{cccccccc} 1 & 0 & 0 & 0 & 0 & 0 & 0 & 0\\ 0& -1 & 0 & 0 & 0 & 0 & 0 & 0\\ 0 & 0 & -1 & 0 & 0 & 0 & 0 & 0\\ 0 & 0 & 0 & 1 & 0 & 0 & 0 & 0\\ 0 & 0 & 0 & 0 & 1 & 0 & 0 & 0\\0 & 0 & 0 & 0 & 0 & -1 & 0 & 0\\ 0 & 0 & 0 & 0 & 0 & 0 & -1 & 0\\ 0 & 0 & 0 & 0 & 0 & 0 & 0 & 1 \end{array}\right]},\g_4 \leftrightarrow \scriptsize{\left[\begin{array}{cccccccc} 0 & 0 & 0 & 1 & 0 & 0 & 0 & 0\\ 0& 0 & 0 & 0 & 0 & -1 & 0 & 0\\ 0 & 0 & 0 & 0 & 0 & 0 & -1 & 0\\ 1 & 0 & 0 & 0 & 0 & 0 & 0 & 0\\ 0 & 0 & 0 & 0 & 0 & 0 & 0 & 1\\0 & -1 & 0 & 0 & 0 & 0 & 0 & 0\\ 0 & 0 & -1 & 0 & 0 & 0 & 0 & 0\\ 0 & 0 & 0 & 0 & 1 & 0 & 0 & 0 \end{array}\right]}$$

$$\g_5 \leftrightarrow \scriptsize{\left[\begin{array}{cccccccc} 0 & -i& 0 & 0 & 0 & 0 & 0 & 0\\ i& 0 & 0 & 0 & 0 & 0 & 0 & 0\\ 0 & 0 & 0 & 0 & -i & 0 & 0 & 0\\ 0 & 0 & 0 & 0 & 0 & -i & 0 & 0\\ 0 & 0 & i & 0 & 0 & 0 & 0 & 0\\0 & 0 & 0 & i & 0 & 0 & 0 & 0\\ 0 & 0 & 0 & 0 & 0 & 0 & 0 & -i\\ 0 & 0 & 0 & 0 & 0 & 0 & i & 0 \end{array}\right]},\g_6 \leftrightarrow \scriptsize{\left[\begin{array}{cccccccc} 0 & -1 & 0 & 0 & 0 & 0 & 0 & 0\\ 1 & 0 & 0 & 0 & 0 & 0 & 0 & 0\\ 0 & 0 & 0 & 0 & -1 & 0 & 0 & 0\\ 0 & 0 & 0 & 0 & 0 & -1 & 0 & 0\\ 0 & 0 & 1 & 0 & 0 & 0 & 0 & 0\\0 & 0 & 0 & 1 & 0 & 0 & 0 & 0\\ 0 & 0 & 0 & 0 & 0 & 0 & 0 & -1\\ 0 & 0 & 0 & 0 & 0 & 0 & 1 & 0 \end{array}\right]},$$
$$\g_{7} \leftrightarrow \scriptsize{\left[\begin{array}{cccccccc} 0 & 0 & -1 & 0 & 0 & 0 & 0 & 0\\ 0& 0 & 0 & 0 & 1 & 0 & 0 & 0\\ 1 & 0 & 0 & 0 & 0 & 0 & 0 & 0\\ 0 & 0 & 0 & 0 & 0 & 0 & -1 & 0\\ 0 & -1 & 0 & 0 & 0 & 0 & 0 & 0\\0 & 0 & 0 & 0 & 0 & 0 & 0 & 1\\ 0 & 0 & 0 & 1 & 0 & 0 & 0 & 0\\ 0 & 0 & 0 & 0 & 0 & -1 & 0 & 0 \end{array}\right]},$$

\item $Cl_{3,4}$ 
\subitem Generating set $\{\g_{14},\g_{25},\g_{3}\}$ 
\subitem $S_{3,4}=span_{\R}\{1,\g_1,\g_2,\g_3,\g_{7},\g_{12},\g_{13},\g_{17},\g_{23},\g_{27},\g_{37},\g_{123},\g_{127},\g_{137},\g_{237},\g_{1237}\}\otimes\{F\}\cong \R^{16}$.
\subitem  $S_{3,4}=span_{\C}\{F,\g_1F,\g_2F,\g_{3}F,\g_{12}F,\g_{13}F,\g_{23}F,\g_{123}F\}$. Identifying $i$ with  $\g_{7}F$.
\subitem Spinor representation in $Mat(8,\C)$.

$$\g_1 \leftrightarrow \scriptsize{\left[\begin{array}{cccccccc} 0 & 1 & 0 & 0 & 0 & 0 & 0 & 0\\ 1& 0 & 0 & 0 & 0 & 0 & 0 & 0\\ 0 & 0 & 0 & 0 & 1 & 0 & 0 & 0\\ 0 & 0 & 0 & 0 & 0 & 1 & 0 & 0\\ 0 & 0 & 1 & 0 & 0 & 0 & 0 & 0\\0 & 0 & 0 & 1 & 0 & 0 & 0 & 0\\ 0 & 0 & 0 & 0 & 0 & 0 & 0 & 1\\ 0 & 0 & 0 & 0 & 0 & 0 & 1 & 0 \end{array}\right]},\g_2 \leftrightarrow \scriptsize{\left[\begin{array}{cccccccc} 0 & 0 & 1 & 0 & 0 & 0 & 0 & 0\\ 0 & 0 & 0 & 0 & -1 & 0 & 0 & 0\\ 1 & 0 & 0 & 0 & 0 & 0 & 0 & 0\\ 0 & 0 & 0 & 0 & 0 & 0 & 1 & 0\\ 0 & -1 & 0 & 0 & 0 & 0 & 0 & 0\\0 & 0 & 0 & 0 & 0 & 0 & 0 & -1\\ 0 & 0 & 0 & 1 & 0 & 0 & 0 & 0\\ 0 & 0 & 0 & 0 & 0 & -1 & 0 & 0 \end{array}\right]},$$

$$\g_3 \leftrightarrow \scriptsize{\left[\begin{array}{cccccccc} 0 & 0 & 0 & 1 & 0 & 0 & 0 & 0\\ 0& 0 & 0 & 0 & 0 & -1 & 0 & 0\\ 0 & 0 & 0 & 0 & 0 & 0 & -1 & 0\\ 1 & 0 & 0 & 0 & 0 & 0 & 0 & 0\\ 0 & 0 & 0 & 0 & 0 & 0 & 0 & 1\\0 & -1 & 0 & 0 & 0 & 0 & 0 & 0\\ 0 & 0 & -1 & 0 & 0 & 0 & 0 & 0\\ 0 & 0 & 0 & 0 & 1 & 0 & 0 & 0 \end{array}\right]},\g_4\leftrightarrow \scriptsize{\left[\begin{array}{cccccccc} 0 & -1 & 0 & 0 & 0 & 0 & 0 & 0\\ 1& 0 & 0 & 0 & 0 & 0 & 0 & 0\\ 0 & 0 & 0 & 0 & -1 & 0 & 0 & 0\\ 0 & 0 & 0 & 0 & 0 & -1 & 0 & 0\\ 0 & 0 & 1 & 0 & 0 & 0 & 0 & 0\\0 & 0 & 0 & 1 & 0 & 0 & 0 & 0\\ 0 & 0 & 0 & 0 & 0 & 0 & 0 & -1\\ 0 & 0 & 0 & 0 & 0 & 0 & 1 & 0 \end{array}\right]}.$$

$$\g_5 \leftrightarrow \scriptsize{\left[\begin{array}{cccccccc} 0 & 0 & -1 & 0 & 0 & 0 & 0 & 0\\ 0 & 0 & 0 & 0 & 1 & 0 & 0 & 0\\ 1 & 0 & 0 & 0 & 0 & 0 & 0 & 0\\ 0 & 0 & 0 & 0 & 0 & 0 & -1 & 0\\ 0 & -1 & 0 & 0 & 0 & 0 & 0 & 0\\0 & 0 & 0 & 0 & 0 & 0 & 0 & 1\\ 0 & 0 & 0 & 1 & 0 & 0 & 0 & 0\\ 0 & 0 & 0 & 0 & 0 & -1 & 0 & 0 \end{array}\right]},\g_6\leftrightarrow \scriptsize{\left[\begin{array}{cccccccc} 0 & 0 & 0 & -1 & 0 & 0 & 0 & 0\\ 0& 0 & 0 & 0 & 0 & 1 & 0 & 0\\ 0 & 0 & 0 & 0 & 0 & 0 & 1 & 0\\ 1 & 0 & 0 & 0 & 0 & 0 & 0 & 0\\ 0 & 0 & 0 & 0 & 0 & 0 & 0 & -1\\0 & -1 & 0 & 0 & 0 & 0 & 0 & 0\\ 0 & 0 & -1 & 0 & 0 & 0 & 0 & 0\\ 0 & 0 & 0 & 0 & 1 & 0 & 0 & 0 \end{array}\right]} ,$$
$$\g_{7} \leftrightarrow \scriptsize{\left[\begin{array}{cccccccc} i & 0 & 0 & 0 & 0 & 0 & 0 & 0\\ 0& -i & 0 & 0 & 0 & 0 & 0 & 0\\ 0 & 0 & -i & 0 & 0 & 0 & 0 & 0\\ 0 & 0 & 0 & -i & 0 & 0 & 0 & 0\\ 0 & 0 & 0 & 0 & i & 0 & 0 & 0\\0 & 0 & 0 & 0 & 0 & i & 0 & 0\\ 0 & 0 & 0 & 0 & 0 & 0 & i & 0\\ 0 & 0 & 0 & 0 & 0 & 0 & 0 & -i \end{array}\right]},$$

\item $Cl_{1,6}$ 
\subitem Generating set $\{\g_{12},\g_{3456},\g_{347}\}$ 
\subitem $S_{1,6}=span_{\R}\{1,\g_1,\g_3,\g_4,\g_{5},\g_{13},\g_{14},\g_{15},\g_{34},\g_{35},\g_{45},\g_{134},\g_{135},\g_{145},\g_{345},\g_{1345}\}\otimes\{F\}\cong \R^{16}$.
\subitem  $S_{1,6}=span_{\C}\{F,\g_1F,\g_3F,\g_{5}F,\g_{13}F,\g_{15}F,\g_{35}F,\g_{135}F\}\cong \C^8$. Identifying $i$ with $\g_{34}F$.
\subitem Spinor representation in $Mat(8,\C)$:

$$\g_1 \leftrightarrow \scriptsize{\left[\begin{array}{cccccccc} 0 & 1 & 0 & 0 & 0 & 0 & 0 & 0\\ 1& 0 & 0 & 0 & 0 & 0 & 0 & 0\\ 0 & 0 & 0 & 0 & 1 & 0 & 0 & 0\\ 0 & 0 & 0 & 0 & 0 & 1 & 0 & 0\\ 0 & 0 & 1 & 0 & 0 & 0 & 0 & 0\\0 & 0 & 0 & 1 & 0 & 0 & 0 & 0\\ 0 & 0 & 0 & 0 & 0 & 0 & 0 & 1\\ 0 & 0 & 0 & 0 & 0 & 0 & 1 & 0 \end{array}\right]},\g_2 \leftrightarrow ,\scriptsize{\left[\begin{array}{cccccccc} 0 & -1 & 0 & 0 & 0 & 0 & 0 & 0\\ 1& 0 & 0 & 0 & 0 & 0 & 0 & 0\\ 0 & 0 & 0 & 0 & -1 & 0 & 0 & 0\\ 0 & 0 & 0 & 0 & 0 & -1 & 0 & 0\\ 0 & 0 & 1 & 0 & 0 & 0 & 0 & 0\\0 & 0 & 0 & 1 & 0 & 0 & 0 & 0\\ 0 & 0 & 0 & 0 & 0 & 0 & 0 & -1\\ 0 & 0 & 0 & 0 & 0 & 0 & 1 & 0 \end{array}\right]}$$

$$\g_3 \leftrightarrow \scriptsize{\left[\begin{array}{cccccccc} 0 & 0 & -1 & 0 & 0 & 0 & 0 & 0\\ 0 & 0 & 0 & 0 & 1 & 0 & 0 & 0\\ 1 & 0 & 0 & 0 & 0 & 0 & 0 & 0\\ 0 & 0 & 0 & 0 & 0 & 0 & -1 & 0\\ 0 & -1 & 0 & 0 & 0 & 0 & 0 & 0\\0 & 0 & 0 & 0 & 0 & 0 & 0 & 1\\ 0 & 0 & 0 & 1 & 0 & 0 & 0 & 0\\ 0 & 0 & 0 & 0 & 0 & -1 & 0 & 0 \end{array}\right]},\g_4\leftrightarrow\scriptsize{\left[\begin{array}{cccccccc} 0 & 0 & -i & 0 & 0 & 0 & 0 & 0\\ 0 & 0 & 0 & 0 & i & 0 & 0 & 0\\ -i & 0 & 0 & 0 & 0 & 0 & 0 & 0\\ 0 & 0 & 0 & 0 & 0 & 0 & -i & 0\\ 0 & i & 0 & 0 & 0 & 0 & 0 & 0\\0 & 0 & 0 & 0 & 0 & 0 & 0 & i\\ 0 & 0 & 0 & -i & 0 & 0 & 0 & 0\\ 0 & 0 & 0 & 0 & 0 & i & 0 & 0 \end{array}\right]} .$$

$$\g_5 \leftrightarrow \scriptsize{\left[\begin{array}{cccccccc} 0 & 0 & 0 & -1& 0 & 0 & 0 & 0\\ 0& 0 & 0 & 0 & 0 & 1 & 0 & 0\\ 0 & 0 & 0 & 0 & 0 & 0 & 1 & 0\\ 1 & 0 & 0 & 0 & 0 & 0 & 0 & 0\\ 0 & 0 & 0 & 0 & 0 & 0 & 0 & -1\\0 & -1 & 0 & 0 & 0 & 0 & 0 & 0\\ 0 & 0 & -1 & 0 & 0 & 0 & 0 & 0\\ 0 & 0 & 0 & 0 & 1 & 0 & 0 & 0 \end{array}\right]},\g_6\leftrightarrow \scriptsize{\left[\begin{array}{cccccccc} 0 & 0 & 0 & i& 0 & 0 & 0 & 0\\ 0& 0 & 0 & 0 & 0 & -i & 0 & 0\\ 0 & 0 & 0 & 0 & 0 & 0 & -i & 0\\ i & 0 & 0 & 0 & 0 & 0 & 0 & 0\\ 0 & 0 & 0 & 0 & 0 & 0 & 0 & i\\0 & -i & 0 & 0 & 0 & 0 & 0 & 0\\ 0 & 0 & -i & 0 & 0 & 0 & 0 & 0\\ 0 & 0 & 0 & 0 & i & 0 & 0 & 0 \end{array}\right]} ,$$
$$\g_{7} \leftrightarrow\scriptsize{\left[\begin{array}{cccccccc} -i & 0 & 0 & 0 & 0 & 0 & 0 & 0\\ 0& i & 0 & 0 & 0 & 0 & 0 & 0\\ 0 & 0 & i & 0 & 0 & 0 & 0 & 0\\ 0 & 0 & 0 & i & 0 & 0 & 0 & 0\\ 0 & 0 & 0 & 0 & -i & 0 & 0 & 0\\0 & 0 & 0 & 0 & 0 & -i & 0 & 0\\ 0 & 0 & 0 & 0 & 0 & 0 & -i & 0\\ 0 & 0 & 0 & 0 & 0 & 0 & 0 & i \end{array}\right]},$$

\end{itemize}

\subsubsection{$\mathbb{H}$-spinor representations}
\begin{itemize}
\item $Cl_{0,4}$ 
\subitem Generating set: $\{\g_{1234}\}$ 
\subitem $S_{0,4}=span_{\R}\{F,\g_1F,\g_2F,\g_{3}F,\g_4F,\g_{12}F,\g_{13}F,\g_{14}F\}\cong \R^8$.
\subitem  $S_{0,4}=span_{\mathbb{H}}\{F,\g_1F\}\cong\mathbb{H}^2$.
\subitem Identifying; $i\leftrightarrow\g_{12}F$, $j\leftrightarrow\g_{13}F$, $k\leftrightarrow\g_{14}F$

\subitem Spinor representation in $Mat(2,\mathbb{H})$: 

$$\g_{1}\leftrightarrow \scriptsize{\left[\begin{array}{cccc} 0 & -1 \\ 1& 0\end{array}\right]}\;\;
\g_{2}\leftrightarrow \scriptsize{\left[\begin{array}{cccc} 0& i \\ i & 0 \end{array}\right]}\;\;\g_{3}\leftrightarrow \scriptsize{\left[\begin{array}{cccc} 0 & j \\ j& 0\end{array}\right]}\g_{4}\leftrightarrow \scriptsize{\left[\begin{array}{cccc} 0& k \\ k & 0 \end{array}\right]}$$

\item $Cl_{1,3}$ 
\subitem Generating set: $\{\g_{12}\}$ 
\subitem $S_{1,3}=span_{\R}\{F,\g_1F,\g_{3}F,\g_4F,\g_{13}F,\g_{14}F,\g_{34}F,\g_{134}F\}\cong \R^8$.
\subitem  $S_{1,3}=span_{\mathbb{H}}\{P^+,\g_1P^+\}\cong \mathbb{H}^2$.
\subitem Identifying; $i\leftrightarrow\g_{3}F$, $j\leftrightarrow\g_{4}F$, $k\leftrightarrow\g_{34}F$

\subitem Spinor representation in $Mat(2,\mathbb{H})$: 

$$\g_{1}\leftrightarrow \scriptsize{\left[\begin{array}{cccc} 0 & 1 \\ 1& 0\end{array}\right]}\;\;
\g_{2}\leftrightarrow \scriptsize{\left[\begin{array}{cccc} 0& -1 \\ 1 & 0 \end{array}\right]}\;\;\g_{3}\leftrightarrow \scriptsize{\left[\begin{array}{cccc} i & 0 \\ 0& -i\end{array}\right]}\g_{4}\leftrightarrow \scriptsize{\left[\begin{array}{cccc} j& 0 \\ 0 & -j \end{array}\right]}$$

\item $Cl_{6,0}$ 
\subitem Generating set: $\{\g_{1234},\g_5\}$ 
\subitem $S_{6,0}=span_{\R}\{1,\g_1,\g_2,\g_3,\g_4,\g_6,\g_{12},\g_{13},\g_{14},\g_{16},\g_{26},\g_{36},\g_{46},\g_{126},\g_{136},\g_{146}\}\otimes\{F\}\cong \R^{16}$.
\subitem  $S_{6,0}=span_{\mathbb{H}}\{F,\g_1F,\g_6F,\g_{16}F\}\cong\mathbb{H}^4$.
\subitem Identifying; $i\leftrightarrow\g_{12}F$, $j\leftrightarrow\g_{13}F$, $k\leftrightarrow\g_{14}F$

\subitem Spinor representation in $Mat(4,\mathbb{H})$:

$$\g_{1}\leftrightarrow \scriptsize{\left[\begin{array}{cccc} 0 & 1 & 0 & 0\\ 1& 0 & 0 & 0\\ 0 & 0 & 0 & 1 \\ 0 & 0 & 1 & 0 \end{array}\right]}\;\;
\g_{2}\leftrightarrow \scriptsize{\left[\begin{array}{cccc} 0 & -i & 0 & 0\\ i & 0 & 0 & 0\\ 0 & 0 & 0 & -i \\ 0 & 0 & i & 0 \end{array}\right]}$$

$$\g_{3}\leftrightarrow \scriptsize{\left[\begin{array}{cccc} 0 & -j & 0 & 0\\ j& 0 & 0 & 0\\ 0 & 0 & 0 & -j \\ 0 & 0 & j & 0 \end{array}\right]}\;\;
\g_{4}\leftrightarrow \scriptsize{\left[\begin{array}{cccc} 0 & -k & 0 & 0\\ k & 0 & 0 & 0\\ 0 & 0 & 0 & -k \\ 0 & 0 & k & 0 \end{array}\right]}$$

$$\g_{5}\leftrightarrow \scriptsize{\left[\begin{array}{cccc} 1 & 0 & 0 & 0\\ 0 & -1 & 0 & 0\\ 0& 0 & -1 & 0 \\ 0 & 0 &0 & 1 \end{array}\right]}\;\;
\g_{6}\leftrightarrow \scriptsize{\left[\begin{array}{cccc} 0 & 0 & 1 & 0\\ 0& 0 & 0 & -1\\ 1 & 0 & 0 & 0\\ 0 & -1 & 0 & 0 \end{array}\right]},$$

\item $Cl_{5,1}$ 
\subitem Generating set: $\{\g_{16},\g_2\}$ 
\subitem $S_{5,1}=span_{\R}\{1,\g_1,\g_3,\g_4,\g_5,\g_{13},\g_{14},\g_{15},\g_{34},\g_{35},\g_{45},\g_{134},\g_{135},\g_{145},\g_{345},\g_{1345}\}\otimes\{F\}\cong \R^{16}$.
\subitem  $S_{5,1}=span_{\mathbb{H}}\{F,\g_1F,\g_5F,\g_{15}F\}\cong\mathbb{H}^4$.
\subitem Identifying; $i\leftrightarrow\g_{14}F$, $j\leftrightarrow\g_{13}F$, $k\leftrightarrow\g_{34}F$

\subitem Spinor representation in $Mat(4,\mathbb{H})$.

$$\g_{1}\leftrightarrow \scriptsize{\left[\begin{array}{cccc} 0 & 1 & 0 & 0\\ 1& 0 & 0 & 0\\ 0 & 0 & 0 & 1 \\ 0 & 0 & 1 & 0 \end{array}\right]}\;\;
\g_{2}\leftrightarrow \scriptsize{\left[\begin{array}{cccc} 1 & 0 & 0 & 0\\ 0 & -1 & 0 & 0\\ 0 & 0 & -1 & 0 \\ 0 & 0 & 0 & 1 \end{array}\right]}$$

$$\g_{3}\leftrightarrow \scriptsize{\left[\begin{array}{cccc} 0 & -j & 0 & 0\\ j& 0 & 0 & 0\\ 0 & 0 & 0 & -j \\ 0 & 0 & j & 0 \end{array}\right]}\;\;
\g_{4}\leftrightarrow \scriptsize{\left[\begin{array}{cccc} 0 & -i & 0 & 0\\ i & 0 & 0 & 0\\ 0 & 0 & 0 & -i \\ 0 & 0 & i & 0 \end{array}\right]}$$

$$\g_{5}\leftrightarrow \scriptsize{\left[\begin{array}{cccc} 0 & 0 & 1 & 0\\ 0 & 0 & 0 & -1\\ 1& 0 & 0 & 0 \\ 0 & -1 &0 & 0\end{array}\right]}\;\;
\g_{6}\leftrightarrow \scriptsize{\left[\begin{array}{cccc} 0 & -1 & 0 & 0\\ 1& 0 & 0 & 0\\ 0 & 0 & 0 & -1\\ 0 & 0 & 1 & 0 \end{array}\right]},$$

\item $Cl_{2,4}$ 
\subitem Generating set: $\{\g_{13},\g_{24}\}$ 
\subitem $S_{2,4}=span_{\R}\{1,\g_1,\g_2,\g_5,\g_6,\g_{12},\g_{15},\g_{16},\g_{25},\g_{26},\g_{56},\g_{125},\g_{126},\g_{156},\g_{256},\g_{1256}\}\otimes\{F\}\cong \R^{16}$.
\subitem  $S_{2,4}=span_{\mathbb{H}}\{F,\g_1F,\g_2F,\g_{12}F\}\cong\mathbb{H}^4$.
\subitem Identifying; $i\leftrightarrow\g_{5}F$, $j\leftrightarrow\g_{6}F$, $k\leftrightarrow\g_{56}F$

\subitem Spinor representation in $Mat(4,\mathbb{H})$:

$$\g_{1}\leftrightarrow \scriptsize{\left[\begin{array}{cccc} 0 & 1 & 0 & 0\\ 1& 0 & 0 & 0\\ 0 & 0 & 0 & 1 \\ 0 & 0 & 1 & 0 \end{array}\right]}\;\;
\g_{2}\leftrightarrow \scriptsize{\left[\begin{array}{cccc} 0 & 0 & 1 & 0\\ 0 & 0 & 0 & -1\\ 1 & 0 & 0 & 0 \\ 0 & -1 & 0 & 0 \end{array}\right]}$$

$$\g_{3}\leftrightarrow \scriptsize{\left[\begin{array}{cccc} 0 & -1 & 0 & 0\\ 1& 0 & 0 & 0\\ 0 & 0 & 0 & -1 \\ 0 & 0 & 1 & 0 \end{array}\right]}\;\;
\g_{4}\leftrightarrow \scriptsize{\left[\begin{array}{cccc} 0 & 0 & -1 & 0\\ 0 & 0 & 0 & 1\\ 1 & 0 & 0 & 0 \\ 0 & -1 & 0 & 0 \end{array}\right]}$$

$$\g_{5}\leftrightarrow \scriptsize{\left[\begin{array}{cccc} i & 0 & 0 & 0\\ 0 & -i & 0 & 0\\ 0& 0 & -i & 0 \\ 0 & 0 &0 & i\end{array}\right]}\;\;
\g_{6}\leftrightarrow \scriptsize{\left[\begin{array}{cccc} j & 0 & 0 & 0\\ 0& -j & 0 & 0\\ 0 & 0 & -j & 0\\ 0 & 0 & 0 & j \end{array}\right]},$$

\item $Cl_{1,5}$ 
\subitem Generating set: $\{\g_{12},\g_{345}\}$ 
\subitem $S_{1,5}=span_{\R}\{1,\g_1,\g_3,\g_4,\g_6,\g_{13},\g_{14},\g_{16},\g_{34},\g_{36},\g_{46},\g_{134},\g_{136},\g_{146},\g_{346},\g_{1346}\}\otimes\{F\}\cong \R^{16}$.
\subitem  $S_{1,5}=span_{\mathbb{H}}\{F,\g_1F,\g_6F,\g_{16}F\}\cong\mathbb{H}^4$.
\subitem Identifying; $i\leftrightarrow\g_{3}F$, $j\leftrightarrow\g_{4}F$, $k\leftrightarrow\g_{34}F$

\subitem Spinor representation in $Mat(4,\mathbb{H})$:

$$\g_{1}\leftrightarrow \scriptsize{\left[\begin{array}{cccc} 0 & 1 & 0 & 0\\ 1& 0 & 0 & 0\\ 0 & 0 & 0 & 1 \\ 0 & 0 & 1 & 0 \end{array}\right]}\;\;
\g_{2}\leftrightarrow \scriptsize{\left[\begin{array}{cccc} 0 & -1 & 0 & 0\\ 1 & 0 & 0 & 0\\ 0 & 0 & 0 & -1 \\ 0 & 0 & 1 & 0 \end{array}\right]}$$

$$\g_{3}\leftrightarrow \scriptsize{\left[\begin{array}{cccc} i & 0 & 0 & 0\\ 0& -i & 0 & 0\\ 0 & 0 & -i & 0 \\ 0 & 0 & 0 & i \end{array}\right]}\;\;
\g_{4}\leftrightarrow \scriptsize{\left[\begin{array}{cccc} j & 0 & 0 & 0\\ 0 & -j & 0 & 0\\ 0 & 0 & -j & 0 \\ 0 & 0 & 0 & j \end{array}\right]}$$

$$\g_{5}\leftrightarrow \scriptsize{\left[\begin{array}{cccc} -k & 0 & 0 & 0\\ 0 & k & 0 & 0\\ 0& 0 & -k & 0 \\ 0 & 0 &0 & k\end{array}\right]}\;\;
\g_{6}\leftrightarrow \scriptsize{\left[\begin{array}{cccc} 0 & 0 & -1 & 0\\ 0& 0 & 0 & 1\\ 1 & 0 & 0 & 0\\ 0 & -1 & 0 & 0 \end{array}\right]},$$

\item $Cl_{7,1}$ 
\subitem Generating set: $\{\g_{18},\g_{2345},\g_{2367}\}$ 
\subitem $S_{7,1}=span_{\R}\{1,\g_1,\g_2,\g_3,\g_4,\g_6,\g_{12},\g_{13},\g_{14},\g_{16},\g_{23},\g_{24},\g_{26},\g_{34},\g_{36},\g_{46},\g_{123},\g_{124},\g_{126},\g_{134},$

$\g_{136},\g_{146},\g_{234},\g_{236},\g_{246},\g_{346},\g_{1234},\g_{1236},\g_{1246},\g_{1346},\g_{2346},\g_{12346}\}\otimes\{F\}\cong\R^{32}$.
\subitem  $S_{7,1}=span_{\mathbb{H}}\{F,\g_1F,\g_4F,\g_{6}F,\g_{14}F,\g_{16}F,\g_{46}F,\g_{146}F\}\cong\mathbb{H}^8$ .
\subitem Identifying; $i\leftrightarrow\g_{12}F$, $j\leftrightarrow\g_{23}F$, $k\leftrightarrow\g_{13}F$

\subitem Spinor representation in $Mat(8,\mathbb{H})$: 

$$\g_1 \leftrightarrow \scriptsize{\left[\begin{array}{cccccccc} 0 & 1 & 0 & 0 & 0 & 0 & 0 & 0\\ 1& 0 & 0 & 0 & 0 & 0 & 0 & 0\\ 0 & 0 & 0 & 0 & 1 & 0 & 0 & 0\\ 0 & 0 & 0 & 0 & 0 & 1 & 0 & 0\\ 0 & 0 & 1 & 0 & 0 & 0 & 0 & 0\\0 & 0 & 0 & 1 & 0 & 0 & 0 & 0\\ 0 & 0 & 0 & 0 & 0 & 0 & 0 & 1\\ 0 & 0 & 0 & 0 & 0 & 0 & 1 & 0 \end{array}\right]},\g_2 \leftrightarrow \scriptsize{\left[\begin{array}{cccccccc} 0 & -i & 0 & 0 & 0 & 0 & 0 & 0\\ i& 0 & 0 & 0 & 0 & 0 & 0 & 0\\ 0 & 0 & 0 & 0 & -i & 0 & 0 & 0\\ 0 & 0 & 0 & 0 & 0 & -i & 0 & 0\\ 0 & 0 & i & 0 & 0 & 0 & 0 & 0\\0 & 0 & 0 & i & 0 & 0 & 0 & 0\\ 0 & 0 & 0 & 0 & 0 & 0 & 0 & -i\\ 0 & 0 & 0 & 0 & 0 & 0 & i & 0 \end{array}\right]}$$

$$\g_3 \leftrightarrow \scriptsize{\left[\begin{array}{cccccccc} 0 & -k & 0 & 0 & 0 & 0 & 0 & 0\\ k& 0 & 0 & 0 & 0 & 0 & 0 & 0\\ 0 & 0 & 0 & 0 & -k & 0 & 0 & 0\\ 0 & 0 & 0 & 0 & 0 & -k & 0 & 0\\ 0 & 0 & k & 0 & 0 & 0 & 0 & 0\\0 & 0 & 0 & k & 0 & 0 & 0 & 0\\ 0 & 0 & 0 & 0 & 0 & 0 & 0 & -k\\ 0 & 0 & 0 & 0 & 0 & 0 & k & 0 \end{array}\right]},\g_4 \leftrightarrow \scriptsize{\left[\begin{array}{cccccccc} 0 & 0 & 1 & 0 & 0 & 0 & 0 & 0\\ 0 & 0 & 0 & 0 & -1 & 0 & 0 & 0\\ 1 & 0 & 0 & 0 & 0 & 0 & 0 & 0\\ 0 & 0 & 0 & 0 & 0 & 0 & 1 & 0\\ 0 & -1 & 0 & 0 & 0 & 0 & 0 & 0\\0 & 0 & 0 & 0 & 0 & 0 & 0 & -1\\ 0 & 0 & 0 & 1 & 0 & 0 & 0 & 0\\ 0 & 0 & 0 & 0 & 0 & -1 & 0 & 0 \end{array}\right]}$$

$$\g_5 \leftrightarrow \scriptsize{\left[\begin{array}{cccccccc} 0 & 0 & j & 0 & 0 & 0 & 0 & 0\\ 0 & 0 & 0 & 0 & -j & 0 & 0 & 0\\ -j & 0 & 0 & 0 & 0 & 0 & 0 & 0\\ 0 & 0 & 0 & 0 & 0 & 0 & j & 0\\ 0 & j & 0 & 0 & 0 & 0 & 0 & 0\\0 & 0 & 0 & 0 & 0 & 0 & 0 & -j\\ 0 & 0 & 0 & -j & 0 & 0 & 0 & 0\\ 0 & 0 & 0 & 0 & 0 & j & 0 & 0 \end{array}\right]},\g_6 \leftrightarrow \scriptsize{\left[\begin{array}{cccccccc} 0 & 0 & 0 & 1 & 0 & 0 & 0 & 0\\ 0 & 0 & 0 & 0 & 0 & -1 & 0 & 0\\ 0 & 0 & 0 & 0 & 0 & 0 & -1 & 0\\ 1 & 0 & 0 & 0 & 0 & 0 & 0 & 0\\ 0 & 0 & 0 & 0 & 0 & 0 & 0 & 1\\0 & -1 & 0 & 0 & 0 & 0 & 0 & 0\\ 0 & 0 & -1 & 0 & 0 & 0 & 0 & 0\\ 0 & 0 & 0 & 0 & 1 & 0 & 0 & 0 \end{array}\right]}$$

$$\g_7 \leftrightarrow \scriptsize{\left[\begin{array}{cccccccc} 0 & 0 & 0 & j & 0 & 0 & 0 & 0\\ 0 & 0 & 0 & 0 & 0 & -j & 0 & 0\\ 0 & 0 & 0 & 0 & 0 & 0 & -j & 0\\ -j & 0 & 0 & 0 & 0 & 0 & 0 & 0\\ 0 & 0 & 0 & 0 & 0 & 0 & 0 & j\\0 & j & 0 & 0 & 0 & 0 & 0 & 0\\ 0 & 0 & j & 0 & 0 & 0 & 0 & 0\\ 0 & 0 & 0 & 0 & -j & 0 & 0 & 0 \end{array}\right]},\g_8 \leftrightarrow \scriptsize{\left[\begin{array}{cccccccc} 0 & -1 & 0 & 0 & 0 & 0 & 0 & 0\\ 1& 0 & 0 & 0 & 0 & 0 & 0 & 0\\ 0 & 0 & 0 & 0 & -1 & 0 & 0 & 0\\ 0 & 0 & 0 & 0 & 0 & -1 & 0 & 0\\ 0 & 0 & 1 & 0 & 0 & 0 & 0 & 0\\0 & 0 & 0 & 1 & 0 & 0 & 0 & 0\\ 0 & 0 & 0 & 0 & 0 & 0 & 0 & -1\\ 0 & 0 & 0 & 0 & 0 & 0 & 1 & 0 \end{array}\right]}$$

\item $Cl_{6,2}$ 
\subitem Generating set: $\{\g_{17},\g_{28},\g_{3}\}$ 
\subitem $S_{6,2}=span_{\R}\{1,\g_1,\g_2,\g_4,\g_5,\g_6,\g_{12},\g_{14},\g_{15},\g_{16},\g_{24},\g_{25},\g_{26},\g_{45},\g_{46},\g_{56},\g_{124},\g_{125},\g_{126},\g_{145},$

$\g_{146},\g_{156},\g_{245},\g_{246},\g_{256},\g_{456},\g_{1245},\g_{1246},\g_{1256},\g_{1456},\g_{2456},\g_{12456}\}\otimes\{F\}\cong\R^{32}$.
\subitem  $S_{6,2}=span_{\mathbb{H}}\{F,\g_1F,\g_2F,\g_{6}F,\g_{12}F,\g_{16}F,\g_{26}F,\g_{126}F\}\cong \mathbb{H}^8$.
\subitem Identifying; $i\leftrightarrow\g_{24}F$, $j\leftrightarrow\g_{45}F$, $k\leftrightarrow\g_{25}F$.

\subitem Representation of generators in $M_8(\mathbb{H})$: 
$$\g_1 \leftrightarrow \scriptsize{\left[\begin{array}{cccccccc} 0 & 1 & 0 & 0 & 0 & 0 & 0 & 0\\ 1& 0 & 0 & 0 & 0 & 0 & 0 & 0\\ 0 & 0 & 0 & 0 & 1 & 0 & 0 & 0\\ 0 & 0 & 0 & 0 & 0 & 1 & 0 & 0\\ 0 & 0 & 1 & 0 & 0 & 0 & 0 & 0\\0 & 0 & 0 & 1 & 0 & 0 & 0 & 0\\ 0 & 0 & 0 & 0 & 0 & 0 & 0 & 1\\ 0 & 0 & 0 & 0 & 0 & 0 & 1 & 0 \end{array}\right]},\g_2 \leftrightarrow \scriptsize{\left[\begin{array}{cccccccc} 0 & 0 & 1 & 0 & 0 & 0 & 0 & 0\\ 0 & 0 & 0 & 0 & -1 & 0 & 0 & 0\\ 1 & 0 & 0 & 0 & 0 & 0 & 0 & 0\\ 0 & 0 & 0 & 0 & 0 & 0 & 1 & 0\\ 0 & -1 & 0 & 0 & 0 & 0 & 0 & 0\\0 & 0 & 0 & 0 & 0 & 0 & 0 & -1\\ 0 & 0 & 0 & 1 & 0 & 0 & 0 & 0\\ 0 & 0 & 0 & 0 & 0 & -1 & 0 & 0 \end{array}\right]}$$

$$\g_{3} \leftrightarrow\scriptsize{\left[\begin{array}{cccccccc} 1 & 0 & 0 & 0 & 0 & 0 & 0 & 0\\ 0& -1 & 0 & 0 & 0 & 0 & 0 & 0\\ 0 & 0 & -1 & 0 & 0 & 0 & 0 & 0\\ 0 & 0 & 0 & -1 & 0 & 0 & 0 & 0\\ 0 & 0 & 0 & 0 & 1 & 0 & 0 & 0\\0 & 0 & 0 & 0 & 0 & 1 & 0 & 0\\ 0 & 0 & 0 & 0 & 0 & 0 & 1 & 0\\ 0 & 0 & 0 & 0 & 0 & 0 & 0 & -1 \end{array}\right]},\g_4 \leftrightarrow \scriptsize{\left[\begin{array}{cccccccc} 0 & 0 & -i & 0 & 0 & 0 & 0 & 0\\ 0 & 0 & 0 & 0 & i & 0 & 0 & 0\\ i & 0 & 0 & 0 & 0 & 0 & 0 & 0\\ 0 & 0 & 0 & 0 & 0 & 0 & -i & 0\\ 0 & -i & 0 & 0 & 0 & 0 & 0 & 0\\0 & 0 & 0 & 0 & 0 & 0 & 0 & i\\ 0 & 0 & 0 & i & 0 & 0 & 0 & 0\\ 0 & 0 & 0 & 0 & 0 & -i & 0 & 0 \end{array}\right]}$$

$$\g_5 \leftrightarrow \scriptsize{\left[\begin{array}{cccccccc} 0 & 0 & -k & 0 & 0 & 0 & 0 & 0\\ 0 & 0 & 0 & 0 & k & 0 & 0 & 0\\ k & 0 & 0 & 0 & 0 & 0 & 0 & 0\\ 0 & 0 & 0 & 0 & 0 & 0 & -k & 0\\ 0 & -k & 0 & 0 & 0 & 0 & 0 & 0\\0 & 0 & 0 & 0 & 0 & 0 & 0 & k\\ 0 & 0 & 0 & k & 0 & 0 & 0 & 0\\ 0 & 0 & 0 & 0 & 0 & -k & 0 & 0 \end{array}\right]},\g_6 \leftrightarrow \scriptsize{\left[\begin{array}{cccccccc} 0 & 0 & 0 & 1 & 0 & 0 & 0 & 0\\ 0 & 0 & 0 & 0 & 0 & -1 & 0 & 0\\ 0 & 0 & 0 & 0 & 0 & 0 & -1 & 0\\ 1 & 0 & 0 & 0 & 0 & 0 & 0 & 0\\ 0 & 0 & 0 & 0 & 0 & 0 & 0 & 1\\0 & -1 & 0 & 0 & 0 & 0 & 0 & 0\\ 0 & 0 & -1 & 0 & 0 & 0 & 0 & 0\\ 0 & 0 & 0 & 0 & 1 & 0 & 0 & 0 \end{array}\right]}$$

$$\g_7 \leftrightarrow \scriptsize{\left[\begin{array}{cccccccc} 0 & -1 & 0 & 0 & 0 & 0 & 0 & 0\\ 1& 0 & 0 & 0 & 0 & 0 & 0 & 0\\ 0 & 0 & 0 & 0 & -1 & 0 & 0 & 0\\ 0 & 0 & 0 & 0 & 0 & -1 & 0 & 0\\ 0 & 0 & 1 & 0 & 0 & 0 & 0 & 0\\0 & 0 & 0 & 1 & 0 & 0 & 0 & 0\\ 0 & 0 & 0 & 0 & 0 & 0 & 0 & -1\\ 0 & 0 & 0 & 0 & 0 & 0 & 1 & 0 \end{array}\right]},\g_8 \leftrightarrow \scriptsize{\left[\begin{array}{cccccccc} 0 & 0 & -1 & 0 & 0 & 0 & 0 & 0\\ 0 & 0 & 0 & 0 & 1 & 0 & 0 & 0\\ 1 & 0 & 0 & 0 & 0 & 0 & 0 & 0\\ 0 & 0 & 0 & 0 & 0 & 0 & -1 & 0\\ 0 & -1 & 0 & 0 & 0 & 0 & 0 & 0\\0 & 0 & 0 & 0 & 0 & 0 & 0 & 1\\ 0 & 0 & 0 & 1 & 0 & 0 & 0 & 0\\ 0 & 0 & 0 & 0 & 0 & -1 & 0 & 0 \end{array}\right]}$$

\item $Cl_{3,5}$ 
\subitem Generating set: $\{\g_{14},\g_{25},\g_{36}\}$ 
\subitem $S_{3,5}=span_{\R}\{1,\g_1,\g_2,\g_3,\g_7,\g_8,\g_{12},\g_{13},\g_{17},\g_{18},\g_{23},\g_{27},\g_{28},\g_{37},\g_{38},\g_{78},\g_{123},\g_{127},\g_{128},$

$\g_{137},\g_{138},\g_{178},\g_{237},\g_{238},\g_{278},\g_{378},\g_{1237},\g_{1238},\g_{1278},\g_{1378},\g_{2378},\g_{12378}\}\otimes\{F\}\cong \R^{32}$.
\subitem  $S_{3,5}=span_{\mathbb{H}}\{F,\g_1F,\g_2F,\g_{3}F,\g_{12}F,\g_{13}F,\g_{23}F,\g_{123}F\}\cong\mathbb{H}^8$.
\subitem Identifying; $i\leftrightarrow\g_{7}F$, $j\leftrightarrow\g_{8}F$, $k\leftrightarrow\g_{78}F$

\subitem Spinor representation in $Mat(8,\mathbb{H})$ 

$$\g_1 \leftrightarrow \scriptsize{\left[\begin{array}{cccccccc} 0 & 1 & 0 & 0 & 0 & 0 & 0 & 0\\ 1& 0 & 0 & 0 & 0 & 0 & 0 & 0\\ 0 & 0 & 0 & 0 & 1 & 0 & 0 & 0\\ 0 & 0 & 0 & 0 & 0 & 1 & 0 & 0\\ 0 & 0 & 1 & 0 & 0 & 0 & 0 & 0\\0 & 0 & 0 & 1 & 0 & 0 & 0 & 0\\ 0 & 0 & 0 & 0 & 0 & 0 & 0 & 1\\ 0 & 0 & 0 & 0 & 0 & 0 & 1 & 0 \end{array}\right]},\g_2 \leftrightarrow \scriptsize{\left[\begin{array}{cccccccc} 0 & 0 & 1 & 0 & 0 & 0 & 0 & 0\\ 0 & 0 & 0 & 0 & -1 & 0 & 0 & 0\\ 1 & 0 & 0 & 0 & 0 & 0 & 0 & 0\\ 0 & 0 & 0 & 0 & 0 & 0 & 1 & 0\\ 0 & -1 & 0 & 0 & 0 & 0 & 0 & 0\\0 & 0 & 0 & 0 & 0 & 0 & 0 & -1\\ 0 & 0 & 0 & 1 & 0 & 0 & 0 & 0\\ 0 & 0 & 0 & 0 & 0 & -1 & 0 & 0 \end{array}\right]}$$
$$\g_3 \leftrightarrow \scriptsize{\left[\begin{array}{cccccccc} 0 & 0 & 0 & 1 & 0 & 0 & 0 & 0\\ 0 & 0 & 0 & 0 & 0 & -1 & 0 & 0\\ 0 & 0 & 0 & 0 & 0 & 0 & -1 & 0\\ 1 & 0 & 0 & 0 & 0 & 0 & 0 & 0\\ 0 & 0 & 0 & 0 & 0 & 0 & 0 & 1\\0 & -1 & 0 & 0 & 0 & 0 & 0 & 0\\ 0 & 0 & -1 & 0 & 0 & 0 & 0 & 0\\ 0 & 0 & 0 & 0 & 1 & 0 & 0 & 0 \end{array}\right]},\g_4 \leftrightarrow \scriptsize{\left[\begin{array}{cccccccc} 0 & -1 & 0 & 0 & 0 & 0 & 0 & 0\\ 1& 0 & 0 & 0 & 0 & 0 & 0 & 0\\ 0 & 0 & 0 & 0 & -1 & 0 & 0 & 0\\ 0 & 0 & 0 & 0 & 0 & -1 & 0 & 0\\ 0 & 0 & 1 & 0 & 0 & 0 & 0 & 0\\0 & 0 & 0 & 1 & 0 & 0 & 0 & 0\\ 0 & 0 & 0 & 0 & 0 & 0 & 0 & -1\\ 0 & 0 & 0 & 0 & 0 & 0 & 1 & 0 \end{array}\right]}$$

$$\g_5 \leftrightarrow \scriptsize{\left[\begin{array}{cccccccc} 0 & 0 & -1 & 0 & 0 & 0 & 0 & 0\\ 0 & 0 & 0 & 0 & 1 & 0 & 0 & 0\\ 1 & 0 & 0 & 0 & 0 & 0 & 0 & 0\\ 0 & 0 & 0 & 0 & 0 & 0 & -1 & 0\\ 0 & -1 & 0 & 0 & 0 & 0 & 0 & 0\\0 & 0 & 0 & 0 & 0 & 0 & 0 & 1\\ 0 & 0 & 0 & 1 & 0 & 0 & 0 & 0\\ 0 & 0 & 0 & 0 & 0 & -1 & 0 & 0 \end{array}\right]},\g_6 \leftrightarrow \scriptsize{\left[\begin{array}{cccccccc} 0 & 0 & 0 & -1 & 0 & 0 & 0 & 0\\ 0 & 0 & 0 & 0 & 0 & 1 & 0 & 0\\ 0 & 0 & 0 & 0 & 0 & 0 & 1 & 0\\ 1 & 0 & 0 & 0 & 0 & 0 & 0 & 0\\ 0 & 0 & 0 & 0 & 0 & 0 & 0 & -1\\0 & -1 & 0 & 0 & 0 & 0 & 0 & 0\\ 0 & 0 & -1 & 0 & 0 & 0 & 0 & 0\\ 0 & 0 & 0 & 0 & 1 & 0 & 0 & 0 \end{array}\right]}$$

$$\g_{7} \leftrightarrow\scriptsize{\left[\begin{array}{cccccccc} i & 0 & 0 & 0 & 0 & 0 & 0 & 0\\ 0& -i & 0 & 0 & 0 & 0 & 0 & 0\\ 0 & 0 & -i & 0 & 0 & 0 & 0 & 0\\ 0 & 0 & 0 & -i & 0 & 0 & 0 & 0\\ 0 & 0 & 0 & 0 & i & 0 & 0 & 0\\0 & 0 & 0 & 0 & 0 & i & 0 & 0\\ 0 & 0 & 0 & 0 & 0 & 0 & i & 0\\ 0 & 0 & 0 & 0 & 0 & 0 & 0 & -i \end{array}\right]},\g_{8} \leftrightarrow\scriptsize{\left[\begin{array}{cccccccc} j & 0 & 0 & 0 & 0 & 0 & 0 & 0\\ 0& -j & 0 & 0 & 0 & 0 & 0 & 0\\ 0 & 0 & -j & 0 & 0 & 0 & 0 & 0\\ 0 & 0 & 0 & -j & 0 & 0 & 0 & 0\\ 0 & 0 & 0 & 0 & j & 0 & 0 & 0\\0 & 0 & 0 & 0 & 0 & j & 0 & 0\\ 0 & 0 & 0 & 0 & 0 & 0 & j & 0\\ 0 & 0 & 0 & 0 & 0 & 0 & 0 & -j \end{array}\right]}$$

\item $Cl_{2,6}$ 
\subitem Generating set: $\{\g_{13},\g_{24},\g_{5678}\}$ 
\subitem $S_{2,6}=span_{\R}\{1,\g_1,\g_2,\g_5,\g_6,\g_7,\g_{12},\g_{15},\g_{16},\g_{17},\g_{25},\g_{26},\g_{27},\g_{56},\g_{57},\g_{67},\g_{125},\g_{126},\g_{127},$

$\g_{156},\g_{157},\g_{167},\g_{256},\g_{257},\g_{267},\g_{567},\g_{1236},\g_{1257},\g_{1267},\g_{1567},\g_{2567},\g_{12567}\}\otimes\{F\}\cong \R^{32}$.
\subitem  $S_{2,6}=span_{\mathbb{H}}\{F,\g_1F,\g_2F,\g_{7}F,\g_{12}F,\g_{17}F,\g_{27}F,\g_{127}F\}\cong\mathbb{H}^8$ .
\subitem Identifying; $i\leftrightarrow\g_{5}F$, $j\leftrightarrow\g_{6}F$, $k\leftrightarrow\g_{56}F$

\subitem Spinor representation in $Mat(8,\mathbb{H})$.
$$\g_1 \leftrightarrow \scriptsize{\left[\begin{array}{cccccccc} 0 & 1 & 0 & 0 & 0 & 0 & 0 & 0\\ 1& 0 & 0 & 0 & 0 & 0 & 0 & 0\\ 0 & 0 & 0 & 0 & 1 & 0 & 0 & 0\\ 0 & 0 & 0 & 0 & 0 & 1 & 0 & 0\\ 0 & 0 & 1 & 0 & 0 & 0 & 0 & 0\\0 & 0 & 0 & 1 & 0 & 0 & 0 & 0\\ 0 & 0 & 0 & 0 & 0 & 0 & 0 & 1\\ 0 & 0 & 0 & 0 & 0 & 0 & 1 & 0 \end{array}\right]},\g_2 \leftrightarrow \scriptsize{\left[\begin{array}{cccccccc} 0 & 0 & 1 & 0 & 0 & 0 & 0 & 0\\ 0 & 0 & 0 & 0 & -1 & 0 & 0 & 0\\ 1 & 0 & 0 & 0 & 0 & 0 & 0 & 0\\ 0 & 0 & 0 & 0 & 0 & 0 & 1 & 0\\ 0 & -1 & 0 & 0 & 0 & 0 & 0 & 0\\0 & 0 & 0 & 0 & 0 & 0 & 0 & -1\\ 0 & 0 & 0 & 1 & 0 & 0 & 0 & 0\\ 0 & 0 & 0 & 0 & 0 & -1 & 0 & 0 \end{array}\right]}$$

$$\g_3 \leftrightarrow \scriptsize{\left[\begin{array}{cccccccc} 0 & -1 & 0 & 0 & 0 & 0 & 0 & 0\\ 1& 0 & 0 & 0 & 0 & 0 & 0 & 0\\ 0 & 0 & 0 & 0 & -1 & 0 & 0 & 0\\ 0 & 0 & 0 & 0 & 0 & -1 & 0 & 0\\ 0 & 0 & 1 & 0 & 0 & 0 & 0 & 0\\0 & 0 & 0 & 1 & 0 & 0 & 0 & 0\\ 0 & 0 & 0 & 0 & 0 & 0 & 0 & -1\\ 0 & 0 & 0 & 0 & 0 & 0 & 1 & 0 \end{array}\right]},\g_4 \leftrightarrow \scriptsize{\left[\begin{array}{cccccccc} 0 & 0 & -1 & 0 & 0 & 0 & 0 & 0\\ 0 & 0 & 0 & 0 & 1 & 0 & 0 & 0\\ 1 & 0 & 0 & 0 & 0 & 0 & 0 & 0\\ 0 & 0 & 0 & 0 & 0 & 0 & -1 & 0\\ 0 & -1 & 0 & 0 & 0 & 0 & 0 & 0\\0 & 0 & 0 & 0 & 0 & 0 & 0 & 1\\ 0 & 0 & 0 & 1 & 0 & 0 & 0 & 0\\ 0 & 0 & 0 & 0 & 0 & -1 & 0 & 0 \end{array}\right]}$$
$$\g_{5} \leftrightarrow\scriptsize{\left[\begin{array}{cccccccc} i & 0 & 0 & 0 & 0 & 0 & 0 & 0\\ 0& -i & 0 & 0 & 0 & 0 & 0 & 0\\ 0 & 0 & -i & 0 & 0 & 0 & 0 & 0\\ 0 & 0 & 0 & -i & 0 & 0 & 0 & 0\\ 0 & 0 & 0 & 0 & i & 0 & 0 & 0\\0 & 0 & 0 & 0 & 0 & i & 0 & 0\\ 0 & 0 & 0 & 0 & 0 & 0 & i & 0\\ 0 & 0 & 0 & 0 & 0 & 0 & 0 & -i \end{array}\right]},\g_{6} \leftrightarrow\scriptsize{\left[\begin{array}{cccccccc} j & 0 & 0 & 0 & 0 & 0 & 0 & 0\\ 0& -j & 0 & 0 & 0 & 0 & 0 & 0\\ 0 & 0 & -j & 0 & 0 & 0 & 0 & 0\\ 0 & 0 & 0 & -j & 0 & 0 & 0 & 0\\ 0 & 0 & 0 & 0 & j & 0 & 0 & 0\\0 & 0 & 0 & 0 & 0 & j & 0 & 0\\ 0 & 0 & 0 & 0 & 0 & 0 & j & 0\\ 0 & 0 & 0 & 0 & 0 & 0 & 0 & -j \end{array}\right]}$$

$$\g_7 \leftrightarrow \scriptsize{\left[\begin{array}{cccccccc} 0 & 0 & 0 & -1 & 0 & 0 & 0 & 0\\ 0 & 0 & 0 & 0 & 0 & 1 & 0 & 0\\ 0 & 0 & 0 & 0 & 0 & 0 & 1 & 0\\ 1 & 0 & 0 & 0 & 0 & 0 & 0 & 0\\ 0 & 0 & 0 & 0 & 0 & 0 & 0 & -1\\0 & -1 & 0 & 0 & 0 & 0 & 0 & 0\\ 0 & 0 & -1 & 0 & 0 & 0 & 0 & 0\\ 0 & 0 & 0 & 0 & 1 & 0 & 0 & 0 \end{array}\right]},\g_8 \leftrightarrow \scriptsize{\left[\begin{array}{cccccccc} 0 & 0 & 0 & k & 0 & 0 & 0 & 0\\ 0 & 0 & 0 & 0 & 0 & -k & 0 & 0\\ 0 & 0 & 0 & 0 & 0 & 0 & -k & 0\\ k & 0 & 0 & 0 & 0 & 0 & 0 & 0\\ 0 & 0 & 0 & 0 & 0 & 0 & 0 & k\\0 & -k & 0 & 0 & 0 & 0 & 0 & 0\\ 0 & 0 & -k & 0 & 0 & 0 & 0 & 0\\ 0 & 0 & 0 & 0 & k & 0 & 0 & 0 \end{array}\right]}$$

\end{itemize}

\section{Octonion Algebra  $\mathbb{O}$}
\subsection{Cayle-Dickson doubling process}
we begin with an  important lemma found in [Har].
\begin{lem}
Suppose that $A$ is a sub-algebra of a normed algebra , and $\epsilon$ is a unit vector orthogonal to $A$ and $||\epsilon||=\pm 1$, then $A\epsilon$ is orthogonal to $A$ and multiplication in $A\oplus A\epsilon$ is defined as  $(x+y\epsilon)(z+w\epsilon)=(xz-\bar{w}y)+(wx+y\bar{z})\epsilon$ if $||\epsilon||=1$, and $(x+y\epsilon)(z+w\epsilon)=(xz+\bar{w}y)+(wx+y\bar{z})\epsilon$ if $||\epsilon||=-1$.

\end{lem}

The \textbf{Cayle-Dickson doubling process} is a direct consequence of this lemma [Har].
\begin{defn}
Suppose $A$ is normed algebra and $\epsilon$ a unit vector described in the previous lemma, then we define $\hat{A}=A\oplus A\epsilon$, where $||\epsilon||=1$ , and $\tilde{A}=A\oplus A\epsilon$, where $||\epsilon||=-1$ . $\hat{A}$ and $\tilde{A}$ are algebras with the multiplicative unit $1=(1_A,0)$.

\end{defn}
The constructions $\hat{A}$ and $\tilde{A}$ of a normed algebra $A$ are also normed algebras, where  $||x+y\epsilon||=||x||+||y||$ for $\hat{A}$, and $||x+y\epsilon||=||x||-||y||$ for $\tilde{A}$. The doubling process for $A$ also inherits the conjugation map $\overline{x+y\epsilon}=\bar{x}-y\epsilon$[Har]. The division algebras $\C$ and $\mathbb{H}$, are direct results of the doubling process for the algebras $\R$ and $\C$ respectively, that is ;

$$\C=\hat{\R}=\R\oplus \R i$$,

$$\mathbb{H}=\hat{\C}=\C\oplus \C j$$.

For $\R$ , the Lorentz numbers are   $\tilde{\R}$ and $\tilde{\C}=Mat(2,\R)$. When we apply the doubling process to the quaternions we get the division algebra known as the octonions denoted $\mathbb{O}$.
\subsection{Octonions}

In this section we focus on some elementary properties of the \textbf{Octonion algebra}. The octonion algebra is a  division algebra that can be thought of as a direct sum of two copies of the quaternions, 
$$\mathbb{O}=\mathbb{H}\oplus\mathbb{H}l$$, where $l$ is an imaginary unit that is orthogonal to $i,j,k$ and anti commutes with all the imaginary quaternions. Any number $x\in\mathbb{O}$ is the sum of two quaternions, that is $x=q_1+q_2l$, or the sum of $4$ complex numbers , $x=z_1+z_2i+z_3j+z_4k$, viewing the imaginary unit as $l$. We can also view ,as is more commonly done, an octonion as a sum of $8$ real numbers, that is $x=x_1+x_2i+x_3j+x_4k+x_5l+x_6il+x_7jl+x_8kl$. Thus the octonion algebras  contains $\R$ ,$\C$ and $\mathbb{H}$ as sub algebras, and  can be viewed as ;

$$\mathbb{O}=\R\oplus \R i\oplus \R j\oplus \R k\oplus \R l \oplus \R li\oplus \R lj\oplus \R lk$$,
$$\mathbb{O}=\C\oplus \C i\oplus \C j\oplus \C k$$,
$$\mathbb{O}=\mathbb{H}\oplus \mathbb{H}l$$.

$\mathbb{O}$ is commonly viewed as $\R\oplus \R^7$, since the imaginary octonions are of real dimension $7$,that is $im(\mathbb{O})= \R i\oplus \R j\oplus \R k\oplus \R l \oplus \R li\oplus \R lj\oplus \R lk\cong\R^7$, as a vector space. Multiplication in the octonions is obviously anti commutative but it is also anti associative, this can be seen if we view multiplication in terms of $$\mathbb{O}\times\mathbb{O}\rightarrow \mathbb{O}$$, where $x\in\mathbb{O}$ is an ordered pair of quaternions , that is $x=(q,r)$. So given $x,y\in\mathbb{O}$ such that $x=(q,r)$ and  $y=(s,t)$, the multiplication rule in the octonions is  $xy=(qs-\bar{t}r,tq+r\bar{s})$.Non associativity can easily been seen with this multiplication rule and is illustrated by the example   $(ij)l=kl$ and $i(jl)=-kl$ [DM]. We can view the standard basis described above as an ordered pair of quaternions as well, established in the following manner ;

$$1=(1,0),i=(i,0),j=(j,0),k=(k,0),l=(0,1),il=(0,i),jl=(0,j),kl=(0,k).$$  The cayle-table for octonion multiplication  is given as;

\begin{tabular}{c | c c c c c c c}
    $\mathbb{O}$ &  $i$ &  $j$ &  $k$ &  $l$ &  $il$ &  $jl$ &  $kl$ \\
    \cline{1-8}
    $i$ & $-1$ &  $k$ &  $-j$ &  $il$&  $-l$ &  $-kl$ &  $jl$\\
     $j$ &  $-k$ & $-1$ &  $i$ &  $jl$&  $kl$ &  $-l$ &  $jl$\\
     $k$ &  $j$ &  $-i$& $-1$ &  $kl$ &  $-jl$ & $il$ & $-l$ \\
     $l$ &  $-il$ &  $-jl$ & $-kl$ & $-1$&  $i$& $j$ &  $k$ \\
     $il$ &  $l$ &  $-kl$ &  $jl$ &  $-i$ & $-1$ & $-k$ & $j$ \\
      $jl$&  $kl$ &  $l$ &  $-il$ &  $-j$ &  $k$ &$-1$ & $-i$ \\
       $kl$ &  $-jl$&  $il$ &  $l$ &  $-k$&  $-j$ &  $i$&$-1$ \\
\end{tabular}

Where as matrices in $Mat(8,\R)$ we can view the left multiplication representations of the generators as ;

$$i \leftrightarrow \scriptsize{\left[\begin{array}{cccccccc} 0 & -1 & 0 & 0 & 0 & 0 & 0 & 0\\ 1& 0 & 0 & 0 & 0 & 0 & 0 & 0\\ 0 & 0 & 0 & -1 & 0 & 0 & 0 & 0\\ 0 & 0 & 1 & 0 & 0 & 0 & 0 & 0\\ 0 & 0 & 0 & 0 & 0 & -1 & 0 & 0\\0 & 0 & 0 & 0 & 1 & 0 & 0 & 0\\ 0 & 0 & 0 & 0 & 0 & 0 & 0 & 1\\ 0 & 0 & 0 & 0 & 0 & 0 & -1 & 0 \end{array}\right]},j \leftrightarrow \scriptsize{\left[\begin{array}{cccccccc} 0 & 0 & -1 & 0 & 0 & 0 & 0 & 0\\ 0 & 0 & 0 & 1 & 0 & 0 & 0 & 0\\ 1 & 0 & 0 & 0 & 0 & 0 & 0 & 0\\ 0 & -1 & 0 & 0 & 0 & 0 & 0 & 0\\ 0 & 0 & 0 & 0 & 0 & 0 & -1 & 0\\0 & 0 & 0 & 0 & 0 & 0 & 0 & -1\\ 0 & 0 & 0 & 0 & 1 & 0 & 0 & 0\\ 0 & 0 & 0 & 0 & 0 & 1 & 0 & 0 \end{array}\right]},$$

$$k \leftrightarrow \scriptsize{\left[\begin{array}{cccccccc} 0 & 0 & 0 & -1 & 0 & 0 & 0 & 0\\ 0& 0 & -1 & 0 & 0 & 0 & 0 & 0\\ 0 & 1 & 0 & 0 & 0 & 0 & 0 & 0\\ 1 & 0 & 0 & 0 & 0 & 0 & 0 & 0\\ 0 & 0 & 0 & 0 & 0 & 0 & 0 & -1\\0 & 0 & 0 & 0 & 0 & 0 & 1 & 0\\ 0 & 0 & 0 & 0 & 0 & -1 & 0 & 0\\ 0 & 0 & 0 & 0 & 1 & 0 & 0 & 0 \end{array}\right]},l\leftrightarrow \scriptsize{\left[\begin{array}{cccccccc} 0 & 0 & 0 & 0 & -1 & 0 & 0 & 0\\ 0& 0 & 0 & 0 & 0 & 1 & 0 & 0\\ 0 & 0 & 0 & 0 & 0 & 0 & 1 & 0\\ 0 & 0 & 0 & 0 & 0 & 0 & 0 & 1\\ 1 & 0 & 0 & 0 & 0 & 0 & 0 & 0\\0 & -1 & 0 & 0 & 0 & 0 & 0 & 0\\ 0 & 0 & -1 & 0 & 0 & 0 & 0 & 0\\ 0 & 0 & 0 & -1 & 0 & 0 & 0 & 0 \end{array}\right]}.$$

$$il \leftrightarrow \scriptsize{\left[\begin{array}{cccccccc} 0 & 0 & 0 & 0 & 0 & -1 & 0 & 0\\ 0 & 0 & 0 & 0 & -1 & 0 & 0 & 0\\ 0 & 0 & 0 & 0 & 0 & 0 & 0 & 1\\ 0 & 0 & 0 & 0 & 0 & 0 & -1 & 0\\ 0 & 1 & 0 & 0 & 0 & 0 & 0 & 0\\1 & 0 & 0 & 0 & 0 & 0 & 0 & 0\\ 0 & 0 & 0 & 1 & 0 & 0 & 0 & 0\\ 0 & 0 & -1 & 0 & 0 & 0 & 0 & 0 \end{array}\right]},jl\leftrightarrow \scriptsize{\left[\begin{array}{cccccccc} 0 & 0 & 0 & 0 & 0 & 0 & -1 & 0\\ 0& 0 & 0 & 0 & 0 & 0 & 0 & -1\\ 0 & 0 & 0 & 0 & -1 & 0 & 0 & 0\\ 0 & 0 & 0 & 0 & 0 & 1 & 0 & 0\\ 0 & 0 & 1 & 0 & 0 & 0 & 0 & 0\\0 & 0 & 0 & -1 & 0 & 0 & 0 & 0\\ 1 & 0 & 0 & 0 & 0 & 0 & 0 & 0\\ 0 & 1 & 0 & 0 & 0 & 0 & 0 & 0 \end{array}\right]} ,$$
$$kl \leftrightarrow \scriptsize{\left[\begin{array}{cccccccc} 0 & 0 & 0 & 0 & 0 & 0 & 0 & -1\\ 0& 0 & 0 & 0 & 0 & 0 & 1 & 0\\ 0 & 0 & 0 & 0 & 0 & -1 & 0 & 0\\ 0 & 0 & 0 & 0 & -1 & 0 & 0 & 0\\ 0 & 0 & 0 & 1 & 0 & 0 & 0 & 0\\0 & 0 & 1 & 0 & 0 & 0 & 0 & 0\\ 0 & -1 & 0 & 0 & 0 & 0 & 0 & 0\\ 1 & 0 & 0 & 0 & 0 & 0 & 0 & 0 \end{array}\right]},$$

 Now if we view $x\in\mathbb{O}$ as $x=a+u$, where $a\in\R$ and $u\in\R^7$(imaginary part of $\mathbb{O}$), the \textbf{conjugate} of $\mathbb{O}$ is just $\bar{x}=a-u$, with the inner product defined as  $<x,y>=Re(x\bar{y})$, or stated in vector form $<x,y>=ab-Re(uv)$, where $x=a+u$ and $y=b+v$. This inner product defines a norm $|x|=\sqrt{<x,x>}$, which allows us to define the inverse of any non-zero octonion.

 That inverse being 
$$x^{-1}=\dfrac{\bar{x}}{|x|^2}$$. The norm is obviously multiplicative and establishes $\mathbb{O}\cong\R^8$ as a normed space. The imaginary units of norm one are isomorphic to $S^6$, that is $\{u\in\mathbb{O}:|u|=1\}\cong S^6\subset\R^7$, since the imaginary units with norm one  satisfy the equation $x_2^2+x_3^2+x_4^2+x_5^2+x_6^2+x_7^2+x_8^2=1$. Also the octonions provide three equivalent ways to view rotations in $\R^8$

$$SO(8)=\{x\rightarrow pxp:p,x\in\mathbb{O}, |p|=1\}$$,
$$SO(8)=\{x\rightarrow px:p,x\in\mathbb{O}, |p|=1\}$$,
$$SO(8)=\{x\rightarrow xp:p,x\in\mathbb{O}, |p|=1\}$$
,
[DM]. 
\subsection{Split octonions}

The \textbf{split octonions} are the direct sum $$\tilde{\mathbb{O}}=\mathbb{H}\oplus \mathbb{H}\epsilon$$, where $\epsilon^2=1$. Multiplication in the split octonions has the following structure 

$$(q+r\epsilon)(s+t\epsilon)=(qs+\bar{t}r)+(tq+r\bar{s})\epsilon$$.  

 With the following Cayle table [DM].

\bigskip

\begin{tabular}{c | c c c c c c c}
    $\mathbb{\tilde{O}}$ &  $i$ &  $j$ &  $k$ &  $\epsilon$ &  $i\epsilon$ &  $j\epsilon$ &  $k\epsilon$ \\
    \cline{1-8}
    $i$ & $-1$ &  $k$ &  $-j$ &  $i\epsilon$&  $-\epsilon$ &  $-k\epsilon$ &  $j\epsilon$\\
     $j$ &  $-k$ & $-1$ &  $i$ &  $j\epsilon$&  $k\epsilon$ &  $-\epsilon$ &  $-i\epsilon$\\
     $k$ &  $j$ &  $-i$& $-1$ &  $k\epsilon$ &  $-j\epsilon$ & $i\epsilon$ & $-\epsilon$ \\
     $\epsilon$ &  $-i\epsilon$ &  $-j\epsilon$ & $-k\epsilon$ & $1$&  $-i$& $-j$ &  $-k$ \\
     $i\epsilon$ &  $\epsilon$ &  $-k\epsilon$ &  $j\epsilon$ &  $i$ & $1$ & $k$ & $-j$ \\
      $j\epsilon$&  $k\epsilon$ &  $\epsilon$ &  $-i\epsilon$ &  $j$ &  $-k$ &$1$ & $i$ \\
       $k\epsilon$ &  $-j\epsilon$&  $i\epsilon$ &  $\epsilon$ &  $k$&  $j$ &  $-i$&$1$ \\
\end{tabular}
\bigskip

Where we obtain the following left mutiplication representations in $Mat(8,\R)$ ;

$$i \leftrightarrow \scriptsize{\left[\begin{array}{cccccccc} 0 & -1 & 0 & 0 & 0 & 0 & 0 & 0\\ 1& 0 & 0 & 0 & 0 & 0 & 0 & 0\\ 0 & 0 & 0 & -1 & 0 & 0 & 0 & 0\\ 0 & 0 & 1 & 0 & 0 & 0 & 0 & 0\\ 0 & 0 & 0 & 0 & 0 & -1 & 0 & 0\\0 & 0 & 0 & 0 & 1 & 0 & 0 & 0\\ 0 & 0 & 0 & 0 & 0 & 0 & 0 & 1\\ 0 & 0 & 0 & 0 & 0 & 0 & -1 & 0 \end{array}\right]},j \leftrightarrow \scriptsize{\left[\begin{array}{cccccccc} 0 & 0 & -1 & 0 & 0 & 0 & 0 & 0\\ 0 & 0 & 0 & 1 & 0 & 0 & 0 & 0\\ 1 & 0 & 0 & 0 & 0 & 0 & 0 & 0\\ 0 & -1 & 0 & 0 & 0 & 0 & 0 & 0\\ 0 & 0 & 0 & 0 & 0 & 0 & -1 & 0\\0 & 0 & 0 & 0 & 0 & 0 & 0 & -1\\ 0 & 0 & 0 & 0 & 1 & 0 & 0 & 0\\ 0 & 0 & 0 & 0 & 0 & 1 & 0 & 0 \end{array}\right]},$$

$$k \leftrightarrow \scriptsize{\left[\begin{array}{cccccccc} 0 & 0 & 0 & -1 & 0 & 0 & 0 & 0\\ 0& 0 & -1 & 0 & 0 & 0 & 0 & 0\\ 0 & 1 & 0 & 0 & 0 & 0 & 0 & 0\\ 1 & 0 & 0 & 0 & 0 & 0 & 0 & 0\\ 0 & 0 & 0 & 0 & 0 & 0 & 0 & -1\\0 & 0 & 0 & 0 & 0 & 0 & 1 & 0\\ 0 & 0 & 0 & 0 & 0 & -1 & 0 & 0\\ 0 & 0 & 0 & 0 & 1 & 0 & 0 & 0 \end{array}\right]},\epsilon\leftrightarrow \scriptsize{\left[\begin{array}{cccccccc} 0 & 0 & 0 & 0 & 1 & 0 & 0 & 0\\ 0& 0 & 0 & 0 & 0 & -1 & 0 & 0\\ 0 & 0 & 0 & 0 & 0 & 0 & -1 & 0\\ 0 & 0 & 0 & 0 & 0 & 0 & 0 & -1\\ 1 & 0 & 0 & 0 & 0 & 0 & 0 & 0\\0 & -1 & 0 & 0 & 0 & 0 & 0 & 0\\ 0 & 0 & -1 & 0 & 0 & 0 & 0 & 0\\ 0 & 0 & 0 & -1 & 0 & 0 & 0 & 0 \end{array}\right]}.$$

$$i\epsilon \leftrightarrow \scriptsize{\left[\begin{array}{cccccccc} 0 & 0 & 0 & 0 & 0 & 1 & 0 & 0\\ 0 & 0 & 0 & 0 & 1 & 0 & 0 & 0\\ 0 & 0 & 0 & 0 & 0 & 0 & 0 & -1\\ 0 & 0 & 0 & 0 & 0 & 0 & 1 & 0\\ 0 & 1 & 0 & 0 & 0 & 0 & 0 & 0\\1 & 0 & 0 & 0 & 0 & 0 & 0 & 0\\ 0 & 0 & 0 & 1 & 0 & 0 & 0 & 0\\ 0 & 0 & -1 & 0 & 0 & 0 & 0 & 0 \end{array}\right]},j\epsilon\leftrightarrow \scriptsize{\left[\begin{array}{cccccccc} 0 & 0 & 0 & 0 & 0 & 0 & 1 & 0\\ 0& 0 & 0 & 0 & 0 & 0 & 0 & 1\\ 0 & 0 & 0 & 0 & 1 & 0 & 0 & 0\\ 0 & 0 & 0 & 0 & 0 & -1 & 0 & 0\\ 0 & 0 & 1 & 0 & 0 & 0 & 0 & 0\\0 & 0 & 0 & -1 & 0 & 0 & 0 & 0\\ 1 & 0 & 0 & 0 & 0 & 0 & 0 & 0\\ 0 & 1 & 0 & 0 & 0 & 0 & 0 & 0 \end{array}\right]} ,$$
$$k\epsilon \leftrightarrow \scriptsize{\left[\begin{array}{cccccccc} 0 & 0 & 0 & 0 & 0 & 0 & 0 & 1\\ 0& 0 & 0 & 0 & 0 & 0 & -1 & 0\\ 0 & 0 & 0 & 0 & 0 & 1 & 0 & 0\\ 0 & 0 & 0 & 0 & 1 & 0 & 0 & 0\\ 0 & 0 & 0 & 1 & 0 & 0 & 0 & 0\\0 & 0 & 1 & 0 & 0 & 0 & 0 & 0\\ 0 & -1 & 0 & 0 & 0 & 0 & 0 & 0\\ 1 & 0 & 0 & 0 & 0 & 0 & 0 & 0 \end{array}\right]},$$
\subsection{ Elementary facts about the exceptional lie group $G_2$}

If we let $A$ be a finite normed algebra , then its automorphism group is defined as ;
$$Aut(A)=\{g\in Gl(A):g(xy)=g(x)g(y),\forall x,y\in A\}$$. For  automorphisms of normed algebras we have $Aut(A)\subset O(Im(A))$[Har], where $Im(A)$ is the imaginary part of the normed algebra $A$. Two well known results are $Aut(\C)\cong\Z_2$ , and $Aut(\mathbb{H})=SO(3)$[Har]. The automorphisms in $\C$ are the identity and conjugation maps. For $\mathbb{O}$ , its automorphism group is called \textbf{the exceptional Lie group, $G_2$}, that is $G_2=Aut(\mathbb{O})$. $G_2$ is a compact Lie group that is a closed subgroup of $O(8)$.If $\phi\in G_2$ then we have the property that $<\phi(x),\phi(y)>=<x,y>$ for all $x,y\in \mathbb{O}$ [Y]. When we view the group action on $S^7$ by the $Spin(7)$, we can define $G_2$ as a the isotropy subgroup of $Spin(7)$. Given a representation $\rho:Cl_{0,7}\rightarrow End_{\R}(\R^8)$ , which can be viewed as an extension $\rho_v\in End_{\R}(\R^8)$, for $v\in Im(\mathbb{O})$, where $\rho_v(x)=vx$ for $x\in\R^8$. Then we can view  $G_2=\{g\in Spin(7):\rho_g(e)=e\}$, where $e$ is the unit element in $S^7$[LM].When we quotient $Spin(7)$ which is of dimension 21 by $G_2$ we have the well known diffeomorphism $Spin(7)/G_2\cong S^7$ , thus the dimension of $G_2$ is clearly 14.[LM]. 

      
       

Focusing on  $Spin(8)$, we have the isomorphism $Spin(8)/G_2\cong S^7\times S^7$ which can be established by the following commutative diagram [Po].

\[\begin{tikzcd}
G_2 \arrow{r} \arrow{d} & Spin(7)\arrow{r} \arrow{d} & S^7 \arrow{d}{=} \\
Spin(7)\arrow{d}\arrow{r} & Spin(8)\arrow{d} \arrow{r}& S^7\\ S^7 & S^7 
\end{tikzcd}
\]
.

The Lie algebra of $G_2$ , $\mathfrak{g}_2$, is a Lie sub-algebra of $\mathfrak{so}(7)$, and it is given by $\mathfrak{g}_2=\{Y\in End_{\R}(\mathbb{O}):Y(a)b+aY(b)=Y(ab), \forall a,b\in\mathbb{O}\}$ [Y]. When it comes to the concept of triality in $SO(8)$ we will use the following theorem in [Y]. 
\begin{thm}
For any $Z\in SO(8)$ there exist $X,Y\in SO(8)$ such that $X(a)Y(b)=Z(ab)$ , $a,b\in \mathbb{O}$, moreover $X$, $Y$ are uniquely determined up to sign by $Z$.

\end{thm}

In [Y] the group  $D=\{(X,Y,Z)\in SO(8)^3:X(a)Y(b)=\overline{Z(\overline{ab})}, a,b\in\mathbb{O}\}$ is introduced and  it is a  compact group that contains $Spin(7)$ as a  subgroup
where $D\cong Spin(8)$. Thus we can identify a triple in $D$ with elements in the $Spin(8)$, where the identification is useful in defining $G_2$. If we identify the automorphism $\nu:Spin(8)\rightarrow Spin(8)$ , with elements in $D$ such that $\nu(X,Y,Z)=(Y,Z,X)$, we can see that $\nu^3=id$. In terms of $Spin(8)$ we define  $G_2=\{\alpha\in Spin(8):\nu(\alpha)=\alpha\}$ [Y]. This section is nothing more than a mere introduction of ways to use the Spin groups, that can be calculated as restrictions of the representations of the appropriate Clifford algebra , in viewing the exceptional lie group $G_2$.  More details about triality can be found in [Po][Lo][LM]{Har][Me][Y]. 

\subsubsection{Special Thanks}
Id like to thank Dr.Anna Fino , Dr.Brian Sittinger , and Dr.Ivona Grzegorczyk for encouraging me to write this preprint.


\newpage

\begin{thebibliography}{EMG}
\addcontentsline{toc}{chapter}{Referencias}

\bibitem[Ab1]{Ab1} R. Ablamowicz, \emph{Lectures on Clifford (Geometric) Algebras and Applications}. Birkhäuser, Boston, 2004.

\bibitem[Ab2]{Ab2} R. Ablamowicz, \emph{Spinor representations of Clifford algebras: a symbolic approach}. Computer Physics Communications \textbf{115} (1998), 510-535.
\bibitem[Ab3]{Ab 3}  R. Ablamowicz, \emph{On the structure theorem of Clifford Algebras }, arXiv : 1610.0241bv1, Oct,7,2016.
\bibitem[BKU]{BKU} A. Bilge, S. Koçak, and S. Uğuz, \emph{Canonical Bases for Real Representations of Clifford Algebras}, Linear Algebra and Its Applications \textbf{2} (2008), 417-439.

\bibitem[Bu]{Bu} M. Burrow, \emph{Representation Theory of Finite Groups}. Academic Press, New York, 1971.
\bibitem[Bum]{Bum} D.Bump, \emph{Lie Groups}. Springer , GTM 225, New York ,2013.
\bibitem[Di]{Di} A. Dimakis, \emph{A New Representation for Spinors in Real Clifford Algebras}, Clifford Algebras and Their Applications in Mathematical Physics, Chisholm, J. S. R. Springer, Netherlands, 49-60, 1986.

\bibitem[DF]{DF} D. Dummit and R. Foote, \emph{Abstract Algebra, Third Edition}. Wiley, Hoboken, 2003.
\bibitem[DM]{DM} T.Dray and C.Manogue ,\emph{The Geometry of the Octonions}. World Scientific ,2015
\bibitem[Ga]{Ga} D. Garling, \emph{Clifford Algebras: An Introduction}. Cambridge University Press, Cambridge, 2011.
\bibitem[Har] R.Harvey \emph{Spinors and Calibrations}, Academic Press,1990   
\bibitem[Ha]{Ha} B. Hall,\emph{Lie Groups Lie Algebras and Representations},Springer ,2003, .

\bibitem[He]{He} D. Hestenes, \emph{Clifford algebra and the Interpretation
of Quantum Mechanics}, Conference Lecture published in Clifford Algebras and their Applications in Mathematical Physics. Reidel, Dordrecht/Boston (1986), 321-346.

\bibitem[HL]{HL} G. Hile and P. Lounesto, \emph{Matrix Representations of Clifford Algebras}, Linear Algebra and Its Applications 128 (1990), 51-63.
\bibitem[HYG] S.Huang, Y.Ying Qiao, G.Chun Wen \emph{Real and Complex Clifford analysis}, Springer , 2006.
\bibitem[La] J.M Landsberg \emph{Tensors :Geometry and Applications}, GSM, AMS, Volume 128, 2012
\bibitem[LM] H.B Lawson JR, ML Michelsohn, \emph{Spin Geometry}, Princeton University Press, 1989.
\bibitem[LW]{LW} P. Lounesto and G. Wene, \emph{Idempotent Structure of Clifford Algebras}, Acta Applicandae Mathematica,(1987), 9.
    
\bibitem[Lo]{Lo} P. Lounesto, \emph{Clifford Algebras and Spinors, Second Edition}.
Cambridge University Press, Cambridge, 2001.

\bibitem[Me]{Me} E. Meinrenken, \emph{Clifford Algebras and Lie Theory}. Springer, Berlin, 2013.

\bibitem[Po]{Po} Ian R. Porteous, \emph{Clifford Algebras and the Classical Groups}. Cambridge University Press, Cambridge, 1995.
\bibitem[R]{R} Wulf Rossmann , \emph{Lie Groups,An introduction through linear groups},Oxford graduate text in mathematics, Oxford University Press University Press, Oxford, 2002. 
\bibitem[Se]{Se} J. Serre, \emph{Linear representations of finite groups}. Springer-Verlag, New York, 1977.
 \bibitem[Sep]{Sep} M.Sepanski \emph{Compact Lie Groups} Springer, 2007.   
\bibitem[Sn]{Sn} J. Snygg, \emph{Clifford Algebra: A Computational Tool for Physicists}. Oxford University Press, New York, 1997.
 \bibitem[T]{T} M.Taylor \emph{The Octonions} : Lecture Notes.
 \bibitem[Ti]{Ti}, Y.Tian \emph{Matrix Representations of Octonions and their applications},arXiv:math:0003166v2, 1 April 2000.   
\bibitem[To]{To} I. Todorov, \emph{Clifford Algebras and Spinors}, Bulg. J. Phys, 38 (2011), 3-28.

\bibitem[Tr]{Tr} A. Trautman, \emph{Clifford Algebras and Their Representations}, Encyclopedia of Mathematical Physics (2006), 518-530.

\bibitem[Y]{Y} I. Yokota, \emph{Exceptional Lie Groups}, Arxiv:0902.0431v1  (2009), 
\end{thebibliography}
\end{document}